\documentclass[11pt,a4paper,leqno]{article}
\usepackage[latin1]{inputenc} 
\usepackage{amsmath}
\usepackage{amsfonts}
\usepackage{amssymb}
\usepackage{MnSymbol}
\usepackage{stmaryrd}

\usepackage{graphicx}
\usepackage{float}
\restylefloat{figure}

\newtheorem{satz}{Satz}[section]
\newtheorem{theorem}[satz]{Theorem}
\newtheorem{lemma}[satz]{Lemma}
\newtheorem{prop}[satz]{Proposition}
\newtheorem{cor}[satz]{Corollary}
\newtheorem{defin}[satz]{Definition}
\newtheorem{example}[satz]{Example}

\newcommand{\abs}[1]{\left|{#1}\right|}
\newcommand{\rund}[1]{\left(#1\right)}
\newcommand{\spitz}[1]{\left\langle{#1}\right\rangle}
\newcommand{\eckig}[1]{\left[{#1}\right]}
\newcommand{\schweif}[1]{\left\{#1\right\}}
\newcommand{\floor}[1]{\left\lfloor{#1}\right\rfloor}

\def\zz{\mathbb{Z}}
\def\cz{\mathbb{C}}
\def\nz{{\rm I\kern-.20em N}}

\def\D{\mathcal{D}}
\def\B{\mathcal{B}}
\def\M{\mathcal{M}}
\def\N{\mathcal{N}}
\def\O{\mathcal{O}}
\def\T{\mathcal{T}}
\def\S{\mathcal{S}}
\def\V{\mathcal{V}}
\def\X{\mathcal{X}}
\def\Y{\mathcal{Y}}
\def\W{\mathcal{W}}

\def\E{\mathfrak{E}}
\def\H{\mathfrak{H}}
\def\K{\mathfrak{K}}

\def\m{\mathfrak m}

\def\Res{{\rm Res}}
\def\Im{{\rm Im\ }}
\def\Mod{{\rm Mod}}
\def\hol{{\rm hol}}
\def\Aut{{\rm Aut}}
\def\Id{{\rm Id}}
\def\id{{\rm id}}
\def\ord{{\rm ord}}
\def\sdim{{\rm sdim\ }}
\def\rel{{\rm rel}}
\def\pr{{\rm pr}}
\def\Pr{{\rm Pr}}
\def\strong{{\rm strong}}
\def\versal{{\rm versal}}
\def\SUSY{{\rm SUSY}}

\def\eps{\varepsilon}

\begin{document}

\begin{center}

{\Huge \bf Versal Families of Compact Super

\vskip 0.4 true cm

Riemann Surfaces}

\vskip 1.0 true cm

Roland Knevel \\

Bar-Ilan University RAMAT GAN

Department of Mathematics 

ISRAEL
\end{center}

\vskip 2.0 true cm



{\bf Abstract:} We call every complex connected $(1,1)$-dimensional supermanifold a super Riemann surface and construct versal super families of compact ones, where the base spaces are allowed to be certain ringed spaces including all complex su\-per\-ma\-ni\-folds. Furthermore we choose maximal supersymmetric sub super families which turn out to be versal among all supersymmetric super families. In the cases where special divisors occur we prove the non-existence of versal super families and instead construct locally complete ones. For an accurate study of supersymmetric super families we introduce the duality functor, a covariant involution of the category of super families of compact super Riemann surfaces, and show that the supersymmetric super families are essentially the self-dual ones. As an application of the classification results it is shown that on a supersymmetric super fa\-mi\-ly of compact super Riemann surfaces locally in the base the supersymmetry is uniquely determined up to pullback by automorphisms with identity as body.

\section*{Introduction}

In 1857 Riemann was the first to observe that the complex structure on a compact connected $2$-dimensional smooth orientable manifold of genus $g \geq 2$ is determined by precisely $3 g - 3$ `independent' complex parameters, which he called moduli. This was the starting point of a lot of investigation and finally classification of complex structures on compact manifolds. It turned out that such geometric classification problems can be approached~- roughly speaking~- in two different ways:

\begin{itemize}
\item[(i)] construction of a coarse {\bf moduli space} $\Mod$~. Heuristically, its `points' are in $1$-$1$-correspondence with the isomorphy classes of the objects being classified. Precisely, $\Mod$ is equipped with a receipt how to assign to every family $\pi_M: M \twoheadrightarrow W$ of such objects a morphism $\xi_M: W \rightarrow \Mod \ $ in such a way that, given two families $M \twoheadrightarrow W$ and $N \twoheadrightarrow W$~, there exists a family isomorphism

\[
\begin{array}{ccc}
M & \simeq & N \\
\pi_M \twoheadsearrow & \circlearrowleft & \twoheadswarrow \pi_N \phantom{12} \\
& W &
\end{array} \,~,
\]

at least locally in $W$~, iff $\xi_M = \xi_N$~.

\item[(ii)] construction of a {\bf versal family} $M^\versal \twoheadrightarrow \T$~. A family $M^\versal$ is called {\bf versal} if first it is {\bf complete}, which means that for every family $N \twoheadrightarrow W$ locally in $W$ there exists a(n in general {\bf not} unique) morphism $\xi: W \rightarrow \T$ such that $N$ is isomorphic to the pullback of $M^\versal$ by $\xi$~:

\[
\begin{array}{ccc}
N & \mathop{\mathop{\longrightarrow}\limits^\sim}\limits^{(*)} & \xi^* M^\versal \\
\pi_N \twoheadsearrow & \circlearrowleft & \twoheadswarrow \pi_{\xi^* M^\versal}  \phantom{1} \\
& W &
\end{array} \,~,
\]

and second it is {\bf infinitesimally universal}, which means that for every point $u \in \T$ the tangent space of $\T$ at $u$ is in $1$-$1$-correspondence with the infinitesimal deformations of the fiber at $u$~. Versality often implies a so-called {\bf anchoring property}: $\xi$ is infact uniquely determined by prescribing the isomorphism $(*)$ on a certain sub\-fa\-mi\-ly of $N$~, and also {\bf local completeness}: for eve\-ry family $N \twoheadrightarrow W$ and isomorphism $\sigma: \pi_N^{- 1}\rund{w_0} \mathop{\rightarrow}\limits^\sim \pi_{M^\versal}^{- 1}\rund{v_0}$~, $v_0 \in \T$ and $w_0 \in W$~, locally in $W$ around $w_0$ there exists a morphism $\xi: W \rightarrow \T$ with $\xi\rund{w_0} = v_0$ and an isomorphism

\[
\begin{array}{ccc}
N & \mathop{\longrightarrow}\limits^\sim & \xi^* M^\versal \\
\pi_N \twoheadsearrow & \circlearrowleft & \twoheadswarrow \pi_{\xi^* M^\versal}  \phantom{1} \\
& W &
\end{array}
\]

extending $\sigma$~.
\end{itemize}

Are these two approaches disjoint? Not at all! In lucky cases the morphism $\xi$ in (ii) is completely determined by $N$~, and then $\T$ also serves as moduli space for the classification problem and the morphism $\xi$ in (ii) as $\xi_N$ in (i). In this case we say that $\T$ is a fine moduli space and call $M^\versal$ universal.

\begin{example} Given a compact Riemann surface $X$ of genus $g$~, its Jacobian variety is a $g$-dimensional torus, and for all $D \in \zz$ it serves as fine moduli space of all holomorphic line bundles $E \twoheadrightarrow X$ of degree $D$~.
\end{example}

In less lucky cases the morphism $\xi$ in (ii) is {\bf not} uniquely determined by $N$~. Then one tries to find the so-called modular group $\Gamma \sqsubset \Aut \ \T$ identifying all `points' with `isomorphic fiber' since then the quotient $\T / \Gamma$ can be taken as moduli space for the classification problem, and given $\xi$ by (ii), the composition $W \mathop{\longrightarrow}\limits^{\xi} \T \twoheadrightarrow \T / \Gamma$ serves as $\xi_N$ in (i).

\begin{example} The moduli space of all compact Riemann surfaces of genus $1$ is $\left.\T_1 \right/ P SL(2, \zz)$~, where $\T_1$ denotes the Teichmüller space for genus $1$~, so the upper half plane, while there exists {\bf no} versal family of compact super Riemann surfaces of genus $1$ on $\left.\T_1 \right/ P SL(2, \zz)$~. Such a versal family $M_1$ exists on $\T_1$ itself, and $P SL(2, \zz)$ is precisely its modular group.

Same procedure for genus $g \geq 2$ : There exists a versal family $M_g$ over the Teichmüller space $\T_g$~, and its modular group $\Gamma_g$ acts properly discontinuous on $\T_g$~. A versal family on the moduli space $\T_g \left/ \Gamma_g\right.$ does {\bf not} exist.
\end{example}

Generalizing the classical theory of ordinary compact Riemann surfaces, this article is devoted to classifying compact super Riemann surfaces. They are complex $(1, 1)$-dimensional supermanifolds (so locally equipped with one even and one odd complex coordinate) with an ordinary compact Riemann surface as body. Why study these objects? \\

The original motivation comes from physics, more precisely from the combination of ordinary string theory with the idea of supersymmetry. As in detail explained in \cite{BaraShvar}, \cite{CraneRab}, and \cite{Friedan}, it is widely believed that passing from ordinary to super string theory will kill the divergences when computing amplitudes in Polyakov's bosonic string theory and that compact super Riemann surfaces serve as a good model for world sheets of super strings. Therefore quantization of super string theory involves integration over all `possible complex structures' on a given world sheet, so in mathematical language integration over the super moduli space. However, in general the super moduli space as a quotient will {\bf not} have nice analytic properties, but the super Teichmüller space has - it will turn out to be a supermanifold if it exists. For mathematicians the world of complex supermanifolds and their deformation theory has become of constant interest: odd (anticommuting) variables are a way to increase - without changing the topology - the dimension of classical objects, and typical phenomena known and unknown from higher even dimensions occur. In section \ref{superfam} we will see a symmetry break which has no counterpart in the classical theory. \\

The obvious general advantages of a versal family $M^\versal$ are:

\begin{itemize}
\item it gives an overview over all possible objects of given type. As discribed before, it can become necessary to integrate over the `space' of all `isomorphy classes' of objects of given type.
\item often $M^\versal$ is described in an `optimal way', for example well adapted to symmetries and additional structures of $M^\versal$~. Now writing an arbitrary family $N$ as pullback of $M^\versal$ one has this optimal description at hand also for $N$~. As an application of this principle we show that every supersymmetric super family of compact super Riemann surfaces up to pullback by super family automorphisms with identity as body and locally in the base possesses only one supersymmetry, see corollaries \ref{SUSYunique} and \ref{SUSYunique1}.
\end{itemize}

The construction of general versal super families of compact super Riemann surfaces will be done by sheaf cohomological methods and seems to be completely new. Super moduli spaces of supersymmetric compact super Riemann surfaces have already been constructed by Crane and Rabin in 1988, \cite{CraneRab} section 3, and Natanzon in 1999, \cite{Natanzon} sections 10 - 12, by a completely different method, namely by classifying all lattices in a suitable super Lie group acting on the universal cover, so the super upper half plane for $g \geq 2$ resp. $\cz^{1|1}$ for $g = 1$~. Of course the results about the structure of the super moduli spaces of supersymmetric compact super Riemann surfaces as super orbifolds here, see sections \ref{superfamSUSY} and \ref{remain}, totally fit together with those by Crane, Rabin, and Natanzon. \\

{\bf Acknowledgement:} I would like to thank the Bar-Ilan University and in particular Andre Reznikov for the support of my research and scientific career. The research was financed by the center of excellence grant of Israel Science Foundation (grant No.1691/10).

\section{Basic concepts} \label{basics}

Throughout this article we have to deal with the category of ringed spaces. In particular complex supermanifolds will be described as ringed spaces with an ordinary complex manifold $M$ as base and unital $\zz_2$-graded structure sheaf locally $\simeq \O_M \otimes \bigwedge \cz^q$~. Given an ordinary complex manifold $M$ and $q \in \nz$ we denote by $M^{|q} := M \times \cz^{|q}$ the supermanifold whose structure sheaf globally equals $\O_M \otimes \bigwedge \cz^q$~. \\

Crucial will be the functor from the category of ringed spaces to the category of topological spaces by assigning to a ringed space $\V :=(V, \S)$ (we call $\S$ the structure sheaf of $\V$ and denote it also by $\O_\V$ sometimes) its underlying topological space $\V^\# := V$ and to a morphism $\Phi$ between ringed spaces its underlying continuous map $\Phi^\#$~. Adopting the terminology from supermanifolds we call it the body functor. Furthermore given a ringed space $\V = (V, \S)$ and an open subset $U \subset V$ we write $\V|_U := \rund{U, \S|_U}$~. In this article, without further mentionning, all ringed spaces are assumed to have the following properties:

\begin{itemize}
\item the body $V$ is a complex manifold,
\item the structure sheaf $\S = \S_0 \oplus \S_1$ is associative, unital and ( $\zz_2$-)graded commutative, and there exists a graded unital projection ${}^\#: \S \twoheadrightarrow \O_V$~, which so gives a natural embedding of the body $V$ into the whole ringed space,
\item every morphism $\Phi$ between two ringed spaces $\V = (V, \S)$ and $\W = (W, \T)$ preserves the $\zz_2$-grading of the structure sheaves. Since $\Phi^\# = \Phi|_V$ it is automatically a holomorphic map from $V$ to $W$~.
\end{itemize}

\begin{defin} Let $\M = \rund{M, \O_\M}$ and $\V = \rund{V, \S}$ be ringed spaces and $\m \lhd \S$ the nilradical.

\item[(i)] $\V$ is called an admissible base iff $\S / \m \simeq \O_V$ and all $\m^r \left/ \m^{r + 1}\right.$~, $r \in \nz \setminus \{0\}$~, are coherent as $\O_V$-modules. Locally of course $\m$ is nilpotent. If it is even globally nilpotent the smallest $k \in \nz \setminus \schweif{0}$ such that $\m^k = 0$ is called the nilpotency index of $\V$~. Obviously $k = 1$ iff $\V$ is an ordinary complex manifold. If $k \geq 2$ then $\V^\natural := \rund{V, \S \left/ \m^{k - 1}\right.}$ is an admissible base of nilpotency index $k - 1$~, and the natural projection $\S \twoheadrightarrow \S \left/ \m^{k - 1}\right.$ induces an embedding $\V^\natural \hookrightarrow \V$ with body $\id_V$~.


Of course, all complex supermanifolds are admissible bases.

\item[(ii)] A projection $\pi_\M: \M \twoheadrightarrow \V$~, $\V$ an admissible base, of ringed spaces is called a super family of complex supermanifolds of super dimension $(m, n)$ iff locally

\[
\begin{array}{ccc}
\M & \simeq & \V \times \cz^{m|n} \\
\pi_\M \twoheadsearrow & \circlearrowleft & \twoheadswarrow \Pr_\V \phantom{12} \\
& \V &
\end{array} \,~.
\]

Such a local isomorphism is called a local super chart of the super family $\M$~, and throughout the article $z$ and $\zeta$ will denote the even resp. odd local super coordinates in fiber direction $\cz^{m|n}$~.

Automatically $\pi_\M^\#: \M^\# \twoheadrightarrow V$ is a holomorphic family of complex manifolds of dimension $m$~. If $\V~=~V = \rund{V, \O_V}$ is an ordinary complex manifold then we call $\pi_\M: \M \twoheadrightarrow \V$ a(n ordinary) family of complex supermanifolds.

\item[(iii)] Given a second super family $\M \twoheadrightarrow \W$ of supermanifolds, a pair $(\Xi, \xi)$ of morphisms, where

\[
\begin{array}{ccc}
\phantom{123,} \M & \mathop{\longrightarrow}\limits^\Xi & \N \phantom{123} \\
\pi_\M \twoheaddownarrow & \circlearrowleft & \twoheaddownarrow \pi_\N \\
\phantom{123,} \V & \mathop{\longrightarrow}\limits_\xi & \W \phantom{123}
\end{array} \,~,
\]

or just $\Xi$ for short is called a super family morphism. $\Xi$ is called fiberwise biholomorphic iff $\Xi^\#: \M^\# \rightarrow \N^\#$ is fiberwise bijective, and in local super charts $\Xi$ can be expressed as $\rund{\xi, \Id_{\cz^{m|n}} }$~. $\Xi$ is called strong iff $\V = \W$ and $\xi = \Id_\W$~. The group of all super family automorphisms of $\M$ is denoted by $\Aut \ \M$~.

\item[(iv)] Every morphism $\xi: \V \rightarrow \W$~, $\W$ a second admissible base, induces a covariant functor $\xi^*$~, called the pullback under $\xi$~, from the category of super families of complex su\-per\-ma\-ni\-folds over $\W$ together with strong super family morphisms to the ones over $\V$ : Given a super family $\M \twoheadrightarrow \W$ of complex $(m, n)$-dimensional supermanifolds, $\xi^* \M$ is defined as the sub ringed space of $\V \times \M$ given by $\xi \circ \Pr_\V = \pi_\M \circ \Pr_\M$~. It is again a super family of complex $(m, n)$-dimensional supermanifolds via

\[
\pi_{\xi^* \M} := \left.\Pr_\V\right|_{\xi^* \M} : \xi^* \M \twoheadrightarrow \V \,~.
\]

The fiberwise biholomorphic super family morphism $\rund{\breve \xi_\M, \xi}: \xi^* \M \rightarrow \M$~, where \\
$\breve \xi_\M := \left.\Pr_\M\right|_{\xi^* \M} : \xi^* \M \rightarrow \M$~, is called the canonical lift of $\xi$ to $\xi^* \M$~. Given a strong super family morphism $\Phi: \M \rightarrow \N$~, $\N \twoheadrightarrow \W$ a second super family of complex supermanifolds, $\xi^* \Phi := \left.\rund{\Id_\V, \Phi}\right|_{\xi^* \M}: \xi^* \M \rightarrow \xi^* \N$ is the unique strong super family morphism making the diagram

\[
\begin{array}{ccc}
\phantom{123} \xi^* \M & \mathop{\longrightarrow}\limits^{\xi^* \Phi} & \xi^* \N \phantom{123} \\
\breve \xi_\M \downarrow & & \downarrow \breve \xi_\N \\
\phantom{123} \M & \mathop{\longrightarrow}\limits_\Phi & \N \phantom{123}
\end{array}
\]

commutative. The assignment $\xi \mapsto \xi^*$ itself is contravariant. In the special case where $\xi: \V \hookrightarrow \W$ is a natural embedding we call $\M|_{\W} := \xi^* \M$ the restriction of the super family $\M$ to $\W$~.

\item[(v)] Given a super family $\M \twoheadrightarrow \V$ of complex $(m, n)$-dimensional supermanifolds, \\
$sT^\rel \M := \ker d \pi_\M$ is called the relative tangent space of $\M$~. It is a super vector bundle on $\M$ of super rank $(m, n)$~. Furthermore the graded $\O_\V$-module morphism \\
$d \M : sT \V \rightarrow \rund{\pi_\M^\#}_{(1)} sT^\rel \M$ given by $sT \V \hookrightarrow \rund{\pi_\M^\#}_* \rund{\pi_\M^* sT \V} \rightarrow \rund{\pi_\M^\#}_{(1)} sT^\rel \M$~, where the second map is induced by the short exact sequence

\[
0 \hookrightarrow sT^\rel \M \hookrightarrow sT \M \mathop{\twoheadrightarrow}\limits^{d \pi_\M}  \pi_\M^* sT \V \twoheadrightarrow 0 \,~,
\]

is called the super differential of $\M$~. Take an atlas of $\M$ with transition morphisms $\rund{\Id_\V, \Phi_{i j}}$ and a section $\delta$ of $sT \V$~. Then $(d \M) \delta = \eckig{\rund{\eps_{i j}} }$ with the sections $\eps_{i j}$ of $sT^\rel \M$~, each on the overlap of the local super charts $i$ and $j$ and expressed in the local super chart $i$ by $\rund{d \Phi_{i j}}^{- 1} \Pr \rund{d \Phi_{i j}} \delta$~, where $\Pr: sT \rund{\V \times \cz^{m|n}} \twoheadrightarrow sT^\rel \rund{\V \times \cz^{m|n}}$ denotes the natural projection.

Given a second admissible base $\W$ and a morphism $\xi:~\W \rightarrow~\V$ we have a chain rule: $d \rund{\xi^* \M} = (\xi^* (d \M)) \ d \xi$ with the pullback $\xi^* (d \M): \xi^* sT^\rel \V \rightarrow \rund{\pi_{\xi^* \M}^\#}_{(1)} sT^\rel \rund{\xi^* \M}$ of $d \M$ under $\xi$ and $\breve \xi_\M$~.
\end{defin}

Finally, given super families $\M \twoheadrightarrow \V$ and $\N \twoheadrightarrow \W$ and a morphism $\xi: \V \rightarrow~\W$~, there is a $1$-$1$-correspondence between fiberwise biholomorphic super family morphisms $(\Xi, \xi): \M \rightarrow \N$ and strong super family isomorphisms $\Phi: \M \mathop{\rightarrow}\limits^\sim \xi^* \N$ given by

\begin{eqnarray*}
(\Xi, \xi) &\mapsto& \rund{\pi_\M, \Xi} \\
\rund{\breve \xi_\N \circ \Phi, \xi} &\mapsfrom& \Phi \,~.
\end{eqnarray*}

We will be mainly interested in the case of super families $\pi_\M:~\M~\twoheadrightarrow \V$ of super Riemann surfaces, where $(m, n) = (1, 1)$ and $\pi_\M^\#: \M^\# \twoheadrightarrow \V^\#$ is a holomorphic family of Riemann surfaces. Assume $\V = V$ is a complex manifold. Then there is an equivalence of categories


\[
\begin{array}{l}
\schweif{\text{pairs } (M, L) \, , \, M \twoheadrightarrow V \text{ a hol. family of Riemann surfaces and } X \twoheadrightarrow M \text{ a hol. line bundle}} \\
\phantom{1234} \leftrightarrow \schweif{\text{ord. families of super Riemann surfaces } \M \twoheadrightarrow V}
\end{array}
\]

\begin{eqnarray*}
(M, L) &\mapsto& \rund{M, \Gamma^\hol\rund{\bigwedge L^*}} \\
\rund{\M^\#, \left.\rund{sT \M}_1\right|_{\M^\#}} &\mapsfrom& \M \phantom{123456789012345} .
\end{eqnarray*}

So given a family $\M = \rund{M, \Gamma^\hol\rund{\bigwedge L^*}}$ of super Riemann surfaces, first we have an action of the bundle $\pr_V: V \times \cz^\times \twoheadrightarrow V$ on $\M$ via $\rund{\id_V, z | a \zeta}$ in every local super chart of $\M \times \cz^\times$~, $a$ denoting the coordinate on $\cz^\times$~, in other words by multiplying with $a$ in the line bundle $L$~, and so we obtain an exact sequence

\begin{eqnarray} \label{exact}
1 \hookrightarrow \O(V)^\times &\hookrightarrow& \Aut \ \M \mathop{\longrightarrow}\limits^{{}^\#} \Aut \ M \,~, \\
\notag f &\mapsto& \Phi_{\M, f}
\end{eqnarray}

where of course all $\Phi_{\M, f} \in \Aut_\strong \ \M$~. Second, $\rund{sT^\rel \M}_0$ splits as an exact sequence of $\O_M$-modules

\begin{eqnarray} \label{exacttangent}
0 \hookrightarrow \O_{M} &\hookrightarrow& \rund{sT^\rel \M}_0 \mathop{\twoheadrightarrow}\limits^{|_{\O_M}} T^\rel M \twoheadrightarrow 0 \,~, \\
\notag 1 &\mapsto& \zeta \partial_\zeta
\end{eqnarray}

and $\rund{sT^\rel \M}_1$ as a direct sum of $\O_M$-modules

\[
\rund{sT^\rel \M}_1 = T^\rel M \otimes L^* \oplus L \,~,
\]

the first summand generated by $\zeta \partial_z$ and the second one by $\partial_\zeta$~.

\begin{defin}
Let $\M \twoheadrightarrow \V$ be a super family of compact super Riemann surfaces, $V := \V^\#$ connected. Write $\M|_V = \rund{M, \Gamma^\hol\rund{\bigwedge L^*}}$ as above. Then the genus $g$ of $\rund{\pi_\M^\#}^{- 1}(v)$ and $D := \deg L|_{\rund{\pi_\M^\#}^{- 1}(v)}$ are independent of $v \in~V$~. The pair $(g, D)$ is called the type of $\M$~.
\end{defin}

\section{Duality} \label{duality}

In this section we construct an involutive covariant functor ${}^\vee$~, the duality, from the ca\-te\-go\-ry of super families of super Riemann surfaces together with fiberwise locally biholomorphic super family morphisms to itself such that the supersymmetric super families are essentially the self-dual ones, theorem \ref{charactSUSY}. Throughout this section let $\V = (V, \S)$ and $\W$ be admissible bases, $U, \Omega \subset \cz$ open and $\M \twoheadrightarrow \V$ and $\N \twoheadrightarrow \W$ super families of super Riemann surfaces. Let $\m \lhd \S$ denote the nilradical.

\begin{defin}[Supersymmetry]
A sub super vector bundle $\D~\subset~sT^\rel \M$ of super rank $(0, 1)$ is called a supersymmetry on $\M$ and the pair $(\M, \D)$ a supersymmetric super family iff the graded commutator

\[
[\phantom{1}, \phantom{1}]: \D \otimes \left.\D \rightarrow sT^\rel \M \right/ \D
\]

is non-degenerate, therefore surjective. If $\N$ is also equipped with a supersymmetry $\D'$~, a fiberwise locally biholomorphic super family morphism $\Phi: \M \rightarrow \N$ is called supersymmetric iff $\Phi^* \D' = (d \Phi) \D$~.
\end{defin}

The supersymmetry on $\V \times U^{|1}$ generated by $\zeta \partial_z + \partial_\zeta$ will be called standard, and we call an atlas of a supersymmetric super family $\M$ standard iff in the local super charts of this atlas the supersymmetry is given by the standard one. Easy calculations show that there is a $1$-$1$-correspondence

\[
\O\rund{\V \times U}_1 \times \O\rund{\V \times U}_0^\times \leftrightarrow \schweif{\text{supersymmetries on } \V \times U^{|1}}
\]

assigning to $(\lambda, F)$ the supersymmetry generated by $(\zeta + \lambda) \partial_z + F \partial_\zeta$~, and so every supersymmetric super family $\M$ admits a standard atlas. \\

\begin{defin} A super family morphism $\Phi: \V \times U^{|1} \rightarrow \W \times \Omega^{|1}$ is called split iff it is of the form $\Phi = (\varphi | A \eta)$ with a super family morphism $\varphi: \V \times U \rightarrow \W \times \Omega$ and $A \in \O(\V \times U)_0$~.
\end{defin}

We can already deduce:

\begin{cor} \label{SUSYunique0} Let $V$ be simply connected and $\M \twoheadrightarrow \V$ a supersymmetric super family of compact super Riemann surfaces admitting a standard atlas with only split transition morphisms. Then its supersymmetry is uniquely determined up to pullback by strong super family automorphisms with body $\id_{\M^\#}$~.
\end{cor}

In particular every ordinary family $\M \twoheadrightarrow V$ of compact super Riemann surfaces admits, locally in its base, at most one supersymmetry up to pullback by strong super family automorphisms with the identity as body, so by $\Phi_{\M, f}$~, $f \in \O(V)^\times$~. Later we will see (corollaries \ref{SUSYunique} and \ref{SUSYunique1}) that this is also true for every super family of compact super Riemann surfaces. \\

{\it Proof:} Let

\[
\Phi_{i j} = \rund{\varphi_{i j} \, \left| \, A_{i j} \zeta\right.}: \left.\rund{\V \times \cz^{|1}}\right|_{U_{i j}} \mathop{\rightarrow}\limits^\sim \left.\rund{\V \times \cz^{|1}}\right|_{U_{j i}} \,~,
\]

$U_i \subset V \times \cz$ open, $\varphi_{i j}: (\V \times \cz)|_{U_{i j}} \mathop{\rightarrow}\limits^\sim (\V \times \cz)|_{U_{j i}}$ strong and biholomorphic and $A_{i j}~\in~\O\rund{(\V \times \cz)|_{U_{i j}} }_0^\times$~, denote the split transition morphisms of the atlas of $\M$~, standard w.r.t. the supersymmetry $\D$ on $\M$~. Since all $\Phi_{i j}$ are supersymmetric w.r.t. the standard supersymmetry we have $A_{i j}^2 = \partial_z \varphi_{i j}$~. Let $\D'$ be a second supersymmetry on $\M$~, in the local super chart $i$ generated by $\rund{\zeta + \lambda_i} \partial_z + F_i \partial_\zeta$~, $\lambda_i \in \O\rund{(\V \times \cz)|_{U_i}}_1$~, $F_i \in \O\rund{(\V \times \cz)|_{U_i}}_0^\times$~. Then an easy calculation shows that $\lambda_j = \rund{A_{i j} \lambda_i} \circ \varphi_{i j}^{- 1}$ and $F_j = F_i \circ \varphi_{i j}^{- 1}$ on the overlaps of the local super charts $i$ and $j$~. Therefore and since all fibers of $\M^\#$ are compact all $F_i^\#$ glue together to a global function $f \in \O(V)^\times$~. So since moreover $V$ is simply connected we can choose squareroots $\sqrt{f} \in \O(V)^\times$ and $\sqrt{F_i} \in \O\rund{(\V \times \cz)|_{U_i}}_0^\times$ such that $\sqrt{F_i}^\# = \sqrt{f}$ and so $\sqrt{F_j} = \sqrt{F_i} \circ \varphi_{i j}$~. We see that

\[
\rund{\Id_\V, z \, \left| \, \sqrt{F_i} \zeta - \lambda_i\right.}: \left.\rund{\V \times \cz^{|1}}\right|_{U_i} \mathop{\rightarrow}\limits^\sim \left.\rund{\V \times \cz^{|1}}\right|_{U_i}
\]

glue together to a global strong supersymmetric isomorphism $\Xi:~(\M, \D) \mathop{\rightarrow}\limits^\sim \rund{\M, \D'}$ with $\Xi^\# = \id_{\M^\#}$~. $\Box$ \\

Let $\Phi: \V \times U^{|1} \rightarrow \W \times \Omega^{|1}$ be a fiberwise locally biholomorphic super family morphism. We want to assign a dual $\Phi^\vee: \V \times U^{|1} \rightarrow \W \times \Omega^{|1}$ to $\Phi$ in a covariant and involutive way, and we will do so in three steps. \\

{\it Step I}: Assume $\V = \W$ and $\Phi = (\varphi | A \zeta)$ split and strong, $\varphi: \V \times U \rightarrow \V \times \Omega$ strong and locally biholomorphic and $A \in \O(\V \times U)_0^\times$~. Then we define $\Phi^\vee := \rund{\varphi \left| \frac{\partial_z \varphi}{A} \zeta\right.}$~. \\

As an infinitesimal version of the later ${}^\vee$ we define the graded $\O_\V$-module involution

\begin{eqnarray}
{}^\vee : sT^\rel \rund{\V \times U^{|1}} &\rightarrow& sT^\rel \rund{\V \times U^{|1}} \,~, \notag \\
(a + b \zeta) \partial_z + (c + d \zeta) \partial_\zeta &\mapsto& (a + c \zeta) \partial_z + \rund{b + \rund{\partial_z a - d} \zeta} \partial_\zeta \,~, \notag
\end{eqnarray}

$a, b, c, d \in \O(\V \times U)$~.

\begin{lemma} \label{compatible} Let $\V = \W$ and $\Phi$ be strong.
\item[(i)] ${}^\vee$ is covariant on strong split super family morphisms and respects the graded commutator on $sT^\rel \rund{\V \times U^{|1}}$~.
\item[(ii)] If $\Phi$ is split and $\chi$ a section of $sT^\rel \rund{\V \times U^{|1}}$ then

\[
\rund{\rund{(d \Phi) \chi} \circ \Phi^{- 1}}^\vee = \rund{\rund{d \Phi^\vee} \chi^\vee} \circ \rund{\Phi^\vee}^{- 1}
\]

in $sT^\rel\rund{\V \times \Omega^{|1}}$~. This is the integral version of (i) second part.

\item[(iii)] If $\Phi$ is split and $\delta$ is a section of $sT \V \hookrightarrow sT\rund{\V \times \Omega^{|1}}$ then

\[
\rund{\Pr \rund{\rund{(d \Phi) \delta} \circ \Phi^{- 1}} }^\vee = \Pr \rund{\rund{\rund{d \Phi^\vee} \delta} \circ \rund{\Phi^\vee}^{- 1}} \,~,
\]

where $\Pr: sT \rund{\V \times U^{|1}} \twoheadrightarrow sT^\rel \rund{\V \times U^{|1}}$ denotes the canonical projection. In other words: $\, {}^\vee$ on vector fields is the infinitesimal analogon of $\, {}^\vee$ on strong split super family morphisms since $\Pr \rund{\rund{(d \Phi) \delta} \circ \Phi^{- 1}}$ is the derivation of $\Phi$ in $\delta$-direction!

\item[(iv)] $\Phi$ admits, locally in $\V$~, a representation

\[
\Phi = (\varphi | A \zeta) \circ \exp \eta \,~,
\]

$(\varphi | A \zeta): \V \times U^{|1} \rightarrow \V \times \Omega^{|1}$ a strong split locally biholomorphic super family morphism and $\eta \in H^0 \rund{\m \ sT^\rel\rund{\V \times U^{|1}} }_0$~.
\end{lemma}

In (iv) $\exp (t \eta)$ denotes the integral flow to $\eta \in H^0 \rund{\m \ sT^\rel\rund{\V \times U^{|1}} }_0$~, so the unique solution of the initial value problem $\Theta|_{\schweif{0} \times \V \times U^{|1}} = \Id_{\V \times U^{|1}}$ and $\Pr \rund{\rund{(d \Theta) \partial_t} \circ \Theta^{- 1}} = \eta$~, $\Theta$ a strong super family automorphism of $\cz \times \V \times U^{|1} \twoheadrightarrow \cz \times \V$~, where $t$ denotes the coordinate on $\cz$~. Since $\eta|_{V \times U^{|1}} = 0$ we have $\exp(t \eta)|_{V \times U^{|1}} = \Id_{V \times U^{|1}}$ and so in particular $\exp(t \eta)^\# = \id_{V \times U}$~. By induction on the nilpotency index of $\V$ one easily shows that $\eta$ is uniquely determined by $\exp \eta$~. \\

{\it Proof:} (i) straightforward calculation.

(ii) Let $z$ and $\zeta$ denote the super coordinates on $U^{|1}$ and $u$ and $\vartheta$ the ones on $\Omega^{|1}$~. Write $\Phi = (\varphi | A \zeta)$~. Easy calculation using

\[
\rund{\rund{d \Phi} \partial_z} \circ \Phi^{- 1} = \rund{\rund{\partial_z \varphi} \circ \varphi^{- 1}} \partial_u + \rund{\frac{\partial_z A}{A} \circ \varphi^{- 1}} \partial_\vartheta
\]

and $\rund{\rund{d \Phi} \partial_\zeta} \circ \Phi^{- 1} = \rund{A \circ \varphi^{- 1}} \partial_\vartheta$~.

(iii) Take super coordinates and $\Phi$ as in the proof of (ii). Again an easy calculation using $\eckig{\delta, \partial_z} = 0$ and

\[
\Pr \rund{\rund{(d \Phi) \delta} \circ \Phi^{- 1}} = \rund{(\delta \varphi) \circ \varphi^{- 1}} \partial_u + \rund{\frac{\delta A}{A} \circ \varphi^{- 1}} \vartheta \partial_\vartheta \,~.
\]

(iv) by induction on the nilpotency index $k$ of $\V$ after maybe shrinking $\V$~.

$k = 1$ : Then $\V$ is an ordinary complex domain, and $\Phi$ is already split.

$k \rightarrow k + 1$ : Assume $\V$ is of nilpotency index $k + 1$~. Then by induction hypothesis $\Phi|_{\V^\natural \times U^{|1}} = (\varphi | A \zeta) \circ \exp \eta$ with appropriate $(\varphi | A \zeta): \V^\natural \times U^{|1} \rightarrow \V^\natural \times \Omega^{|1}$ split and $\eta \in H^0 \rund{\m \ sT^\rel\rund{\V^\natural \times U^{|1}} }_0$~. Locally in $\V$ we choose extensions $\rund{\widetilde\varphi \left| \widetilde A \zeta\right.}: \V \times U^{| 1} \rightarrow \V \times \Omega^{| 1}$ of $(\varphi | A \zeta)$ and $\widehat\eta \in H^0\rund{\m \ sT^\rel\rund{\V \times U^{|1}} }_0$ of $\eta$~. Then

\[
\left.\rund{\Phi^{- 1} \circ \rund{\widetilde\varphi \left| \widetilde A \zeta\right.} \circ \exp \widehat\eta}\right|_{\V^\natural \times U^{|1}} = \Id_{\V^\natural \times U^{|1}} \,~,
\]

and so $\Phi^{- 1} \circ \rund{\widetilde\varphi \left| \widetilde A \zeta\right.} \circ \exp \widehat\eta = \Id_{\V \times U^{|1}} + \delta$ with

\[
\delta \in H^0 \rund{\m^k sT^\rel\rund{\V \times U^{|1}} }_0 = H^0 \rund{\m^k \boxtimes_{\O_V} sT^\rel\rund{V \times U^{|1}} }_0
\]

appropriately chosen. Therefore $\exp\rund{\widehat \eta - \delta} = \rund{\exp \widehat \eta} \circ \rund{\Id_{\V \times U^{|1}} - \delta}$~, and so

\[
\Phi^{- 1} \circ \rund{\widetilde\varphi \left| \widetilde A \zeta\right.} \circ \exp\rund{\widehat\eta - \delta} = \Id_{\V \times U^{|1}} \,~. \, \Box
\]

{\it Step II}: For $\V = \W$ and $\Phi = \exp \eta$~, $\eta \in H^0\rund{\m \ sT^\rel\rund{\V \times U^{|1}} }_0$~, we define $\Phi^\vee := \exp \eta^\vee$~. By lemma \ref{compatible} we already obtain a covariant involution of the category of super families of super Riemann surfaces over $\V$ together with strong locally biholomorphic super family morphisms, which commutes with pullbacks. \\

\begin{lemma} \label{compatible2}
Lemma \ref{compatible} (ii) and (iii) are true for general strong $\Phi$~.
\end{lemma}

{\it Proof:} It suffices to consider the case $\Phi = \exp \eta$~, $\eta \in H^0 \rund{\m \ sT^\rel \rund{\V \times U^{|1}} }_0$~. So let \\
$\Theta := \exp(t \eta)$~.

(ii) Let $\chi$ be a section of $sT^\rel \V$~. Then $\eta = \rund{\Pr (d \Theta) \partial_t} \circ \Theta^{- 1}$ yields

\[
\partial_t \rund{\rund{(d \Theta) \chi} \circ \Theta^{- 1}} = \eckig{\rund{(d \Theta) \chi} \circ \Theta^{- 1}, \eta} \,~.
\]

Using $\Theta^\vee = \exp\rund{t \eta^\vee}$ and lemma \ref{compatible} (i) we see that both sections $\rund{\rund{(d \Theta) \chi} \circ \Theta^{- 1}}^\vee$ and $\rund{(d \Theta)^\vee \chi^\vee} \circ \rund{\Theta^\vee}^{- 1}$ of $sT^\rel \rund{\cz \times \V \times U^{|1}}$ fulfill the same initial value problem $X|_{\schweif{0} \times \V \times U^{|1}} =~\eta^\vee$ and $\partial_t X = \eckig{X, \eta^\vee}$ and so are equal.

(iii) Same trick as in the proof of (ii) now using $\left.\Pr \rund{((d \Theta) \delta) \circ \Theta^{- 1}}\right|_{\schweif{0} \times \V \times U^{|1}} = 0$~, $\eckig{\delta, \partial_t} = 0$ and so

\[
\partial_t \rund{\Pr \rund{((d \Theta) \delta) \circ \Theta^{- 1}} } = \delta \eta + \eckig{\Pr \rund{((d \Theta) \delta) \circ \Theta^{- 1}}, \eta} \,~. \, \Box
\]

{\it Final step}: Let $\xi: \V \rightarrow \W$ be a morphism. Then we define

\[
\rund{\breve \xi_\N}^\vee := \breve \xi_{\N^\vee}: \rund{\xi^* \N}^\vee = \xi^* \N^\vee~\rightarrow~\N \,~.
\]

We are done since we can split any fiberwise locally biholomorphic super family morphism $(\Xi, \xi): \M \rightarrow \N$ into $(\Xi, \xi) = \rund{\breve \xi_\N, \xi} \circ \rund{\pi_\M, \Xi}$ with the associated strong locally biholomorphic super family morphism $\rund{\pi_\M, \Xi}: \M \rightarrow \xi^* \N$ to $\Xi$~. \\

Obviously $\rund{\M^\vee}^\# = \M^\#$ and $\rund{\Xi^\vee}^\# =~\Xi^\#$ for all super families $\M$ of super Riemann surfaces and fiberwise locally biholomorphic super family morphisms $\Xi$ between them. Now what does the functor ${}^\vee$ have to do with supersymmetries? Here the answer:

\begin{lemma} \label{respectSUSY} Let $\D$ and $\D'$ denote the standard supersymmetries on $\V \times \cz^{|1}$ resp. $\W \times \cz^{|1}$~.
\item[(i)] Let $\eta$ be a section of $sT^\rel \rund{\V \times U^{|1}}$~. Then $\eckig{\D, \eta} \subset \D$ iff $\eta^\vee = \eta$~.
\item[(ii)] Let $\Xi: \V \times U^{|1} \rightarrow \W \times \Omega^{|1}$ be a fiberwise locally biholomorphic super family morphism. Then $\Xi$ is supersymmetric w.r.t. $\D$ and $\D'$ iff $\Xi^\vee = \Xi$~.
\end{lemma}

{\it Proof:} (i) easy exercise as is (ii) for split $\Xi$~.

(ii) general: It remains to prove the statement for $\V = \W$ and $\Xi = \exp \eta$~, \\
$\eta~\in~H^0 \rund{\m \ sT^\rel \rund{\V \times U^{|1}} }_0$~.

`$\Leftarrow$': Let $\chi$ be a section of $\D$ and $\Theta := \exp(t \eta)$~. Again

\[
\partial_t \rund{( (d \Theta) \chi) \circ \Theta^{- 1}} = \eckig{( (d \Theta) \chi) \circ \Theta^{- 1}, \eta} \,~.
\]

Since $\Xi^\vee = \Xi$ also $\eta^\vee = \eta$~. Therefore by (i) $\eckig{\diamondsuit, \eta}$ factors through the natural projection $\overline{\phantom{1}} : sT^\rel \rund{\cz \times \V \times U^{|1}} \twoheadrightarrow \left.sT^\rel \rund{\cz \times \V \times U^{|1}}\right/ \D$~, and so both sections $\overline{( (d \Theta) \chi) \circ \Theta^{- 1}}$ and $0$ of $\left.sT^\rel \rund{\cz \times \V \times U^{|1}}\right/ \D$ are solutions of the intial value problem $X|_{\schweif{0} \times \V \times U^{|1}} = 0$ and $\partial_t X = \eckig{X, \eta}$ and therefore are equal. \\

`$\Rightarrow$': by induction on the nilpotency index $k$ of $\V$ after maybe shrinking $\V$~.

$k = 1$ : Then $\Xi = \Id_{\V \times \cz^{1|1}}$~.

$k \rightarrow k + 1$ : Assume $\V$ is of nilpotency index $k + 1$~. Then by induction hypothesis $\Xi|_{\V^\natural \times U^{|1}}$ is already invariant under ${}^\vee$~. Define $\Psi~:=~\Xi^{- 1} \circ \Xi^\vee$~. Then $\Psi|_{\V^\natural \times U^{|1}} = \Id_{\V^\natural \times U^{|1}}$ and $\Psi^\vee = \Psi^{- 1}$~. Therefore $\Psi~=~\Id_{\V \times U^{|1}} + \delta$ with an appropriate

\[
\delta \in H^0 \rund{\m^k sT^\rel \M}_0 = H^0 \rund{\m^k \boxtimes_{\O_V} sT^\rel \M|_V}_0
\]

having $\delta^\vee = - \delta$~. On the other hand since $\D^\vee = \D$~, also $\chi \mapsto ((d \Psi) \chi) \circ \Psi^{- 1} = \chi + \eckig{\chi, \delta}$ respects $\D$~. Therefore by (i) $\delta^\vee = \delta$ and so $\delta = 0$~. $\Box$ \\

We see that $\M$ is supersymmetric iff it admits an atlas with self-dual transition morphisms. Moreover:

\begin{theorem} \label{charactSUSY} $\M$ is supersymmetric iff there exists a strong super family isomorphism $\Xi: \M \mathop{\rightarrow}\limits^\sim \M^\vee$ with $\Xi^\# = \id_M$ and $\Xi^\vee = \Xi^{- 1}$~.
\end{theorem}


{\it Proof:} `$\Rightarrow$': Choose local super charts with self-dual transition morphisms and define $\Xi$ to be the identity in these local super charts.

`$\Leftarrow$': Choose an atlas on $\M$ with local super charts $\left.\rund{\V \times \cz^{|1}}\right|_{U_i}$~, $U_i~\subset~V \times~\cz$ open, and transition morphisms $\Phi_{i j}: \left.\rund{\V \times \cz^{|1}}\right|_{U_{i j}} \rightarrow \left.\rund{\V \times \cz^{|1}}\right|_{U_{j i}}$ and write $\Xi$ in these local super charts as $\Xi_i: \left.\rund{\V \times \cz^{|1}}\right|_{U_i} \rightarrow \left.\rund{\V \times \cz^{|1}}\right|_{U_i}$~, $\Xi_i^\vee =~\Xi_i^{- 1}$~. We show that locally in $\M$ there exist $\Theta_i \in \Aut_\strong \ \left.\rund{\V \times \cz^{|1}}\right|_{U_i}$ such that $\Theta_i^\# = \id_{U_i}$ and $\Xi_i = \Theta_i^{- 1} \circ \Theta_i^\vee$ by induction on the nilpotency index $k$ of $\V$ :

\begin{quote}
$k = 1$ : Then $\V$ is an ordinary complex manifold, and $\Xi_i = \rund{\left.\id_{U_i} \right| A \zeta}$ with some 
$A \in \O\rund{U_i}^\times$~. Locally we can define $\Theta_i := \rund{\id_{U_i} \left| \frac{1}{\sqrt{A}} \zeta\right.}$ with a sqareroot $\sqrt{A} \in \O\rund{U_i}^\times$~.

$k \rightarrow k + 1$ : Assume $\V$ is of nilpotency index $k + 1$~. Then by induction hy\-po\-the\-sis already $\left.\Xi_i\right|_{\V^\natural \times U^{|1}} = \Theta_i^{- 1} \circ \Theta_i^\vee$ with some $\Theta_i \in \Aut_\strong \ \left.\rund{\V^\natural \times \cz^{|1}}\right|_{U_i}$ having $\Theta_i^\#~=~\id_{U_i}$~. We choose a local extension $\widehat\Theta \in \Aut_\strong \ \left.\rund{\V \times \cz^{|1}}\right|_{U_i}$ of $\Theta_i$~. We see that $\left.\rund{\widehat\Theta \circ \Xi_i \circ \rund{\widehat\Theta^\vee}^{- 1}}\right|_{\V^\natural \times U^{|1}}~=~\Id_{\V^\natural \times U^{|1}}$ and so $\widehat\Theta \circ \Xi_i \circ \rund{\widehat\Theta^\vee}^{- 1} = \Id_{\V \times U^{|1}} + \delta$ with some

\[
\delta \in H^0 \rund{U_i, \m^k sT^\rel \rund{\V \times \cz^{|1}} }_0 = H^0\rund{U_i, \m^k \boxtimes_{\O_V} sT^\rel \rund{V \times \cz^{|1}} }_0
\]

having $\delta^\vee = - \delta$ since $\rund{\widehat\Theta \circ \Xi_i \circ \rund{\widehat\Theta^\vee}^{- 1}}^\vee = \rund{\widehat\Theta \circ \Xi_i \circ \rund{\widehat\Theta^\vee}^{- 1}}^{- 1}$~. Therefore \\
$\Xi_i = {\widetilde \Theta_i}^{- 1} \circ {\widetilde \Theta_i}^\vee$ with

\[
\widetilde\Theta_i := \rund{\Id_{\left.\rund{\V \times \cz^{|1}}\right|_{U_i}} - \frac{1}{2} \delta} \circ \widehat\Theta \in \Aut_\strong \ \left.\rund{\V \times \cz^{|1}}\right|_{U_i} \,~.
\]

\end{quote}

Overall we obtain, after refining the atlas,

\[
\begin{array}{ccccc}
\phantom{123} \left.\rund{\V \times \cz^{|1}}\right|_{U_{i j}} & \mathop{\longleftarrow}\limits^{\Theta_i} & \left.\rund{\V \times \cz^{|1}}\right|_{U_{i j}} & \mathop{\longrightarrow}\limits^{\Theta_i^\vee} & \left.\rund{\V \times \cz^{|1}}\right|_{U_{i j}} \phantom{12} \\
\Phi_{ij} \downarrow && \circlearrowleft && \downarrow \Phi_{ij}^\vee \\
\phantom{123} \left.\rund{\V \times \cz^{|1}}\right|_{U_{j i}} & \mathop{\longleftarrow}\limits_{\Theta_j} & \left.\rund{\V \times \cz^{|1}}\right|_{U_{j i}} & \mathop{\longrightarrow}\limits_{\Theta_j^\vee} & \left.\rund{\V \times \cz^{|1}}\right|_{U_{j i}} \phantom{12}
\end{array} \,~,
\]

and so we can pass to a new atlas of $\M$ with the self-dual transition morphisms $\Theta_j^{- 1} \circ \Phi_{i j} \circ \Theta_i$~. $\Box$ \\

In the special case where $\V = V$ is an ordinary complex manifold we write \\
$\M = \rund{M, \Gamma^\hol\rund{\bigwedge L^*}}$~, $L \twoheadrightarrow M$ a holomorphic line bundle, and observe that

\begin{itemize}
\item $\Xi^\# = \id_M$ automatically implies $\Xi^\vee = \Xi^{- 1}$~,
\item $\M^\vee = \rund{M, \Gamma^\hol\rund{\bigwedge \rund{T^\rel M \otimes L^*}^*} }$~, and so $\M$ is supersymmetric iff $L^{\otimes 2} \simeq T^\rel M$~, which means that $L^*$ is a relative theta characteristic.
\end{itemize}

In the general case we conclude that

\begin{itemize}
\item if $\M$ is of type $(g, D)$ then $\M^\vee$ is of type $(g, 2 (1 - g) - D)$~, and so if $\M$ is supersymmetric then it must be of type $(g, 1 - g)$ for some $g \in \nz$~,
\item by lemma \ref{compatible2}, ${}^\vee$ taken in the local super charts of a standard atlas glue together to an $\S$-linear involution of $sT^\rel \M$~. Let $sT^\rel \M = sT^{\rel, +} \M \oplus sT^{\rel, -} \M$ be the graded decompomposition into the $(+ 1)$- and $(- 1)$-eigensheaf of ${}^\vee$~. Then by lemma \ref{respectSUSY} (i) $sT^{\rel, +} \M$ is the sheaf of all sections $\chi$ of $sT^\rel \M$ with $\eckig{\D, \chi} \subset \D$~.
\item given a single supersymmetric Riemann surface $\X$~, the space of infinitesimal deformations of $\X$ preserving the supersymmetry is precisely $H^1\rund{sT^+ \X}$~.
\end{itemize}

\begin{lemma} \label{SUSYprop}
\item[(i)] Let $\V = V$ be an ordinary complex manifold and $\M$ be supersymmetric. Write $\M~=~\rund{M, \Gamma^\hol\rund{\bigwedge L^*}}$ with a holomorphic line bundle $L \twoheadrightarrow M$ having $L^{\otimes 2} = T^\rel M$~. Then we have isomorphisms

\[
T^\rel M \mathop{\rightarrow}\limits^\sim \rund{sT^{\rel, +} \M}_0 \,~, \, f \partial_z \mapsto f \partial_z + \frac{1}{2} \rund{\partial_z f} \zeta \partial_\zeta
\]

in any standard standard atlas of $\M$~,

\[
\O_M \mathop{\rightarrow}\limits^\sim \rund{sT^{\rel, -} \M}_0 \,~, \, 1 \mapsto \zeta \partial_\zeta
\]

in any local super chart of $\M$~,

\[
\rund{\id_L, \id_L} : L \mathop{\rightarrow}\limits^\sim \rund{sT^{\rel, +} \M}_1
\]

and 

\[
\rund{\id_L, - \id_L} : L \mathop{\rightarrow}\limits^\sim \rund{sT^{\rel, -} \M}_1 \,~,
\]

the first of $\O_V$-modules, the other ones of $\O_M$-modules.

\item[(ii)] If $\M$ is supersymmetric the image of $d \M$ lies in $\rund{\pi_\M^\#}_{(1)} sT^{\rel, +} \M$~.
\end{lemma}

{\it Proof:} (i) Of course the assignment $f \partial_z \mapsto f \partial_z + \frac{1}{2} \rund{\partial_z f} \zeta \partial_\zeta$ is a local cross section to the projection of (\ref{exacttangent}) in section \ref{basics}. It is uniquely determined by the fact that its image lies in $\rund{sT^{\rel, +} \M}_0$ and is therefore globally defined. The rest is easy.

(ii) follows directly from lemma \ref{compatible2}. $\Box$

\section{Ordinary families} \label{ordfam}

Let $(g, D) \in \nz \times \zz$~. In this section we construct a versal family $\M_{g, D}$ amoung all ordinary families of compact super Riemann surfaces of type $(g, D)$~. In principle we have to classify pairs $(M, L)$~, where $M \twoheadrightarrow V$ is a family of compact Riemann surfaces of genus $g$ and $L$ is a holomorphic line bundle on $M$ of degree $D$~. \\

By classical Teichmüller theory, see \cite{Ares}, \cite{Duma}, \cite{Earle}, \cite{Hurwitz}, \cite{Kodaira}, and \cite{Natanzon}, we already know that there exists a versal family $\pi_{M_g}: M_g \twoheadrightarrow \T_g$ of compact Riemann surfaces $X_u := \pi_{M_g}^{- 1}(u)$~, $u \in \T_g$~, of genus $g$ on the Teichmüller space $\T_g$~, which is a contractible bounded domain of holomorphy, for genus $g$~:

\begin{itemize}
\item {\bf Completeness:} Let $N \twoheadrightarrow W$ be a holomorphic family of compact Riemann surfaces $Y_w$ of genus $g$~, $W$ simply connected if $g \geq 2$~, $W$ a contractible Stein manifold if $g \leq 1$~. Then there exists a holomorphic map $\varphi: W \rightarrow \T_g$ such that $N \simeq \varphi^* M_g$~. More precisely we have the following

{\bf Anchoring property (including local completeness):} For every isomorphism $\sigma: Y_{w_0} \mathop{\rightarrow}\limits^\sim X_{u_0}$~, $w_0 \in W$ and $u_0 \in \T_g$~, there exists a unique $\varphi: W \rightarrow \T_g$ such that $\varphi\rund{w_0} = u_0$ and $\varphi$ can be extended to a fiberwise biholomorphic family morphism $(\Phi, \varphi): N \rightarrow M_g$ with $\Phi|_{Y_{w_0}} = \sigma$~. If $g \geq 2$ then also $\Phi$ is uniquely determined by $\sigma$~.
\item {\bf Infinitesimal universality:} $\rund{d M_g}_u: T_u \T_g \rightarrow H^1\rund{T X_u}$ is an isomorphism for all $u \in \T_g$~.
\end{itemize}

$\T_0$ is a single point and $M_0$ the Riemann sphere.

For $g = 1$ we have a standard realization: $\T_1 \subset \cz$ the upper half plane and 

\[
\M_1 := \rund{\T_1 \times \cz} / \spitz{S, T}
\]

with $S, T \in \Aut_\strong \ \rund{\T_1 \times \cz}$ given by $(u, z) \mapsto (u, z + 1)$ resp. $(u, z + u)$~.

For $g \geq 2$~, $\T_g$ is of dimension $3(g - 1)$~. \\

As in detail explained in \cite{Earle}, the Jacobian varieties of all $X_u$~, $u \in \T_g$~, glue together to a holomorphic tori bundle $\B_g \twoheadrightarrow \T_g$~. Take a global frame $\rund{\omega_1, \dots, \omega_g}$ of $\rund{\pi_{M_g}}_* \rund{T^\rel M_g}^*$ and a set $\gamma_1, \dots, \gamma_{2 g}$ of piecewise smooth generators of $\pi_1\rund{X_u}$~, which we choose to be independent of the point $u \in \T_g$ in a trivialization of $M_g$ as a smooth family of compact real manifolds. Then we obtain an exact sequence

\[
\begin{array}{rcl}
0 \hookrightarrow \T_g \times \zz^{2 g} & \hookrightarrow & \T_g \times \cz^g \twoheadrightarrow \B_g \twoheadrightarrow 0 \\
\rund{u, e_\mu} & \mapsto & \rund{u, \lambda_\mu(u)}
\end{array}
\]

of holomorphic Lie group bundles over $\T_g$~, where $\lambda_\mu \in \O\rund{\T_g}^{\oplus g}$ is given by \\
$\lambda_\mu(u) := \int_{\gamma_\mu} \left.\rund{\omega_1, \dots, \omega_g}\right|_{X_u}$~. For $u \in \T_g$ let $\Lambda_u \sqsubset \cz^g$ be the lattice generated by $\lambda_\mu(u)$~, $\mu = 1, \dots, 2 g$~, so $\cz^g \left/ \Lambda_u\right.$ is the Jacobian variety of $X_u$~. It is well known, see \cite{Forster} sections 2.21 and 3.29, that for every $u \in \T_g$ we have a group isomorphism









\begin{eqnarray*}
&& \aleph_u: \schweif{\text{isom. classes of hol. line bundles on } X_u \text{ of degree } 0} \mathop{\rightarrow}\limits^\sim \cz^g \left/ \Lambda_u\right. \,~, \\ 
&& \phantom{1234567890} [L] \mapsto \sum_{\kappa = 1}^k \int_{n_\kappa}^{p_\kappa} \left.\rund{\omega_1, \dots, \omega_g}\right|_{X_u} \,~,
\end{eqnarray*}

where $\sum_{\kappa = 1}^k \rund{n_\kappa - p_\kappa}$ is the divisor of an arbitrary non-zero meromorphic section of~$L$~.

\begin{lemma} \label{versallb}
\item[(i)] There exists a holomorphic line bundle $L_g \twoheadrightarrow M_g \times \cz^g$ of degree $0$ such that for every $u \in \T_g$ the map

\[
\cz^g \rightarrow \left.\cz^g \right/ \Lambda_u \,~, \, b \mapsto \aleph_u \eckig{\left.L_g\right|_{X_u \times \schweif{b}}}
\]

is the canonical projection.

For every $(u, b) \in \T_g \times \cz^g$ the differential $\rund{d \ L_g |_{X_u \times \cz^g}}_b : T_b \ \cz^g \rightarrow H^1\rund{X_u, \O}$ is an isomorphism.

\item[(ii)] Let $W$ be a complex manifold, $\varphi: W \rightarrow \T_g$ holomorphic and $E \twoheadrightarrow \varphi^* M_g$ a holomorphic line bundle of degree $0$~. Then

\[
w \mapsto \aleph_{\varphi(w)} \eckig{E|_{X_{\varphi(w)} }}
\]

gives a holomorphic section of $\varphi^* \B_g$~, and if $\psi \in \O(W)^{\oplus g}$ is a lift of this section and $W$ a contractible Stein manifold then $E \simeq \rund{\breve \varphi_{M_g}, \psi \circ \pi_{\varphi^* M_g}}^* L_g$~.
\end{lemma}

{\it Proof:} (ii) {\it First part:} Write $N := \varphi^* M_g$ and $Y_w := \pi_N^{- 1}(w) = X_{\varphi(w)}$ for all $w \in W$~. Let $w_0 \in W$ be arbitrary and $s$ a meromorphic section of $E|_{Y_{w_0}}$~. After maybe shrinking $W$ as an open neighbourhood of $w_0$ we may assume that there exist holomorphic sections $a_\kappa, b_\kappa$~, and $p_\mu \in \Gamma^\hol(W, N)$~, $\kappa = 1, \dots, k$~, $\mu = 1, \dots, 2 g$~, of $N$ such that $\sum_{\kappa = 1}^k \rund{a_\kappa\rund{w_0} - b_\kappa\rund{w_0}}$ is the divisor of $s$ and $p_1\rund{w_0}, \dots, p_{2 g}\rund{w_0} \in Y_{w_0}$ are distinct. Let $F \twoheadrightarrow N$ be the holomorphic line bundle to the divisor

\[
\sum_{\kappa = 1}^k \rund{b_\kappa(W) - a_\kappa(W)} + \sum_{\mu = 1}^{2 g} p_\mu(W)
\]

on $N$~. Then $s$ as holomorphic section of $\left.(E \otimes F)\right|_{Y_{w_0}}$ has the divisor $\sum_{\mu = 1}^{2 g} p_\mu\rund{w_0}$~. On the other hand $h^0\rund{Y_w, (E \otimes F)|_{Y_w}} =~g + 1$ independent of $w \in W$~, so \cite{GrauRem} theorem 10.5.5 tells us that the evaluation map

\[
\rund{\rund{\pi_N}_* \rund{E \otimes F}}_{w_0} \rightarrow H^0\rund{Y_{w_0}, \rund{E \otimes F}|_{Y_{w_0}} }
\]

is surjective, and so after replacing $W$ by a smaller open neighbourhood of $w_0$ there exists \\
$S \in H^0\rund{W, E \otimes F}$ such that $s = S|_{Y_{w_0}} $~. Again after maybe replacing $W$ by a smaller open neighbourhood of $w_0$~, by the implicit function theorem there exist $n_\mu \in \Gamma^\hol(W, N)$~, $\mu = 1, \dots, 2 g$~, such that $S$ as holomorphic section of $E \otimes F$ has the divisor $\sum_{\mu = 1}^{2 g} n_\mu(W)$~, so as a meromorphic section of $E$

\[
\sum_{\kappa = 1}^k \rund{a_\kappa(W) - b_\kappa(W)} + \sum_{\mu = 1}^{2 g} \rund{n_\mu(W) - p_\mu(W)} \,~.
\]

Therefore for all $w \in W$

\[
\aleph_{\varphi(w)} \eckig{E|_{Y_w}} = \rund{\sum_{\kappa = 1}^k \int_{a_\kappa(w)}^{b_\kappa(w)} + \sum_{\mu = 1}^{2 g} \int_{n_\mu(w)}^{p_\mu(w)} } \left.{\breve \varphi_{M_g}}^* \rund{\omega_1, \dots, \omega_g}\right|_{Y_w} \,~,
\]

which is holomorphic in $w$~.

{\it Second part:} Assume $L_g$ to be already constructed on $\left.\rund{M_g \times \cz^g}\right|_V$ and $(\varphi, \psi)(W) \subset V$ for some $V \subset \T_g \times \cz^g$ open. Define the holomorphic line bundle $C := E \otimes \rund{\breve \varphi_{M_g}, \psi \circ \pi_N}^* L_g^*$ on $N$~. For every $w \in W$ seperately $C|_{Y_w}$ is trivial. So by \cite{GrauRem} theorem 10.5.5 and Grauert's theorem since $W$ is contractible and Stein $\rund{\pi_{\varphi^* M_g}}_* C \simeq \O_W$~, which yields a global trivialization of the bundle $C$~. \\

(i) {\it Local construction:} Let $\rund{u_0, b_0} \in \T_g \times \cz^g$ be arbitrary. Take a holomorphic line bundle $E \twoheadrightarrow \left.M_g\right|_U$~, $U \subset \T_g$ a suitable open neighbourhood of $u_0$~, such that $\aleph_{u_0} \eckig{E|_{X_{u_0}} } = b_0$~. By \cite{Forster} lemma 2.21.3 we may assume that after shrinking $U$ there are $n_1, \dots, n_g \in \Gamma^\hol(U, M_g)$ such that $\det\rund{\left.\omega_\nu\right|_{X_{u_0}}\rund{n_\rho\rund{u_0}} }$ is invertible. Let $M_g^{\times_{\T_g} g} \twoheadrightarrow \T_g$ denote the fiberwise $g$-fold direct product (it is the obvious family over $\T_g$ with fiber $\underbrace{X_u \times \cdots \times X_u}_{g \text{ times}}$ at the point $u \in \T_g$~). Then the strong holomorphic morphism

\begin{eqnarray*}
\Phi: \left.M_g^{\times_{\T_g} g}\right|_U &\rightarrow& \left.\B_g\right|_U \,~, \\
\rund{u, z_1, \dots, z_r} &\mapsto& \rund{u, \aleph_u \eckig{E|_{X_u}} + \sum_{\rho = 1}^g \int_{n_\rho(u)}^{z_\rho} \left.\rund{\omega_1, \dots, \omega_g}\right|_{X_u}} \,~,
\end{eqnarray*}


of families over $U$~, mapping $\rund{n_1\rund{u_0}, \dots, n_1\rund{u_0}} \in X_{u_0}^g$ to $\overline{\rund{u_0, b_0}}$~, has invertible dif\-fe\-ren\-tial at the point $\rund{n_1\rund{u_0}, \dots, n_g\rund{u_0}}$~. So let $\Psi~=~\rund{q_1, \dots, q_g}: V \rightarrow M_g^{\times_{\T_g} g}$~, $V \subset U \times \cz^g$ an open neighbourhood of $\rund{u_0, b_0}$ and $q_\mu: V \rightarrow M_g$ holomorphic, be a local inverse. Then $q_\mu(u, b) \in X_u$ for every $(u, b) \in V$~. Finally define $L_g := F \otimes \pr_{M_g}^* E$ where $F \twoheadrightarrow \left.\rund{M_g \times \cz^g}\right|_V$ denotes the holomorphic line bundle associated to the divisor 

\[
\sum_{\mu = 1}^g \rund{n_\mu(U) \times \cz^g - \rund{q_\mu, \pr_{\cz^g}}(V)}
\]

on $\left.\rund{M_g \times \cz^g}\right|_V$~. Then for $(u, b) \in V$ 

\[
\aleph_u \eckig{\left.L_g\right|_{X_u \times \schweif{b}} } = \Phi\rund{q_1(u, b), \dots, q_g(u, b)} = \overline{(u, b)} \,~.
\]

$d L_g |_{X_{u_0} \times \cz^g}$ is an isomorphism by standard Riemann surface theory.

{\it Global construction:} By the local construction we obtain a cover $\T_g \times \cz^g = \bigcup_i V_i$ by open balls and holomorphic line bundles $L_i \twoheadrightarrow \left.\rund{M_g \times \cz^g}\right|_{V_i}$ such that $\aleph_u \eckig{\left.L_i\right|_{X_u \times \schweif{b}} } = \overline{(u, b)}$ for all $i$ and $(u, b) \in V_i$~. By the proof of (ii) second part $L_i \simeq L_j$ on the overlaps $\left.\rund{M_g \times \cz^g}\right|_{V_i \cap V_j}$~, and so the obstructions for the $L_i$ glueing together to a holomorphic line bundle $L_g \twoheadrightarrow M_g \times \cz^g$ lie in $H^2\rund{\T_g \times \cz^g, \O^\times}$~, which is trivial by Cartan's theorem B since $\T_g \times \cz^g$ is contractible and Stein. $\Box$






\begin{lemma} \label{deg1} There exists a holomorphic line bundle $E_g \twoheadrightarrow M_g$ of degree $1$ such that $E_g^{\otimes 2 (1 - g)} = T^\rel M_g$~. \end{lemma}

{\it Proof: Local construction:} Let $u_0 \in \T_g$ and choose an open neighbourhood $U \subset \T_g$ of $u_0$ and a holomorphic line bundle $E \twoheadrightarrow \left.M_g\right|_U$ of degree $1$~. If $g = 1$ define $E_g := E$~, and we are done since $T^\rel M_1$ is trivial. For $g \not= 1$ define the strong fiberwise affine linear bundle epimorphism

\[
\Phi: U \times \cz^g \rightarrow \left.\B_g\right|_U \,~, \, (u, b) \mapsto \aleph_u \eckig{\rund{\left.L_g\right|_{X_u \times \schweif{b}} \otimes E|_{X_u}}^{\otimes 2 (g - 1)} \otimes T X_u} \,~.
\]

After shrinking $U$ to a small enough open ball there exists $\rho~\in~\O(U)^{\oplus g}$ such that \\
$\Phi \circ \rund{\id_U, \rho} = \rund{\id_U, 0}$~, and so by lemma \ref{versallb} (ii)

\[
\rund{(\id, \rho)^* L_g \otimes E}^{\otimes 2 (g - 1)} \otimes \left.T^\rel M_g\right|_U
\]

is trivial. Define $E_g := (\id_U, \rho)^* L_g \otimes E$~.

{\it Global construction:} By the local construction we obtain a cover $\T_g = \bigcup_i U_i$ by open balls and holomorphic line bundles $E_i \twoheadrightarrow \left.M_g \right|_{U_i}$ of degree $1$ such that $E_i^{\otimes 2 (1 - g)} = T^\rel M_g$~. Define $\alpha_{i j} \in H^0\rund{U_i \cap U_j, \B_g}$ as

\[
u \mapsto \aleph_u\eckig{\left.\rund{E_i \otimes E_j^*}\right|_{X_u}} \,~.
\]

$\rund{\alpha_{i j}}$ is a $1$-cocycle in $\B_g$~. $H^1\rund{\T_g, \B_g} = 0$ by Cartan's theorem B since $\T_g$ is contractible and Stein, and so by passing to a refinement of the cover we obtain $\beta_i \in \O\rund{U_i}^{\oplus g}$ such that $\alpha_{i j} = \overline{\beta_i - \beta_j}$~. For $g \not= 1$~, $\rund{\alpha_{i j}}$ is infact a $1$-cocycle in the subbundle \\
$\dot\bigcup_{u \in \T_g} \left.\frac{1}{2 (g - 1)} \Lambda_u \right/ \Lambda_u \simeq \T_g \times \zz_{2 (g - 1)}^{2 g}$ of $\B_g$~. Since $\T_g$ is contractible also $H^1\rund{\T_g, \zz_{2 (g - 1)} } = 0$~, and so we may assume that infact all $\beta_i(u) \in \frac{1}{2 (g - 1)} \Lambda_u$~.

By taking $E_i \otimes \rund{\id_{\left.M_g\right|_{U_i}}, \beta_i}^* L_g^*$ instead of $E_i$ we may assume that $\aleph_u \eckig{\left.\rund{E_i \otimes E_j^*}\right|_{X_u}} = 0$ for all $i, j$ and $u \in U_i \cap U_j$~, and so by theorem \ref{versallb} (ii) there exist isomorphisms $\varphi_{i j}: E_i \mathop{\rightarrow}\limits^\sim E_j$ on the overlaps $\left.M_g\right|_{U_i \cap U_j}$ with $\varphi_{i j}^{\otimes 2 (g - 1)} = \id_{T^\rel M_g}$~.

We see that the obstructions for the $E_i$ glueing together to a holomorphic line bundle $E_g \twoheadrightarrow M_g$ lie in $H^2\rund{\T_1, \O^\times}$ if $g = 1$ resp. $H^2\rund{\T_g, \zz_{2 (g - 1)} }$ if $g \not= 1$~, which are both trivial since $\T_g$ is contractible and Stein, the first by Cartan's theorem B. $\Box$ \\

We are now prepaired for defining the family $\M_{g, D} \twoheadrightarrow V_{g, D}$~, and we have to distinguish two cases:

\begin{itemize}

\item For $g \not= 1$ or $D = 0$ we define the holomorphic family

\[
\pi_{M_{g, D}} := \rund{\pi_{M_g}, \id_{\cz^g}}: M_{g, D} := M_g \times \cz^g \twoheadrightarrow V_{g, D} := \T_g \times \cz^g
\]

of compact Riemann surfaces of genus $g$ and the holomorphic line bundle \\
$L_{g, D} := L_g \otimes \pr_{M_g}^* E_g^{\otimes D} \twoheadrightarrow M_{g, D}$ of degree $D$~. Of course $V_{0, D}$ is just a single point.

\item For $D \not= 0$ we define the holomorphic familiy $M_{1, D} := M_1 \twoheadrightarrow V_{1, D} := \T_1$ of compact Riemann surfaces of genus $1$ and the holomorphic line bundle $L_{1, D} := E_1^{\otimes D} \twoheadrightarrow M_{1, D}$ of degree $D$~.

\end{itemize}

In both cases we define the family

\[
\pi_{\M_{g, D}}: \M_{g, D} := \rund{M_{g, D}, \Gamma^\hol\rund{\bigwedge L_{g, D}^*}} \twoheadrightarrow V_{g, D}
\]

of compact super Riemann surfaces $\X_v$ of type $(g, D)$~. \\


Why do we have to distinguish these two cases? Well, as a holomorphic Lie group bundle $M_1$ acts on itself by translations $\tau_a: (u, z) \mapsto (u, z + a)$ and so on the isomorphy classes of holomorphic line bundles of degree $D$ over its fibers by pullback. The crucial point is that this action is trivial if $D = 0$ and fiberwise transitive if $D \not= 0$~:

\begin{lemma} \label{genus1}
\item[(i)] $\tau_a^* \left.L_1\right|_{X_u \times \schweif{b}} \simeq \left.L_1\right|_{X_u \times \schweif{b}}$ for all $u \in \T_1$~, $b \in \cz$ and $a \in X_u$~.

\item[(ii)] $M_1 \rightarrow \B_1 \,~, \, \overline{(u, a)} \mapsto \aleph_u \eckig{\rund{\tau_a^* \left.E_1\right|_{X_u}} \otimes \left.E_1^*\right|_{X_u}}$ is a strong Lie group bundle isomorphism.

\item[(iii)] Let $\X$ be a compact super Riemann surface of type $(1, D)$ with body $X_u$~, $u \in \T_1$~. Then in (\ref{exacttangent}) of section \ref{basics} with $\M := \X$ the vectorfield $\partial_z \in H^0\rund{T X_u}$ has a global preimage in $H^0\rund{sT \X}_0$ iff $D = 0$~.
\end{lemma}

{\it Proof:} (i) Take a meromorphic section $S$ of $\left.L_1\right|_{X_u \times \schweif{b}}$ with one single pole, and let $z_0 - w_0$ be its divisor. Then $\rund{S \circ \tau_a} S^{- 1}$ is a meromorphic section of $\rund{\tau_a^* \left.L_1\right|_{X_u \times \schweif{b}} } \otimes \left.L_1\right|_{X_u \times \schweif{b}}$ with divisor $\rund{z_0 - a} - \rund{w_0 - a} - z_0 + w_0$~, and so

\[
\aleph_u \eckig{\rund{\tau_a^* \left.L_1\right|_{X_u \times \schweif{b}} } \otimes \left.L_1\right|_{X_u \times \schweif{b}}} = \rund{\int_{z_0 - a}^{z_0} + \int_{w_0}^{w_0 - a}} d z = 0 \,~.
\]

(ii) Let $u \in \T_1$~, $S \in H^0\rund{X_u, \left.E_1\right|_{X_u}}$ and $z_0 \in X_u$ be the divisor associated to $S$~. Then $\rund{S \circ \tau_a} S^{- 1}$ is a meromorphic section of $\rund{\tau_a^* \left.E_1\right|_{X_u}}~\otimes~\left.E_1^*\right|_{X_u}$ with divisor $\rund{z_0 - a} - z_0$~. So

\[
\aleph_u \eckig{\rund{\tau_a^* \left.E_1\right|_{X_u}} \otimes \left.E_1^*\right|_{X_u}} = \int_{z_0 - a}^{z_0} d z = a \,~.
\]

(iii) Write $\X = \rund{X_u, \Gamma^\hol\rund{\bigwedge L^*}}$ with a holomorphic line bundle $L \twoheadrightarrow X$ of degree $D$~.

`$\Rightarrow$': Let $\delta \in H^0\rund{sT \X}_0$ such that $\delta|_{\O_{X_u}} = \partial_z$ and $S$ be a global meromorphic section of $L^*$ (therefore away from its poles an odd holomorphic function on $\X$ ) such that $\delta S \not= 0$~. For every $a \in X_u$ we can find a local super chart $U^{|1}$~, $U \subset \cz$ open, of $\X_u$ around $a$ in which $S = (z - a)^{\ord_a S} \zeta$ and $\delta = \partial_z + h \zeta \partial_\zeta$ for some $h \in \O(U)$~. Therefore in this super chart

\[
\delta S = \rund{\frac{\ord_a \ S}{z - a} + h} S~.
\]

We see that $\delta S$ is again a global meromorphic section of $L^*$~, so $f := \rund{\delta S} S^{- 1} \in \M\rund{X_u}$~. Furthermore $\Res_a \ f = \ord_a \ S$ for all $a \in X_u$~, and so

\[
0 = \sum_{a \in X_u} \Res_{a} f = - D \,~.
\]

`$\Leftarrow$': Assume $D = 0$~. Then $L$ admits an everywhere flat connection, which gives local super charts of $\X$ with transition morphisms of the form $(z + a, A \zeta)$~, $a \in \cz$~, $A \in \cz^\times$~. Therefore $\partial_z$ in these local super charts glue together to an element of $H^0\rund{sT \X}_0$~. $\Box$ \\

\begin{theorem}[Versality of $\M_{g, D}$ ] \label{versalfam}

\item[(i)] {\bf Completeness:} Let $W$ be a contractible Stein manifold and $\N \twoheadrightarrow W$ a family of compact super Riemann surfaces $\Y_w$ of type $(g, D)$~. Then there exists a fiberwise biholomorphic family morphism $(\Xi, \xi): \N \rightarrow \M_{g, D}$~. More precisely we have the following

{\bf Anchoring property (including local completeness):} For every isomorphism \\
$\sigma: \Y_{w_0} \mathop{\rightarrow}\limits^\sim \X_{v_0}$~, $w_0 \in W$ and $v_0 \in V_{g, D}$~, there exists $\xi: W \rightarrow~V_{g, D}$ holomorphic such that $\xi\rund{w_0} = v_0$ and $\xi$ can be extended to a fiberwise biholomorphic family morphism \\
$(\Xi, \xi): \N \rightarrow \M_{g, D}$ with $\Xi|_{\Y_{w_0}} = \sigma$~. $\xi$ is uniquely determined by $\sigma^\#$~.

\item[(ii)] {\bf Infinitesimal universality:} $\rund{d \M_{g, D}}_v : T_v V_{g, D} \rightarrow H^1\rund{sT \X_v}_0$ is an isomorphism for all $v \in V_{g, D}$~.
\end{theorem}

Observe that for every ordinary family $\M \twoheadrightarrow V$ of compact complex supermanifolds $\X_v$ and $v \in V$ the image of $(d \M)_v$ automatically lies in $H^1\rund{sT \X_v}_0$~. So (ii) is the best infinitesimal universality we can expect in the framework of ordinary families. In section \ref{superfam} we will construct super families $\widetilde \M_{g, D} \twoheadrightarrow V_{g, D}^{|q}$~, $q$ appropriate, such that $\rund{d \widetilde \M_{g, D}}_v: sT_v V_{g, D}^{|q} \rightarrow H^1\rund{sT \X_v}$ is an isomorphism for every $v \in V_{g, D}$~. \\

{\it Proof:} (i) For proving the {\it anchoring property} write $\N = \rund{N, \Gamma^\hol\rund{\bigwedge F^*}}$ and $v_0 = \rund{u_0, b_0}$~, $u_0 \in \T_g$~, $b_0 \in \cz^g$~, if $g \not= 1$ or $D = 0$ resp. $u_0 := v_0 \in \T_1$ if $g = 1$ and $D \not= 0$~. By the anchoring property of $M_g$ there exists a unique $\varphi: W \rightarrow \T_g$ holomorphic such that $\varphi\rund{w_0} = u_0$ and $\varphi$ can be extended to a fiberwise biholomorphic family morphism $\rund{\Phi, \varphi}: N \rightarrow M_g$ with $\Phi|_{\Y_{w_0}^\#} = \sigma^\#$~. So without restriction we may assume that $N = \varphi^* M_g$ and $\Phi = \breve \varphi_{M_g}$~. By lemma \ref{versallb} (ii)

\[
w \mapsto \aleph_{\varphi(w)} \eckig{\left.\rund{F \otimes \rund{\Phi^* E_1}^{\otimes (- D)} }\right|_{X_{\varphi(w)} }}
\]

gives a holomorphic section $f \in H^0\rund{W, \varphi^* \B_g}$ with $f\rund{w_0} = \overline{\rund{u_0, b_0}}$ if $g \not= 1$ or $D = 0$ resp. $f\rund{w_0} = \overline{\rund{u_0, 0}}$ if $g = 1$ and $D \not= 0$~. Let $\psi \in \O(W)^{\oplus g}$ be a holomorphic lift of $f$ with $\psi\rund{w_0} = b_0$ if $g \not= 1$ or $D = 0$ resp. $\psi\rund{w_0} = 0$ if $g = 1$ and $D \not= 0$~. \\

{\it Existence:} If $g \not= 1$ or $D = 0$ by lemma \ref{versallb} (ii) since $W$ is contractible and Stein $F \simeq \rund{\Phi, \psi \circ \pi_N}^* L_{g, D}$~, which yields a fiberwise biholomorphic family morphism $(\Xi, \xi): \N \rightarrow \M_{g, D}$~, $\xi := \rund{\varphi, \psi}$~, having $\Xi^\# = (\Phi, \psi)$~, and so $\rund{\Xi|_{\Y_{w_0}} }^\# = \sigma^\#$~. \\

If $g = 1$ and $D \not= 0$ by lemma \ref{genus1} (ii) for all $w \in W$

\begin{eqnarray*}
&& \aleph_{\varphi(w)} \eckig{\left.\rund{F \otimes \rund{\Phi \circ \tau_{\frac{\psi}{D}} }^* E_1^{\otimes (- D)} }\right|_{X_{\varphi(w)} }} \\
&& \phantom{12} = f(w) - D \ \aleph_{\varphi(w)} \eckig{\rund{\left.\tau_{\frac{\psi(w)}{D}}^* E_1\right|_{X_{\varphi(w)} }} \otimes \left.E_1^*\right|_{X_{\varphi(w)} }} = 0
\end{eqnarray*}

and so $F \simeq \rund{\Phi \circ \tau_{\frac{\psi}{D}} }^* E_1^{\otimes D}$ by lemma \ref{versallb} (ii), which yields a fiberwise biholomorphic family morphism $(\Xi, \xi): \N \rightarrow \M_{1, D}$~, $\xi := \varphi$~, having $\Xi^\# = \Phi \circ \tau_{\frac{\psi}{D}}$~. $\rund{\Xi|_{\Y_{w_0}} }^\# = \sigma^\#$ since $\psi\rund{w_0} = 0$~. \\

Finally in both cases $\sigma = \Phi_{\X_{w_0}, a} \circ \Xi|_{\Y_{w_0}}$ for some $a \in \cz^\times$ by (\ref{exact}) in section \ref{basics}, so take $\Phi_{\N, a} \circ \Xi$ instead of $\Xi$ to obtain even $\sigma = \Xi|_{\Y_{w_0}}$~. \\

{\it Uniqueness:} trivial if $g = 0$~. So assume $g \geq 1$~. Let $\rund{\Xi, \xi}: \N \rightarrow \M_{g, D}$ be a fiberwise biholomorphic family morphism such that $\xi\rund{w_0} = v_0$ and $\rund{\Xi|_{\Y_{w_0}} }^\# = \sigma^\#$~. Then $\rund{\pr_{M_g} \circ \Xi^\#, \pr_{\T_g} \circ \xi}: N \rightarrow M_g$ is a fiberwise biholomorphic family morphism.  \\

If $g = 1$ and $D \not= 0$ we obtain $\xi = \varphi$ by the anchoring property of $M_1$~. \\

If $g \not= 0$ or $D = 0$ then by the anchoring property of $M_g$ we obtain $\xi = \rund{\varphi, \widehat \psi}$ with some $\widehat \psi \in \O(W)^{\oplus g}$~.

If $g \geq 2$ then even $\Xi^\# = \rund{\Phi, \widehat \psi}$ by the anchoring property of $M_g$~. Therefore \\
$F~\simeq~\rund{\Phi, \widehat \psi \circ \pi_N}^*~L_{g, D}$~. 

If $(g, D) = (1, 0)$ we obtain $\Xi^\# = \rund{\Phi \circ \tau_h, \widehat \psi}$ with some $h \in \O(W)$~, therefore \\
$F \simeq \tau_h^* \rund{\Phi, \widehat \psi \circ \pi_N}^* L_1$~. But by lemma \ref{genus1} (i) for all $w \in W$

\[
f(w) = \aleph_{\varphi(w)} \eckig{F|_{X_{\varphi(w)} }} = \aleph_{\varphi(w)} \eckig{\tau_{h(w)}^* \left.L_1\right|_{X_{\varphi(w)} \times \schweif{\widehat \psi(w)} }} = \aleph_{\varphi(w)} \eckig{\left.L_1\right|_{X_{\varphi(w)} \times \schweif{\widehat \psi(w)} }} = \overline{\widehat \psi(w)} \,~.
\]

In both cases we see that $\widehat \psi$ and $\psi$ are holomorphic lifts of $f$~, which coincide at $w_0$ and so must be equal. \\

For proving {\it completeness} now it suffices to check that every compact super Riemann surface of type $(g, D)$ is represented by at least one fiber in the family $\M_{g, D}$~, which can be easily done by similar calculations as in the proof of the existence in the anchoring property. \\

(ii) For $g \not= 1$ or $D = 0$ and $v = (u, b) \in V_{g, D}$ we obtain

\[
\begin{array}{ccccc}
H^0\rund{sT \X_v} \mathop{\longrightarrow}\limits^{(*)} H^0\rund{T X_u} \rightarrow H^1\rund{X_u, \O} & \mathop{\longrightarrow}\limits^{(**)} & \phantom{12345} H^1\rund{sT \X_v}_0 & \twoheadrightarrow & H^1\rund{T X_u} \twoheadrightarrow 0 \phantom{12} \\
\phantom{1234567890} \rund{d \ \left.L_{g, D}\right|_{X_u \times \cz^g}}_b \uparrow & \circlearrowleft & \rund{d \M_{g, D}}_v \uparrow \phantom{12345} & \circlearrowleft & \uparrow \rund{d M_g}_u \\
\phantom{1234567890123456,} 0 \phantom{123,} \hookrightarrow \phantom{123} T_b \ \cz^g \phantom{1} & \hookrightarrow & \phantom{12345} T_v V_{g, D} & \twoheadrightarrow & T_u \T_g \phantom{12,} \twoheadrightarrow 0
\end{array}
\]

and for $g = 1$~, $D \not= 0$ and $v = u \in V_{1, D} = \T_1$

\[
\begin{array}{ccc}
H^0\rund{sT \X_u} \mathop{\longrightarrow}\limits^{(*)} H^0\rund{T X_u} \rightarrow H^1\rund{X_u, \O} \mathop{\longrightarrow}\limits^{(**)} H^1\rund{sT \X_u}_0 & \twoheadrightarrow & H^1\rund{T X_u} \twoheadrightarrow 0 \\
\phantom{1234567890123456789012345678} \rund{d \M_{1, D}}_u \uparrow & \circlearrowleft & \phantom{12} \uparrow \rund{d M_1}_u \\
\phantom{12345678901234567890123456789012345678} T_u V_{1, D} & = & T_u \T_1 \phantom{1234}
\end{array} \,~,
\]

where in both cases the rows are exact, $\rund{d \ \left.L_{g, D}\right|_{X_u \times \cz^g}}_b$ and $\rund{d M_g}_u$ are isomorphisms, and the second rows are induced by the short exact sequence (\ref{exacttangent}) in section \ref{basics} with $\M := \X_v$~. \\

If $g = 0$ then $V_{0, D}$ is a point, and $H^1\rund{X_u, \O} = H^1\rund{T X_u} = 0$~. \\

If $g \geq 2$ then $H^0\rund{T X_u} = 0$~. If $g = 1$ and $D = 0$ then $(*)$ is surjective by lemma \ref{genus1} (iii) since $h^0\rund{T X_u} = h^1\rund{X_u, \O} = 1$~. In both cases we see that $(**)$ is injective. Therefore by the five lemma also $\rund{d \M_{g, D}}_v$ is an isomorphism. \\

If $g = 1$ and $D \not= 0$ then by lemma \ref{genus1} (iii) we see that $(*)$ and so also $(**)$ is the zero map. $\Box$ \\

Recall the classical modular group $\Gamma_g \sqsubset \Aut \ \T_g$ for genus $g$~, so the image of the canonical projection $\Aut \ M_g \twoheadrightarrow \Aut \ \T_g \,~, \, (\Gamma, \gamma) \mapsto \gamma$~. It acts properly discontinuously on $\T_g$~, see \cite{Duma}, $\T_g \left/ \Gamma_g\right.$ is the moduli space of compact Riemann surfaces of genus $g$~, and $\Gamma_1 = P SL(2, \zz)$ acts on $\T_1$ by Möbius transformations. By simple calculations, \cite{Ares}, \cite{Baily}, \cite{Duma}, \cite{Hurwitz}, and \cite{Poonen} we know that \\

$\schweif{\tau_f \, \left| \, f \in H^0\rund{\T_1, \B_1}\right.} \subset \Aut_\strong \ M_1$ is a normal subgroup, and $\schweif{\id_{M_1}, J}$~, $J \in \Aut_\strong \ M_1$ given by $(u, z) \mapsto (u, - z)$~, is a set of representatives for $\rund{\Aut_\strong \ M_1} \left/ \schweif{\tau_f \, \left| \, f \in H^0\rund{\T_1, \B_1}\right.}\right.$~, 

$\Aut_\strong \ M_2 = \schweif{\id_{M_2}, J}$ with the hyperelliptic reflection $J$~, and

$\Aut_\strong \ M_g = \schweif{\id_{M_g}}$ for $g \geq 3$~.

\begin{defin}
$V_{g, D}$ is called the Teichmüller space and the image $\Gamma_{g, D}$ of the group homomorphism

\[
\Aut \ \M_{g, D} \rightarrow \Aut \ V_{g, D} \,~, \, (\Gamma, \gamma) \mapsto \gamma
\]

the modular group for type $(g, D)$~.
\end{defin}

As in the classical case with the help of the anchoring property of $\M_{g, D}$ one proves:

\begin{cor} \label{ordfamunique}
\item[(i)] Every fiberwise biholomorphic family endomorphism of $\M_{g, D}$ is infact an automorphism.
\item[(ii)] The family $\M_{g, D} \twoheadrightarrow V_{g, D}$ is up to family isomorphism uniquely determined by the completeness and infinitesimal universality of theorem \ref{versalfam} and $V_{g, D}$ being a contractible Stein manifold.
\item[(iii)] $V_{g, D} / \Gamma_{g, D}$ is the moduli space for type $(g, D)$~. More precisely: Let $W$ be a connected complex manifold and $\xi, \tau: W \rightarrow V_{g, D}$ holomorphic. Then there exists a strong family isomorphism $\xi^* \M_{g, D} \simeq \tau^* \M_{g, D}$ iff there exists $\gamma \in \Gamma_{g, D}$ such that $\tau = \gamma \circ \xi$~.
\end{cor}








The rest of this section is devoted to the investigation of the modular group $\Gamma_{g, D}$~. \\

Let $g \not= 1$ or $D = 0$~. By a first observation there is an embedding $\zz^{2 g} \hookrightarrow \Gamma_{g, D}$~. Indeed, for $i \in \schweif{1, \dots, 2 g}$ by lemma \ref{versallb} (ii) $\rund{\id_{M_g}, \id_{\cz^g} + \lambda_i}^* L_{g, D} \simeq L_{g, D}$~, which leads to a family automorphism $\rund{\Gamma, \id_{V_{g, D}} + \lambda_i}$ of $\M_{g, D}$ with $\Gamma^\# = \rund{\id_{M_g}, \id_{\cz^g} + \lambda_i \circ \pi_{M_g}}$~. How big might be the quotient? $\Aut \ M_g$ acts on the isomorphy classes of holomorphic line bundles of degree $D$ by pullback, so on $\B_g$ by bundle automorphisms: For every $(\Delta, \delta) \in \Aut \ M_g$~, $u \in \T_g$ and holomorphic line bundles $F \twoheadrightarrow X_u$ and $K \twoheadrightarrow X_{\delta(u)}$ of degree $D$

\[
(\Delta, \delta) \aleph_u \eckig{F \otimes \left.E_g^{\otimes (- D)}\right|_{X_u}} = \aleph_{\delta(u)} \eckig{K \otimes \left.E_g^{\otimes (- D)}\right|_{X_\delta(u)} }
\]

iff $F \simeq \Delta|_{X_u}^* K$~. With the help of the line bundle $L_{g, D} \twoheadrightarrow V_{g, D}$ we can write this action as

\[
(\Delta, \delta) \overline{(u, b)} := \aleph_{\delta(u)} \eckig{\rund{\left.\Delta^{- 1}\right|_{X_{\delta(u)}} }^* \rund{\left.L_{g, D}\right|_{X_u \times \schweif{b}} } \otimes \left.E_g^{\otimes (- D)}\right|_{X_{\delta(u)}} }
\]

for all $(\Delta, \delta) \in \Aut \ M_g$ and $(u, b) \in \T_g \times \cz^g$ and observe that it is fiberwise affine linear. Let $Q_{g, D} \sqsubset \Aut_\strong \ M_g$ be the kernel of this action. \\

\begin{theorem} \label{supermodexplicit}
\item[(i)] Let $g \not= 1$ or $D = 0$~. Then $\Gamma_{g, D}$ acts properly discontinuously on $V_{g, D} \twoheadrightarrow \T_g$ by fiberwise affine linear bundle automorphisms, and

\[
\begin{array}{ccccccc}
& e_\mu & \mapsto & \rund{\id_{\T_g}, \id_{\cz^g} + \lambda_\mu} &&& \\
0 \hookrightarrow & \zz^{2 g} & \hookrightarrow & \Gamma_{g, D} & \twoheadrightarrow & \left.\rund{\Aut \ M_g} \right/ Q_{g, D} & \twoheadrightarrow 1 \\
& \phantom{12} \twoheadsearrow & \circlearrowleft & \phantom{1234567} \twoheadsearrow & \circlearrowleft & \twoheadswarrow \pi_{g, D} \phantom{12345} & \\
&& 1 & \hookrightarrow & \Gamma_g &&
\end{array} \,~,
\]

where the first row is exact and the projection onto $\left.\rund{\Aut \ M_g} \right/ Q_{g, D}$ is given by the action of $\Aut \ M_g$ on $\B_g$ described before. The projection $\Gamma_{g, D} \twoheadrightarrow \Gamma_g$ is given by the underlying action on $\T_g$ and $\pi_{g, D}$ by the canonical projection $(\Gamma, \gamma) \mapsto \Gamma$~.

\begin{itemize}
\item $Q_{1, 0} = \schweif{\tau_f \, \left| \, f \in H^0\rund{\T_1, \B_1}\right.}$~, and so $\pi_{1, 0}$ is a double cover,
\item $Q_{2, D} = \schweif{\id_{M_2}}$~, and so $\pi_{2, D}$ is a double cover,
\item $Q_{g, D} = \schweif{\id_{M_g}}$~, and so $\pi_{g, D}$ is an isomorphism for $g \geq 3$~.
\end{itemize}

\item[(ii)] $\Gamma_{1, D} = \Gamma_1 = PSL(2, \zz)$ for $D \not= 0$~.
\end{theorem}

{\it Proof:} (i) Let $(\Gamma, \gamma) \in \Aut \ \M_{g, D}$~. We would like to show that $\gamma \in \Gamma_{g, D}$ projects to the automorphism on $\B_g$ given by some $(\Delta, \delta) \in \Aut \ M_g$~. $\rund{\Gamma^\#, \gamma} \in \Aut \ M_{g, D}$~, and so for every $b \in \cz^g$

\[
\begin{array}{ccccccc}
M_g & \mathop{\hookrightarrow }\limits^{\rund{\id_{M_g}, b}} & M_{g, D} & \mathop{\longrightarrow}\limits^{\Gamma^\#} & M_{g, D} & \mathop{\twoheadrightarrow}\limits^{\pr_{M_g}} & M_g \\
\twoheaddownarrow & \circlearrowleft & \twoheaddownarrow & \circlearrowleft & \twoheaddownarrow & \circlearrowleft & \twoheaddownarrow \\
\T_g & \mathop{\hookrightarrow }\limits_{\rund{\id_{\T_g}, b}} & V_{g, D} & \mathop{\longrightarrow}\limits_{\gamma} & V_{g, D} & \mathop{\twoheadrightarrow}\limits_{\pr_{\T_g}} & \T_g
\end{array} \,~,
\]

where the vertical maps are the family projections. $\rund{\Delta_b, \delta_b} \in \Aut \ M_g$~, given by the first resp. second row, depends holomorphically on $b$~. So if $g \geq 2$ then $(\Delta, \delta) := \rund{\Delta_b, \delta_b}$ infact does not depend on $b$~at all. If $g = 1$ then $\rund{\Delta_b, \delta_b} = \rund{\tau_{f(\diamondsuit, b)} \circ \Delta, \delta}$ with some $(\Delta, \delta) \in \Aut \ M_g$ and $f \in \O\rund{\T_1 \times \cz}$~. In both cases we see that $\gamma$ is a bundle automorphism of $V_{g, D} \twoheadrightarrow \T_g$~, so $\gamma = \rund{\delta, \kappa}$ with some $\kappa \in \O\rund{V_{g, D}}^{\oplus g}$~. Furthermore $\Gamma|_{\X_{u, b}}$ induces an isomorphism

\[
\left.L_{g, D}\right|_{X_u \times \schweif{b}} \simeq \rund{\Delta|_{X_u}}^* \left.L_{g, D}\right|_{X_{\delta(u)} \times \schweif{\kappa(u, b)}}
\] 

if $g \geq 2$~, resp.

\[
\left.L_1\right|_{X_u \times \schweif{b}} \simeq \rund{\Delta|_{X_u}}^* \tau_{f(u, b)}^* \left.L_1\right|_{X_{\delta(u)} \times \schweif{\kappa(u, b)}} \simeq \rund{\Delta|_{X_u}}^* \left.L_1\right|_{X_{\delta(u)} \times \schweif{\kappa(u, b)}}
\] 

if $g = 1$ and $D = 0$~, which in both cases implies $\overline{\gamma(u, b)} = \rund{\Delta, \delta} \overline{(u, b)}$~. \\

Conversely, let $(\Delta, \delta) \in \Aut \ M_g$~. Take a holomorphic lift $\gamma = \rund{\delta, \kappa} \in \Aut \rund{\T_g \times \cz^g}$~, $\kappa \in \O\rund{\T_g \times \cz^g}^{\oplus g}$ of the bundle automorphism of $\B_g$ corresponding to $(\Delta, \delta)$~. Then $L_{g, D} \simeq \rund{\Delta, \kappa \circ \rund{\pi_{M_g}, \id_{\cz^g}} }^* L_{g, D}$ by lemma \ref{versallb} (ii), which yields a family automorphism $(\Gamma, \gamma)$ of $\M_{g, D}$~. \\

Finally easy calculations show that if $g = 1$ or $g = 2$ then $\aleph_u \eckig{\rund{J|_{X_u}}^* \left.L_g\right|_{X_u \times \schweif{b}} } = - \overline{(u, b)}$~, and so $J \notin Q_{1, 0}$ resp. $J \notin Q_{2, D}$~. In the case $g = 2$ one just has to use the realization of an arbitrary $X_u$~, $u \in \T_2$~, as the projective algebraic curve $w^2 t^4 = \rund{z - c_1 t} \cdots \rund{z - c_6 t}$~, $c_1, \dots, c_6 \in \cz$ distinct, see \cite{FarkKra} section V.2: In this realization $\rund{\frac{d z}{w}~, \frac{z \ d z}{t w} }$ is a basis of $H^0\rund{T^* X_u}$~, and $J|_{X_u}$ is given by $[z : w : t] \mapsto [z : - w : t]$~. \\

(ii) Let $(\Gamma, \gamma) \in \Aut \ \M_{1, D}$~. Then $\rund{\Gamma^\#, \gamma} \in \Aut \ M_1$~, and so $\gamma \in \Gamma_g$~. \\

Conversely, let $(\Delta, \delta) \in \Aut \ M_1$ and $u_0 \in \T_1$ be arbitrary. Then by lemma \ref{genus1} (ii) after re\-pla\-cing $\Delta$ by $\tau_a \circ \Delta$~, $a \in \cz$ appropriate, we may assume that we have an isomorphism $\X_u \mathop{\rightarrow}\limits^\sim \X_{\delta\rund{u_0}}$ with body $\Delta|_{X_{u_0}}$~. By the anchoring property of $\M_{1, D}$ and corollary \ref{ordfamunique} (i) there exists $(\Gamma, \gamma) \in \Aut \ \M_{1, D}$ such that $\gamma\rund{u_0} = \delta\rund{u_0}$ and $\left.\Gamma^\#\right|_{X_{u_0}}~=~\Delta|_{X_{u_0}}$~. So $\delta = \gamma$ by the anchoring property of $M_1$~. $\Box$

\begin{cor} \label{autstrong} $\Aut_\strong \ \M_{g, D} = \schweif{\Phi_{\M_{g, D}, f} \, \left| \, f \in \O\rund{V_{g, D}}^\times\right.}$ for $g \geq 2$~.
\end{cor}

\section{Ordinary supersymmetric families} \label{ordfamSUSY}

In this section we write the dual $\M_{g, D}^\vee$ of the versal family $\M_{g, D}$ as pulback of $\M_{g, 2 (1 - g) - D}$~, which of course must be possible by the completeness of $\M_{g, 2 (1 - g) - D}$~, find the biggest supersymmetric subfamilies of $\M_{g, 1 - g}$ and show that they are versal amoung all supersymmetric families. \\

If $D \not= 0$ then $T^\rel M_1 \otimes L_{1, D}^* = L_{1, - D}$ induces a strong family isomorphism $\Omega_{1, D}: \M_{1, D}^\vee \mathop{\rightarrow}\limits^\sim \M_{1, - D}$ with $\Omega_{1, D}^\# = \id_{M_1}$~. \\

Now let $g \not= 1$ or $D = 0$~. For all $\eps \in \frac{1}{2} \zz^{2 g}$

\[
\rund{T^\rel M_{g, D}} \otimes L_{g, D}^* \simeq \rund{\id_{M_g}, 2 \eps^\mu \lambda_\mu - \id_{\cz^g}}^* L_{g, 2 (1 - g) - D}
\]

induces a family isomorphism $\rund{\Omega_{g, D, \eps}, \omega_{g, D, \eps}}: \M_{g, D}^\vee \mathop{\rightarrow}\limits^\sim \M_{g, 2 (1 - g) - D}$ with \\
$\Omega_{g, D, \eps}^\# = \rund{\id_{M_g}, 2 \eps^\mu \lambda_\mu - \id_{\cz^g}}$ and

\[
\omega_{g, D, \eps}: V_{g, D} \mathop{\rightarrow}\limits^\sim V_{g, 2 (1 - g) - D} \,~, \, (u, b) \mapsto \rund{u, 2 \eps^\mu \lambda_\mu - b} \,~.
\]

The image $\T_g^\eps$ of the cross section $\rund{\id_{\T_g}, \eps^\mu \lambda_\mu} : \T_g \hookrightarrow V_{g, 1 - g}$ is precisely the fixed point sub\-ma\-ni\-fold of $\omega_{g, 1 - g, \eps}$~, and so $\M_g^\eps := \M_{g, 1 - g}|_{\T_g^\eps}$ is supersymmetric by theorem \ref{charactSUSY}, while the supersymmetry is even uniquely determined up to pullback by strong family automorphisms with identity as body by corollary \ref{SUSYunique0}. Furthermore $\dot\bigcup_{\eps \in \frac{1}{2} \zz^{2 g}} \M_g^\eps$ is the biggest subfamily of $\M_{g, 1 - g}$ admitting a supersymmetry:

\begin{quote}
Let $\X_{u, b}$ admit a supersymmetry~, $(u, b) \in V_{g, 1 - g}$~. Then $\left.L_{g, 1 - g}^{\otimes 2}\right|_{X_u \times \schweif{b}} \simeq T X_u$ or equivalently $\left.L_g^{\otimes 2}\right|_{X_u \times \schweif{b}}$ is trivial, which implies $(u, b)~\in~\frac{1}{2}~\Lambda_u$~.
\end{quote}

We will see that this family is versal amoung all supersymmetric families of genus $g$~, but do we really need all its components for versality if $g \geq 1$~? The answer is no, only two of them: In the case of $\M \twoheadrightarrow V$ supersymmetric $\Phi_{\M, f} \in \Aut_\SUSY \ \M$ iff $f$ maps to $\schweif{\pm 1}$~. So (\ref{exact}) in section \ref{basics} with $\M := \dot\bigcup_{\eps \in \frac{1}{2} \zz^{2 g}} \M_g^\eps$ reduces to an exact sequence

\[ 
1 \hookrightarrow \schweif{\pm 1}^{\frac{1}{2} \zz^{2 g}} \hookrightarrow \Aut_\SUSY \ \mathop{\dot\bigcup}\limits_{\eps \in \frac{1}{2} \zz^{2 g}} \M_g^\eps \mathop{\longrightarrow}\limits^{{}^\#} \Aut \ \mathop{\dot\bigcup}\limits_{\eps \in \frac{1}{2} \zz^{2 g}} M_g \,~.
\]

Since $\dot\bigcup_{\eps \in \frac{1}{2} \zz^{2 g}} \M_g^\eps$ is the biggest subfamily of $\M_{g, 1 - g}$ admitting a supersymmetry, by corollary \ref{SUSYunique0} restriction to $\dot\bigcup_{\eps \in \frac{1}{2} \zz^{2 g}} \M_g^\eps \twoheadrightarrow \dot\bigcup_{\eps \in \frac{1}{2} \zz^{2 g}} \T_g^\eps$ and choice of suitable representatives for each $\eps \in \frac{1}{2} \zz^{2 g}$ seperately induces homomorphisms

\begin{equation} \label{restrict}
\begin{array}{ccc}
\rund{\Aut \ \M_{g, 1 - g}} \left/ \O\rund{V_{g, 1 - g}}^\times \right. & \longrightarrow & \left.\rund{\Aut_{\SUSY} \ \mathop{\dot\bigcup}\limits_{\eps \in \frac{1}{2} \zz^{2 g}} \M_g^\eps}\right/ \schweif{\pm 1}^{\frac{1}{2} \zz^{2 g}} \\
\twoheaddownarrow & \circlearrowleft & \twoheaddownarrow \\
\Gamma_{g, 1 - g} & \mathop{\longrightarrow}\limits_{| \ {}_{\bigcup_{\eps \in \frac{1}{2} \zz^{2 g}} \T_g^\eps}} & \Aut \ \mathop{\bigcup}\limits_{\eps \in \frac{1}{2} \zz^{2 g}} \T_g^\eps
\end{array} \,~,
\end{equation}

the verticle maps being the canonical projections. The first row is an embedding except for $g = 1$~, and the second row gives a faithful action of $\Gamma_{g, 1 - g}$ on $\bigcup_{\eps \in \frac{1}{2} \zz^{2 g}} \T_g^\eps$~, indeed:

\begin{quote}
for $g = 0$ evident by the exact sequence (\ref{exact}) in section~\ref{basics} and since $\T_0^0$ is a single point.

Let $g \geq 1$ and $\gamma \in \Gamma_{g, 1 - g}$ be in the kernel of the second row. Then $\gamma \in \Aut_\strong \ \B_g$ by theorem \ref{supermodexplicit} (i), but since $\aleph_u \eckig{\rund{J|_{X_u}}^* \left.L_g\right|_{X_u \times \schweif{b}} } = - \overline{(u, b)}$ if $g = 1$ or $g = 2$~, by theorem \ref{supermodexplicit} (i) this implies $\gamma = \id_{V_{g, 1 - g}}$~.

Let $g \geq 2$ and $\Phi \in \Aut \ \M_{g, 1 - g}$ be in the kernel of the first row. Then since the second row is an embedding infact $\Phi \in \Aut_\strong \ \M_{g, 1 - g}$~, and therefore \\
$\Phi = \Phi_{\M_{g, 1 - g}, f}$ for some $f \in \O\rund{V_{g, 1 - g}}^\times$ by corollary \ref{autstrong}.

Finally since $V_{1, 0} = \T_1 \times \cz$ is contractible and Stein there exists $f \in \O\rund{V_{1, 0}} \setminus \{0\}$ vanishing on $\bigcup_{\eps \in \frac{1}{2} \zz^2} \T_1^\eps$~. Let $\tau_f \in \Aut_\strong \ M_{1, 0}$ be the translation with such an $f$~. Then $\tau_f^* L_1 \simeq L_1$~, which induces an automorphism \\
$\Phi \in \Aut_\strong \ \M_{1, 0}$ with $\Phi^\# = \tau_f \not= \id_{M_{1, 0}}$~. Therefore \\
$\Phi \notin~\schweif{\Phi_{\M_{1, 0}, h} \,~, \, h \in \O\rund{V_{1, 0}}^\times}$~, but $\rund{\Phi|_{\dot\bigcup_{\eps \in \frac{1}{2} \zz^2} \M_1^\eps}}^\# = \id_{\dot\bigcup_{\eps \in \frac{1}{2} \zz^2} M_1}$~, which shows that $\Phi$ lies in the kernel of the first row for $g = 1$~.
\end{quote}

The second row induces an action of $\left.\rund{\Aut \ M_g} \right/ Q_{g, 1 - g}$ on $\bigcup_{\eps \in \schweif{0, \frac{1}{2}}^{2 g}} \T_g^\eps$ and so on the set $\schweif{0, \frac{1}{2}}^{2 g}$ of connected components of $\bigcup_{\eps \in \schweif{0, \frac{1}{2}}^{2 g}} \T_g^\eps$~. Obviously $\eps \mapsto \left.L^*_{g, 1 - g}\right|_{\left.M_{g, 1 - g}\right|_{\T_g^\eps}}$ gives a $1$-$1$-correspondence between $\schweif{0, \frac{1}{2}}^{2 g}$ and the set of isomorphy classes of relative theta characteristics on $M_g$~, and by this $1$-$1$-correspondence the last action becomes precisely the pullback. Let $\E_g$ be a set of representatives for the orbits. By an easy calculation for $g = 1$ and by the example in \cite{Natanzon} section 12 for $g \geq 2$ one sees that $\abs{\E_g} = 2$ for $g \geq 1$~, and of course $\E_0 = \schweif{0}$~.


\begin{theorem}[Versality of $\dot\bigcup_{\eps \in \E_g} \M_g^\eps$ ] \label{versalfamSUSY}
\item[(i)] {\bf Completeness:} Let $W$ be a contractible Stein manifold and $\N \twoheadrightarrow W$ a supersymmetric family of compact super Riemann surfaces $\Y_w$ of genus $g$~. Then there exists a supersymmetric fiberwise biholomorphic family morphism $(\Xi, \xi): \N \rightarrow \dot\bigcup_{\eps \in \E_g} \M_g^\eps$~. More precisely we have the

{\bf Anchoring property (including local completeness):} For every supersymmetric isomorphism $\sigma: \Y_{w_0} \mathop{\rightarrow}\limits^\sim \X_{v_0}$~, $w_0 \in W$~, $\eps \in \frac{1}{2} \zz^{2 g}$ and $v_0 \in \T_g^\eps$~, there exists $\xi: W \rightarrow \T_g^\eps$ holomorphic such that $\xi\rund{w_0} = v_0$ and $\xi$ can be extended to a supersymmetric fiberwise biholomorphic family morphism $\rund{\Xi, \xi}: \N \rightarrow \M_g^\eps$ with $\Xi|_{\Y_{w_0}} = \sigma$~. Again $\xi$ is uniquely determined by $\sigma^\#$~.

\item[(ii)] {\bf Infinitesimal universality:} $\rund{d \M_{g, 1 - g}}_v T_v \T_g^\eps = H^1\rund{sT^+ \X_v}_0$ for all $\eps \in \frac{1}{2} \zz^{2 g}$ and $v \in \T_g^\eps$~.
\end{theorem}

{\it Proof:} (i) {\it Anchoring property:} {\it Uniqueness} is obvious by the anchoring property of $\M_{g, 1 - g}$~.

{\it Existence:} Let $(\Xi, \xi): \N \rightarrow \M_{g, 1 - g}$ be given by theorem \ref{versalfam} (i) anchoring property. Then by the maximality of $\dot\bigcup_{\eps \in \frac{1}{2} \zz^{2 g}} \M_g^\eps$ and the connectedness of $W$ we see that $\xi(W) \subset \T_g^\eps$~. By corollary \ref{SUSYunique0} there exists $f \in \O(W)^\times$ such that $\Phi_{\N, f} \circ \Xi$ is supersymmetric. Since $\Xi|_{\Y_{w_0}} = \sigma$ is already supersymmetric we have $f\rund{w_0} \in \schweif{\pm 1}$~, and so finally we have to replace $\Xi$ by $\Phi_{\N, \pm f} \circ \Xi$~.

{\it Completeness} is obvious since by construction of $\E_0$ and corollary \ref{SUSYunique0} for every $\eps \in \frac{1}{2} \zz^{2 g}$ there exists a unique $\eps' \in \E_0$ such that $\M_g^\eps \simeq \M_g^{\eps'}$ as supersymmetric families.

(ii) By lemma \ref{SUSYprop} (ii) the image lies in $H^1\rund{sT^+ \X_v}_0$~. Computing the dimension on both sides using lemma \ref{SUSYprop} (i) gives equality. $\Box$ \\

\begin{cor}
The moduli space of all supersymmetric compact super Riemann surfaces of genus $g$ is given by

\[
\left.\rund{\mathop{\bigcup}\limits_{\eps \in \frac{1}{2} \zz^{2 g}} \T_g^\eps} \right/ \Gamma_{g, 1 - g} = \mathop{\dot\bigcup}\limits_{\eps \in \E_g} \rund{\left.\T_g^\eps \right/ \Gamma_g^\eps} \,~,
\]

where $\Gamma_g^\eps$ denotes the normalizer of $\T_g^\eps$ in $\Gamma_{g, 1 - g}$ or equivalently in $\left.\rund{\Aut \ M_g} \right/ Q_{g, 1 - g}$~.
\end{cor}

\section{Super families} \label{superfam}

In this section we construct versal super families $\widetilde \M_{g, D}$ amoung all super families of type $(g, D)$~. Let us start with an oberservation:

\begin{prop} \label{versalitygeneral} Let $\M \twoheadrightarrow V^{|q}$ be a super family of compact complex supermanifolds, $V$ a complex manifold, such that $s\dim H^l\rund{sT \X_v}$~, $\X_v := \pi_\M^{- 1}(v)$~, is independent of $v \in V$ for all $l \in \nz$~. Then equivalent are:

\item[(i)] {\bf infinitesimal universality:} $(d \M)_v: sT_v V^{|q} \rightarrow H^1\rund{sT \X_v}$ is an isomorphism for all $v \in V$~,
\item[(ii)] {\bf a weak anchoring property:} for every super family $\N \twoheadrightarrow \W$ of complex supermanifolds, $\W$ of nilpontency index $\geq 2$ and $\W^\#$ a Stein manifold, and fiberwise biholomorphic super family morphism $(\Sigma, \sigma): \N|_{\W^\natural} \rightarrow \M$ there exists a unique morphism $\xi: \W \rightarrow V^{|q}$ such that $\xi|_{\W^\natural} = \sigma$ and $\xi$ can be extended to a fiberwise biholomorphic super family morphism $(\Xi, \xi): \N \rightarrow \M$ such that $\Xi|_{\N|_{\W^\natural}} = \Sigma$~.
\end{prop}

For the proof we need

\begin{lemma} \label{coherent iso}
Let $\pi: N \twoheadrightarrow W$ be a holomorphic family of compact complex manifolds $Y_w$~, $W$ a complex manifold, $\H$ a coherent sheaf on $W$ and $F \twoheadrightarrow N$ a holomorphic vector bundle. For $r \in \nz$ let

\[
\lambda_{\H, r}: \H \otimes_{\O_W} \pi_{(r)} F \rightarrow \pi_{(r)} \rund{\H \otimes_{\O_W} F} \,~, \, h \otimes \eckig{\rund{S_{i j}} } \mapsto \eckig{\rund{h \otimes S_{i j}} }
\]

be the canonical $\O_W$-module homomorphism, where we write $\H \otimes_{\O_W} F := \rund{\pi^{- 1} \H} \otimes_{\pi^{- 1} \O_W} F$ for short.

\item[(i)] If $H^l\rund{Y_w, F|_{Y_w}}$ is independent of $w \in W$ for all $l \in \nz$ then all $\lambda_{\H, r}$ are isomorphisms.

\item[(ii)] If $\pi_{(l)} F = 0$ for all $l \geq r + 1$ then $\lambda_{\H, r}$ is surjective.
\end{lemma}

{\it Proof:} by induction over the length $n$ of a finite free resolution of $\H$~locally in $W$~. \\


$n = 0$~: Then $\H = 0$~, and the statements are trivial.

$n \rightarrow n + 1$~: Assume $\H$ admits a free resolution

\[
0 \hookrightarrow \O_W^{\oplus s_1} \hookrightarrow \O_W^{\oplus s_2} \rightarrow \dots \rightarrow \O_W^{\oplus s_n} \mathop{\longrightarrow}\limits^\rho \O_W^{\oplus s_{n + 1}} \twoheadrightarrow \H \twoheadrightarrow 0
\]

of length $n + 1$~. Then since $\pi$ as a submersion of complex manifolds is flat, taking tensor product of the exact sequence

\begin{equation} \label{exact coherent}
0 \hookrightarrow \K \hookrightarrow \O_W^{\oplus s_{n + 1}} \twoheadrightarrow \H \twoheadrightarrow 0 \,~,
\end{equation}

$\K := \Im \rho$~, with $F$ gives a short exact sequence $0 \hookrightarrow \K \otimes_{\O_W} F \hookrightarrow F^{\oplus s_{n + 1}} \twoheadrightarrow \H \otimes_{\O_W} F \twoheadrightarrow 0$ of $\O_N$-modules. So for all $r \in \nz$ we obtain

\[
\begin{array}{ccccc}
\phantom{12345} \pi_{(r)} \rund{\K \otimes_{\O_W} F} & \longrightarrow & \phantom{12345} \pi_{(r)} F^{\oplus s_{n + 1}} & \mathop{\longrightarrow}\limits^{(**)} & \pi_{(r)} \rund{\H \otimes_{\O_W} F} \mathop{\longrightarrow}\limits^{(*)} \pi_{(r + 1)} \rund{\K \otimes_{\O_W} F} \\
\lambda_{\K, r} \uparrow & \circlearrowleft & \lambda_{\O_W^{\oplus \rund{s_n + 1}}, r} \uparrow \phantom{12345} & \circlearrowleft & \uparrow \lambda_{\H, r} \phantom{12345678901234} \\
0 \hookrightarrow \phantom{1} \K \otimes_{\O_W} \pi_{(r)} F & \longrightarrow & \phantom{12345} \rund{\pi_{(r)} F}^{\oplus s_{n + 1}} & \twoheadrightarrow & \H \otimes_{\O_W} \pi_{(r)} F \phantom{1,} \twoheadrightarrow \phantom{1} 0 \phantom{1234567890123}
\end{array}
\]

and

\[
\begin{array}{ccc}
\pi_{(r)} \rund{\H \otimes_{\O_W} F} \mathop{\longrightarrow}\limits^{(*)} \pi_{(r + 1)} \rund{\K \otimes_{\O_W} F} & \mathop{\longrightarrow}\limits^{(***)} & \pi_{(r + 1)} F^{\oplus s_{n + 1}} \phantom{123456789} \\
\phantom{1234567890} \lambda_{\K, r + 1} \uparrow & \circlearrowleft & \uparrow \lambda_{\O_W^{\oplus \rund{s_n + 1}}, r + 1} \\
\phantom{12345678901} 0 \phantom{1} \hookrightarrow \phantom{1,} \K \otimes_{\O_W} \pi_{(r + 1)} F & \longrightarrow & \rund{\pi_{(r + 1)} F}^{\oplus s_{n + 1}} \phantom{123456789}
\end{array}
\]

with exact first rows. The second rows are given by the tensor product of (\ref{exact coherent}) with $\pi_{(r)} F$ resp. $\pi_{(r + 1)} F$~. Trivially $\lambda_{\O_W^{\oplus \rund{s_n + 1}}, r}$ and $\lambda_{\O_W^{\oplus \rund{s_n + 1}}, r + 1}$ are isomorphisms, and

\[
0 \hookrightarrow \O_W^{\oplus s_1} \hookrightarrow \O_W^{\oplus s_2} \rightarrow \dots \rightarrow \O_W^{\oplus s_n} \mathop{\twoheadrightarrow}\limits^\rho \K \twoheadrightarrow 0
\]

is a free resolution of $\K$ of length $n$~. 

(i) Since $\pi_{(r)} F$ and $\pi_{(r + 1)} F$ are locally free by \cite{GrauRem} theorem 10.5.5 also the second rows become exact. Furthermore by induction hypothesis also $\lambda_{\K, r}$ and $\lambda_{\K, r + 1}$ are isomorphisms. We see that $(***)$ is injective and so $(*)$ is the zero map. Therefore by the five lemma $\lambda_{\H, r}$ is an isomorphism.

(ii) By induction hypothesis $\lambda_{\K, r + 1}$ is surjective, and so $\pi_{(r + 1)} \rund{\K \otimes_{\O_W} F} = 0$~. We see that $(**)$ is surjective, and therefore so must be $\lambda_{\H, r}$~. $\Box$ \\

{\it Proof of proposition \ref{versalitygeneral}:} Let $\M$ be given by local super charts $\left.\rund{V^{|q} \times \cz^{m|n}}\right|_{U_i}$~, $U_i \subset V \times \cz^m$ open, with transition morphisms

\[
\Phi_{i j}: \left.\rund{V^{|q} \times \cz^{m|n}}\right|_{U_{i j}} \mathop{\rightarrow}\limits^\sim \left.\rund{V^{|q} \times \cz^{m|n}}\right|_{U_{j i}} \,~.
\]

(i) $\Rightarrow$ (ii): Write $\W = (W, \S)$ and let $\m \lhd \S$ be the nilradical. Using \cite{GrauRem} theorem 10.5.5, since $\sdim H^1\rund{sT \X_v}$ is independent of $v \in V$~, we see that

\[
\rund{\sigma^\#}^* (d \M): \rund{\sigma^\#}^* sT V^{|q} \rightarrow \rund{\pi_{\sigma^* \M}^\#}_{(1)} sT^\rel \rund{\sigma^\#}^* \M \simeq \rund{\pi_\N^\#}_{(1)} sT^\rel \N|_W
\]

is an isomorphism, the last identification by the strong isomorphism $\N|_W \simeq \rund{\sigma^\#}^* \M$ associated to $\Sigma|_{\N|_W}$~. Let $\N$ be given by local super charts $\left.\rund{\W \times \cz^{m|n}}\right|_{\Omega_r}$~, $\Omega_r \subset W \times \cz$ open, with transition morphisms 

\[
\Psi_{r s}: \left.\rund{\W \times \cz^{m|n}}\right|_{\Omega_{r s}} \mathop{\rightarrow}\limits^\sim \left.\rund{\W \times \cz^{m|n}}\right|_{\Omega_{s r}} \,~.
\]

After refining the atlas of $\N$ we may assume that for every $r$ there exists a local super chart $i(r)$ of $\M$ covering $\Sigma^\# \rund{\Omega_r} \subset \M^\#$~. Therefore if we write $\Sigma$ in the local super charts $r$ of $\N$ and $i(r)$ of $\M$ as $\rund{\sigma, \Sigma_r}: \left.\rund{\W^\natural \times \cz^{m|n}}\right|_{\Omega_r} \rightarrow \left.\rund{V^{|q} \times \cz^{m|n}}\right|_{U_{i(r)} }$ we obtain

\[
\begin{array}{ccc}
\phantom{123456789012} \left.\rund{\W^\natural \times \cz^{m|n}}\right|_{\Omega_{r s}} & \mathop{\longrightarrow}\limits^{\rund{\sigma, \Sigma_r}} & \left.\rund{V^{|q} \times \cz^{m|n}}\right|_{U_{i(r), i(s)} } \phantom{1234567890} \\
\left.\Psi_{r s}\right|_{\W^\natural \times \cz^{m|n}} \downarrow & \circlearrowleft & \downarrow \Phi_{i(r), i(s)} \phantom{1234}~. \\
\phantom{123456789012} \left.\rund{\W^\natural \times \cz^{m|n}}\right|_{\Omega_{s r}} & \mathop{\longrightarrow}\limits_{\rund{\sigma, \Sigma_s}} & \left.\rund{V^{|q} \times \cz^{m|n}}\right|_{U_{i(s), i(r)} } \phantom{1234567890}
\end{array}
\]

{\it Uniqueness:} Assume both $\rund{\Xi, \xi}$ and $\rund{\widehat\Xi, \widehat\xi}$ are fiberwise biholomorphic super family morphisms from $\N$ to $\M$ such that $\xi|_{\W^\natural} = \sigma = \left.\widehat\xi\right|_{\W^\natural}$ and $\Xi|_{\N|_{\W^\natural}} = \Sigma =~\left.\widehat\Xi\right|_{\N|_{\W^\natural}}$~. Then we can write $\widehat\xi = \xi - \delta$~, $\delta \in H^0\rund{W, \m^{k - 1} \boxtimes_{\O_W} \rund{\sigma^\#}^* sT V^{|q}}_0$ suitable, and express $\Xi$ and $\widehat\Xi$ in the local super charts $r$ of $\N$ and $i(r)$ of $\M$ as

\[
\rund{\xi, \Xi_r} \, \text{ resp. } \, \rund{\xi - \delta, \Xi_r \circ \rund{\Id_{\W \times \cz^{m|n}} + \alpha_r}}: \left.\rund{\W \times \cz^{m|n}}\right|_{\Omega_r} \rightarrow \left.\rund{V^{|q} \times \cz^{m|n}}\right|_{U_{i(r)} } \,~,
\]

$\left.\Xi_r\right|_{\W^\natural \times \cz^{m|n}} = \Sigma_r$~, with $\alpha_r$ suitable even sections of $\m^{k - 1} \boxtimes_{\O_W} sT^\rel \N|_W$~, each $\alpha_r$ expressed in the local super chart $r$ of $\N$~. We obtain

\[
\begin{array}{ccc}
\phantom{123} \left.\rund{\W \times \cz^{m|n}}\right|_{\Omega_{r s}} & \mathop{\longrightarrow}\limits^{\rund{\xi, \Xi_r}} & \left.\rund{V^{|q} \times \cz^{m|n}}\right|_{U_{i(r), i(s)} } \phantom{12345} \\
\Psi_{r s} \downarrow & \circlearrowleft & \downarrow \Phi_{i(r), i(s)} \\
\phantom{123} \left.\rund{\W \times \cz^{m|n}}\right|_{\Omega_{s r}} & \mathop{\longrightarrow}\limits_{\rund{\xi, \Xi_s}} & \left.\rund{V^{|q} \times \cz^{m|n}}\right|_{U_{i(s), i(r)} } \phantom{12345}
\end{array} \,~.
\]

and the same with $\xi$ replaced by $\xi - \delta$~, $\Xi_r$ by $\Xi_r \circ \rund{\Id_{\W \times \cz^{m|n}} + \alpha_r}$ and $\Xi_s$ by \\
$\Xi_s \circ \rund{\Id_{\W \times \cz^{m|n}} + \alpha_s}$~. Therefore $\rund{\rund{\sigma^\#}^* (d \M)} \delta = \eckig{d \rund{\alpha_r}} = 0$ as even sections of

\[
\rund{\pi_{\N}^\#}_{(1)} \rund{\m^{k - 1} \boxtimes_{\O_W} sT^\rel \N|_W} = \m^{k - 1} \boxtimes_{\O_W} \rund{\pi_{\N}^\#}_{(1)} \rund{sT^\rel \N|_W}
\]

identified via the canonical graded isomorphism $\lambda_{\m^{k - 1}, 1}$ of lemma \ref{coherent iso} (i) with $\H := \m^{k - 1}$~, $F := sT^\rel \N|_W$~, and $\pi := \pi_\N^\#$~. But $\rund{\sigma^\#}^* (d \M)$ is an isomorphism, and so $\delta = 0$~. \\

{\it Existence: Local construction:} Locally in $\W$ we can take arbitrary extensions $\xi: \W \rightarrow~V^{|q}$ of $\sigma$ and $\rund{\xi, \widehat \Sigma_r}: \left.\rund{\W \times \cz^{m|n}}\right|_{\Omega_r} \rightarrow \left.\rund{V^{|q} \times \cz^{m|n}}\right|_{U_{i(r)} }$ of $\rund{\sigma, \Sigma_r}$~. Then

\[
\begin{array}{ccc}
\phantom{1234567890123456789} \left.\rund{\W \times \cz^{m|n}}\right|_{\Omega_{r s}} & \mathop{\longrightarrow}\limits^{\rund{\xi, \widehat \Sigma_r}} & \left.\rund{V^{|q} \times \cz^{m|n}}\right|_{U_{i(r), i(s)} } \phantom{12345678} \\
\Psi_{r s} \circ \rund{\Id_{\W \times \cz^{m|n}} - \eps_{r s}} \downarrow & \circlearrowleft & \downarrow \Phi_{i(r), i(s)} \phantom{1234} \\
\phantom{1234567890123456789} \left.\rund{\W \times \cz^{m|n}}\right|_{\Omega_{s r}} & \mathop{\longrightarrow}\limits_{\rund{\xi, \widehat \Sigma_s}} & \left.\rund{V^{|q} \times \cz^{m|n}}\right|_{U_{i(s), i(r)} } \phantom{12345678}
\end{array}
\]

with suitable even sections $\eps_{r s}$ of $\m^{k - 1} \boxtimes_{\O_W} sT^\rel \N|_W$~, each $\eps_{r s}$ expressed in the local super chart $r$ of $\N$~. A straight forward calculation shows that $\rund{\eps_{r s}}$ is a $1$-cocycle. Since $\rund{\sigma^\#}^* (d \M)$ is an isomorphism there exists, locally in $W$~, an even section $\delta$ of $\m^{k - 1} \boxtimes_{\O_W}~\rund{\sigma^\#}^*~sT V^{|q}$ such that $\rund{\rund{\sigma^\#}^* (d \M)} \delta~=~\eckig{\rund{\eps_{r s}} }$ as even sections of

\[
\rund{\pi_{\N}^\#}_{(1)} \rund{\m^{k - 1} \boxtimes_{\O_W} sT^\rel \N|_W} = \m^{k - 1} \boxtimes_{\O_W} \rund{\pi_{\N}^\#}_{(1)} \rund{sT^\rel \N|_W} \,~.
\]

So after passing from $\xi$ to $\xi + \delta$ we may assume that $\rund{\eps_{r s}}$ is a coboundary. Therefore after refining the atlas of $\N$ there exist even sections $\alpha_r$ of $\m^{k - 1} \boxtimes_{\O_W} sT^\rel \N|_W$~, each $\alpha_r$ on the range of the local super chart $r$ of $\N$~, such that $\eps_{r s} = \alpha_r - \alpha_s$~. So by expressing each $\alpha_r$ in the local super chart $r$ we obtain

\[
\begin{array}{ccc}
\phantom{1234} \left.\rund{\W \times \cz^{m|n}}\right|_{\Omega_r} & \mathop{\longrightarrow}\limits^{\rund{\xi, \widehat\Sigma_r} \circ \rund{\Id_{\W \times \cz^{m|n}} +  \alpha_r}} & \left.\rund{V^{|q} \times \cz^{m|n}}\right|_{U_{i(r)} } \phantom{1234} \\
\Psi_{r s} \downarrow & \circlearrowleft & \downarrow \Phi_{i(r), i(s)} \\
\phantom{1234} \left.\rund{\W \times \cz^{m|n}}\right|_{\Omega_s} & \mathop{\longrightarrow}\limits_{\rund{\xi, \widehat\Sigma_s} \circ \rund{\Id_{\W \times \cz^{m|n}} +  \alpha_s}} & \left.\rund{V^{|q} \times \cz^{m|n}}\right|_{U_{i(s)} } \phantom{1234}
\end{array} \,~,
\]

and so all $\rund{\xi, \widehat\Sigma_r} \circ \rund{\Id_{\W \times \cz^{m|n}} +  \alpha_r}$ glue together to a fiberwise biholomorphic super family morphism $(\Xi, \xi): \N \rightarrow \M$~. \\

{\it Global construction:} By the uniqueness $\xi: \W \rightarrow V^{|q}$ is already globally defined on $\W$~. So on their overlaps the locally constructed fiberwise biholomorphic super family morphisms $(\Xi, \xi)$ from $\N$ to $\M$ differ by sections of

\[
\rund{\pi_\N^\#}_* \rund{\m^{k - 1} \boxtimes_{\O_W} sT^\rel \N|_W}_0 = \rund{\m^{k - 1} \boxtimes_{\O_W} \rund{\pi_\N^\#}_* sT^\rel \N|_W}_0
\]

identified via the canonical graded isomorphism $\lambda_{\m^{k - 1}, 0}$ of lemma \ref{coherent iso} (i), which is coherent by the direct image theorem 10.4.6 of \cite{GrauRem} since $\pi_{\N^\#}$ is proper. We see that the obstructions to define $\Xi$ globally lie in $H^1\rund{W, \m^{k - 1} \boxtimes_{\O_W} \rund{\pi_\N^\#}_* sT^\rel \N|_W}_0$~, which vanishes by Cartan's theorem B since $W$ is Stein. \\




(ii) $\Rightarrow$ (i): Let $v \in V$~. Then $\X_v$ is given by local super charts $\Omega_i^{|n}$~, $\Omega_i := U_i \cap \rund{\schweif{v} \times \cz^m} \subset \cz^m$ open, with transition morphisms $\Psi_{i j}:= \Phi_{i j}(v, \diamondsuit): \Omega_{i j}^{|n} \mathop{\rightarrow}\limits^\sim \Omega_{j i}^{|n}$~. \\

{\it Surjectivity:} Let $\beta \in H^1\rund{sT \X_v}_1$~. After refining the atlas of $\M$ we may assume that $\beta$ is represented by the $1$-cocycle $\rund{\eps_{i j}}$~, each $\eps_{i j}$ an odd section of $sT \X_v$ on the overlap of the local super charts $i$ and $j$~. So

\[
\rund{\Id_{\cz^{0|1}}, \Psi_{i j}} \circ \rund{\Id_{\cz^{0|1} \times \Omega_{i j}^{|n}} + \eta \eps_{i j}} : \cz^{0|1} \times \Omega_{i j}^{|n} \mathop{\rightarrow}\limits^\sim \cz^{0|1} \times \Omega_{j i}^{|n} \,~,
\]

$\eps_{i j}$ expressed in the local super chart $i$~, are the transition morphisms of a super family $\N \twoheadrightarrow (\schweif{0}, \bigwedge \cz) = \cz^{0|1}$~, $\eta$ being the odd super coordinate on $\cz^{0|1}$~, with $\N|_{\schweif{0}} = \X_v$ and $(d \N)_0 \partial_\eta = \beta$~. By (ii) there exists a fiberwise biholomorphic super family morphism $\rund{\Xi, \xi}: \N \rightarrow \M$ such that $\xi(0) = v$ and $\Xi|_{\X_v} = \Id_{\X_v}$~. Therefore $\xi = v + \eta \delta$ with some $\delta \in \rund{sT_v V^{|q}}_1$~, and so $(d \M)_v \delta = (d \N)_0 \partial_\eta = \beta$~. Same calculation for $\beta \in H^1\rund{sT \X_v}_0$ regarding $\bigwedge \cz$ as purely even. \\

{\it Injectivity:} Let $(d \M)_v \delta = 0$ for some $\delta \in \rund{sT_v V^{|q}}_1$~. Define $\xi:= v + \eta \delta: \cz^{0|1} \rightarrow V^{|q}$~. Then $\xi^* \M$ is given by local super charts $\cz^{0|1} \times \Omega_i^{|n}$ with transition morphisms

\[
\rund{\Id_{\cz^{0|1}}, \Psi_{i j}} \circ \rund{\Id_{\cz^{0|1} \times \Omega_i^{|n}} + \eta \eps_{i j}}: \cz^{0|1} \times \Omega_{i j}^{|n} \mathop{\rightarrow}\limits^\sim \cz^{0|1} \times \Omega_{j i}^{|n} \,~,
\]

where $\eps_{i j} := \rund{d \Psi_{i j}}^{- 1} \Pr \rund{\rund{d \Phi_{i j}}(v, \diamondsuit)} \delta \in H^0\rund{sT \ \Omega_{i j}^{|n}}_1$ and \\
$\Pr: \left.\rund{sT \rund{V^{|q} \times \Omega_{i j}^{|n}} }\right|_{\{v\} \times \Omega_{i j}^{|n}} \twoheadrightarrow sT \ \Omega_{i j}^{|n}$ denotes the canonical projection.

$\rund{\breve \xi_\M, \xi}: \xi^* \M \rightarrow \M$ is a fiberwise biholomorphic super family morphism with $\left.\breve \xi_\M\right|_{\X_v} = \Id_{\X_v}$~.

On the other hand $0 = (d \M)_v \delta = \eckig{\rund{\eps_{i j}}}$ in $H^1\rund{sT \X_v}_1$~, which implies that after a refinement of the atlas of $\M$ there exist $\alpha_i \in H^0\rund{sT \ \Omega_i^{|n}}_1$ such that $\eps_{i j} = \alpha_i - \rund{\rund{d \Psi_{i j}}^{- 1} \alpha_j} \circ \Psi_{i j}$~. Therefore $\Pr_{\Omega_i^{|n}} + \eta \alpha_i: \cz^{0|1} \times \Omega_i^{|n} \rightarrow \Omega_i^{|n}$ glue together to another fiberwise biholomorphic super family morphism $(\Phi, v): \xi^* \M \rightarrow \M$ with $\Phi|_{\X_v} = \Id_{\X_v}$~, and so $\delta = 0$ by the uniqueness in (ii). Again same calculation for $\beta \in H^1\rund{sT \X_v}_0$ regarding $\bigwedge \cz$ as purely even. $\Box$ \\

So let us attack the construction of a versal super family by constructing a super family $\widetilde \M_{g, D} \twoheadrightarrow V_{g, D}^{|q}$ of compact super Riemann surface of type $(g, D)$ such that $\left.\widetilde \M_{g, D}\right|_{V_{g, D}} = \M_{g, D}$ and

\[
\rund{d \widetilde \M_{g, D}}_v: sT_v V_{g, D}^{|q} \rightarrow H^1\rund{sT \X_v}
\]

is an isomorphism for all $v \in V$~. Obviously a necessary condition is that $q = h^1\rund{sT \X_v}_1$ is independent of $v \in V_{g, D}$~. As $\O_{M_{g, D}}$-module $\rund{sT^\rel \M_{g, D}}_1$ splits as a direct sum

\[
\rund{sT^\rel \M_{g, D}}_1 = \rund{\rund{T^\rel M_{g, D}} \otimes L_{g, D}^*} \oplus L_{g, D} \,~.
\]

If $g \geq 1$ and $D \in \schweif{4 - 4 g, \dots, 2 - 2 g} \cup \schweif{0, \dots, 2 g - 2}$ then $h^1\rund{T X_v \otimes \rund{\left.L_{g, D}^*\right|_{X_v}} \oplus \left.L_{g, D}\right|_{X_v}}$ is {\bf not} independent of $v \in \T_g$ because of the presence of special divisors. These types will be treated in section \ref{remain}.

{\bf So for the rest of this section let $g = 0$ or $D \notin \schweif{4 - 4 g, \dots, 2 - 2 g} \cup \schweif{0, \dots, 2 g - 2}$~.} Then $r_D := h^1\rund{T X_v \otimes \rund{\left.L_{g, D}^*\right|_{X_v}} }$ and $s_D := h^1\rund{X_v, \left.L_{g, D}\right|_{X_v}}$ are independent of $v \in V_{g, D}$~, and by the Riemann-Roch theorem and Serre duality

\[
s_D = \left\{\begin{array}{l}
0 \text{ if } D \geq 2 g - 1 \\
g - 1 - D \text{ if } D \leq - 1
\end{array}\right.~.
\]

$T^\rel M_g \otimes L_{g, D}^* = \rund{\Omega_{g, D, 0}^\#}^* L_{g, 2 (1 - g) - D}$ if $g \not= 1$ and $T^\rel M_1 \otimes L_{1, D}^* = L_{1, - D}$ if $D \not= 0$~, and so in both cases $r_D = s_{2 (1 - g) - D}$~.

\begin{figure}[H]
\begin{center}
\scalebox{0.18}{\includegraphics{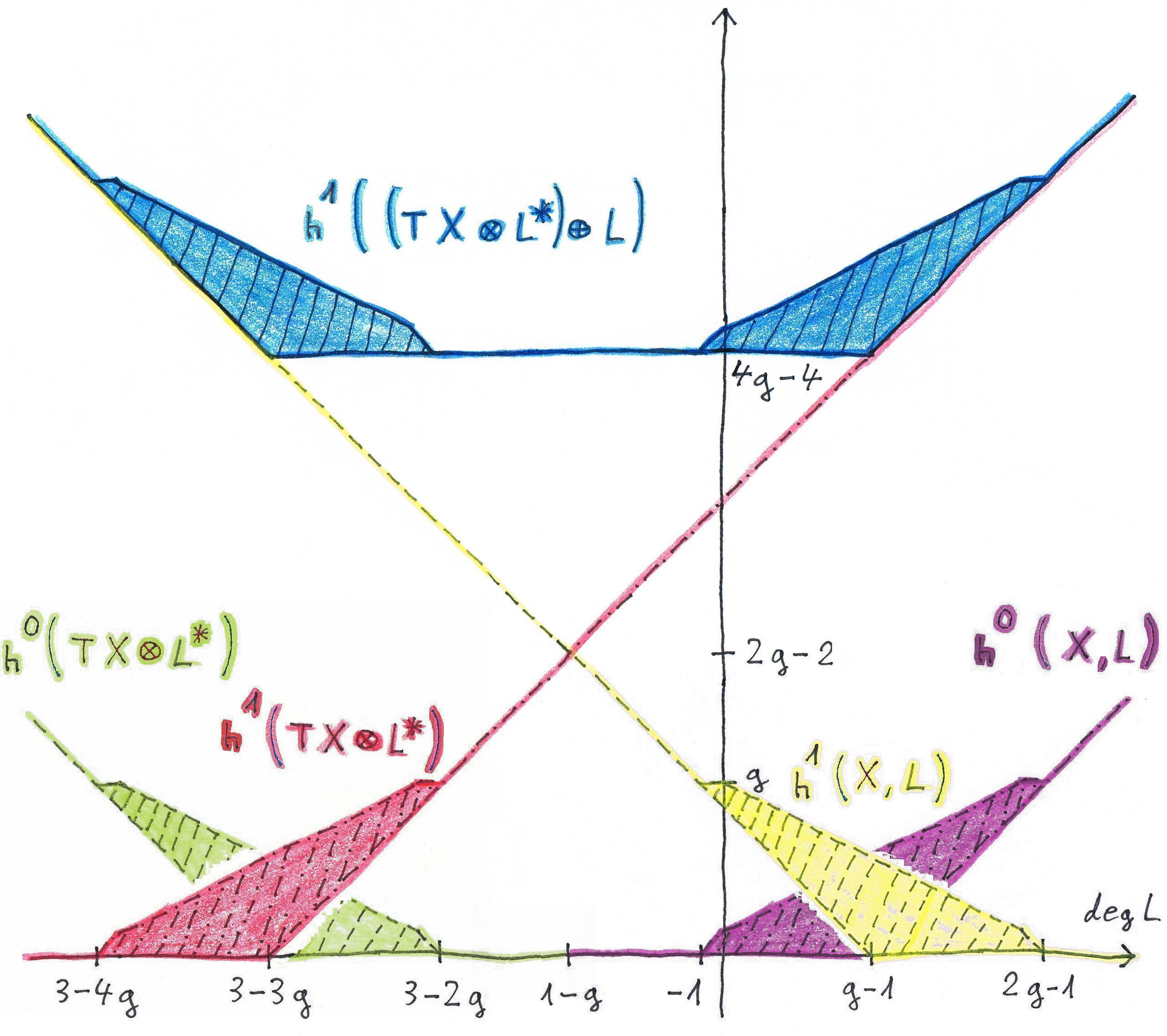}}
\caption{Possible values of $h^0$ and $h^1$ for a holomorphic line bundle $L$ on a compact Riemann surface $X$ of genus $g$~.}
\end{center}
\end{figure}

\begin{defin}
$s\T_{g, D} := V_{g, D}^{| r_D + s_D}$ is called the super Teichmüller space for type $(g, D)$~.
\end{defin}

Observe that $s\T_{0, D}$ is a single point iff $D \in \schweif{- 1, \dots, 3}$~.

\begin{theorem}[Construction of $\widetilde \M_{g, D}$ ] \label{constructsuperfam}
There exist super families $\widetilde \M_{g, D} \twoheadrightarrow s\T_{g, D}$ with $\left.\widetilde \M_{g, D}\right|_{V_{g, D}} = \M_{g, D}$ and local super coordinates $\rund{v, z | \eta, \vartheta, \zeta}$ (which we call the standard local super coordinates of $\widetilde \M_{g, D}$~), where $v$ denote the even and $\eta_1, \dots, \eta_{r_D}$ and $\vartheta_1, \dots, \vartheta_{s_D}$ the odd super coordinates on $s\T_{g, D}$~, such that

\begin{itemize}
\item[(i)] {\bf Infinitesimal universality:} $\rund{d \widetilde \M_{g, D}}_v: sT_v \ s\T_{g, D} \rightarrow H^1\rund{sT \X_v}$ is an isomorphism for every $v \in V_{g, D}$~,

\item[(ii)] every $\widetilde \M_{g, D} \twoheadrightarrow s\T_{g, D}$ admits an action of the trivial bundle $s\T_{g, D} \times \cz^\times \twoheadrightarrow s\T_{g, D}$ in the standard local super coordinates given by $\rund{\id_{V_{g, D}}, z \left| \frac{1}{a} \eta, a \vartheta, a \zeta\right.}$~, $a$ denoting the coordinate in $\cz^\times$-direction, and

\item[(iii)] there exist super family isomorphisms

\[
\rund{\widetilde \Omega_{g, D}, \widetilde \omega_{g, D}}: {\widetilde \M_{g, D}}^\vee \mathop{\rightarrow}\limits^\sim \widetilde \M_{g, 2 (g - 1) - D}
\]

in the standard local super coordinates given by $(u, - b, z | \vartheta, \eta, \zeta)$~, $u$ denoting the coordinates on $\T_g$ and $b$ on $\cz^g$~, having $\left.\widetilde \Omega_{g, D}\right|_{\M_{g, D}} = \Omega_{g, D, 0}$ if $g \not= 1$ or $D = 0$ resp. in the standard local super coordinates given by $\rund{\id_{V_{g, D}}, z | \vartheta, \eta, \zeta}$ and $\left.\widetilde \Omega_{1, D}\right|_{\M_{1, D}} = \Omega_{1, D}$ if $g = 1$ and $D \not= 0$ (in other words, $\widetilde \Omega_{g, D}$ interchanges the $\eta$- and $\vartheta$-coordinates and flips the signs in $\cz^g$-direction).
\end{itemize}

\end{theorem}

{\it Proof:} For $n \in \nz \setminus \{0\}$ we define the ringed spaces

\[
\V_n^D := \rund{V_{g, D}, \O_{V_{g, D}} \otimes \bigwedge \cz^{r_D + s_D} \left/ \rund{\bigwedge\nolimits^n \cz^{r_D + s_D}}\right.} = \rund{\schweif{0}, \bigwedge \cz^{r_D + s_D} \left/ \rund{\bigwedge\nolimits^n \cz^{r_D + s_D}}\right.} \times V_{g, D}  \,~.
\]

Then $\V_1^D = V_{g, D}$~, $\V_{r_D + s_D + 1}^D = s\T_{g, D}$~, and $\rund{\V_{n + 1}^D}^\natural = \V_n^D$ for all $n = 1, \dots, r_D + s_D$~. We construct super families $\M_n^D \twoheadrightarrow \V_n^D$ fulfilling (i) for $n \geq 2$~, (ii) and (iii) and $\left.\M_n^D\right|_{\V_{n- 1}^D} = \M_{n - 1}^D$ by induction on $n = 1, \dots, r_D + s_D + 1$~. Then $\widetilde \M_{g, D} := \M_{r_D + s_D + 1}^D$ will have all the desired properties. \\

$\M_1^D := \M_{g, D}$~. Choose atlasses of $\M_1^D$ with local super charts $U_{i, D}^{|1}$~, $U_{i, D} \subset V_{g, D} \times \cz$ open, and transition morphisms $\Phi^D_{i j} = \rund{\varphi_{i j}^D \left| A^D_{i j} \zeta\right.}: U_{i j, D}^{|1} \mathop{\rightarrow}\limits^\sim U_{j i, D}^{|1}$ with suitable $\varphi_{i j}^D: U_{i j, D} \mathop{\rightarrow}\limits^\sim U_{j i, D}$ and $A^D_{i j} \in \O\rund{U_{i j, D}}^\times$~.

Define the family isomorphism $\rund{\Omega_D, \omega_D}: \rund{\M_1^D}^\vee \mathop{\rightarrow}\limits^\sim \M_1^{2 (1 - g) - D}$ as \\
$\rund{\Omega_D, \omega_D}~:=~\rund{\Omega_{g, D, 0}, \omega_{g, D, 0}}$ if $g \not= 1$ resp. $\rund{\Omega_D, \omega_D} := \rund{\Omega_{1, D}, \id_{\T_1}}$ if $g = 1$~. We may assume a bijection $i \mapsto \widehat i$ between the local super charts of $\M_1^D$ and the ones of $\M_1^{2 (1 - g) - D}$ such that $\Omega_D^\#: M_{g, D} \rightarrow M_{g, 2 (g - 1) - D}$ maps the range of $i$ to the one of $\widehat i$~, $\Omega_D$ in the local super charts $i$ and $\widehat i$ is given by $\rund{\omega_D, \Id_{\cz^{|1}} }: U_{i, D}^{|1} \rightarrow U_{\widehat i, 2 (1 - g) - D}^{|1}$~, and $\rund{\Phi^D_{i j}}^\vee = \omega_D^* \Phi^{2 (1 - g) - D}_{\widehat i \ \widehat j}$~. \\

Since $V_{g, D}$ is contractible and Stein, by \cite{GrauRem} theorem 10.5.5 and Grauert's theorem $\rund{\pi_{M_{g, D}} }_{(1)} L_{g, D}$ is globally free of rank $s_D$~. So let $\rund{\gamma^{D, 1}, \dots, \gamma^{D, s_D}}$ be a global frame. After refining the atlas of $\M_{g, D}$~, $\gamma^{D, \sigma}$ are given by cocycles $\rund{c_{i j}^{D, \sigma} \partial_\zeta}$~, $c_{i j}^{D, \sigma} \in \O\rund{U_{i j, D}}$~:

\begin{quote}
Of course $\gamma^{D, \sigma}$ has such a realization locally in $V_{g, D}$~, and on their overlaps two local realizations differ by a section of the sheaf $\rund{\pi_{M_{g, D}}}_* L_{g, D}$ on $V_{g, D}$~, which is locally free by \cite{GrauRem} theorem 10.5.5. Therefore the obstructions to find a global realization lie in $H^1\rund{V_{g, D}, \rund{\pi_{M_{g, D}}}_* L_{g, D}}$~, which is $0$ by Cartan's theorem B.
\end{quote}


Now $\rund{\beta^{D, 1}, \dots, \beta^{D, r_D}}$~, $\beta^{D, \rho} := \gamma^{2 (1 - g) - D, \rho} \circ \Omega_D^\#$~, is a global frame of $\rund{\pi_{M_{g, D}} }_{(1)} T^\rel M_{g, D} \otimes~L_{g, D}^*$~, and each $\beta^{D, \rho}$ is realized by $\rund{b_{i j} \zeta \partial_z}$~, \\ $b_{i j}^{D, \rho} := c_{\widehat i \ \widehat j}^{2 (1 - g) - D, \rho} \circ \rund{\omega_D, \id_\cz}$~.

We define $\M_2^D$ by the local super charts $\left.\rund{\V_2^D \times \cz^{|1}}\right|_{U_{i, D}}$ with transition morphisms

\begin{eqnarray*}
\Phi_{i j}^{D, 2} &:=& \rund{\Id_{\rund{\schweif{0}, \bigwedge \cz^{r_D + s_D} \left/ \rund{\bigwedge\nolimits^2 \cz^{r_D + s_D}}\right.}}, \Phi_{i j}} \circ \rund{\Id_{\left.\rund{\V_2^D \times \cz^{|1}}\right|_{U_{i j, D}} } + \eta_\rho b_{i j}^{D, \rho} \zeta \partial_z + \vartheta_\sigma c_{i j}^{D, \sigma} \partial_\zeta} \\
&=& \rund{\left.\varphi_{i j} \circ \rund{\id_{V_{g, D}}, z + \eta_\rho b_{i j}^{D, \rho} \zeta} \, \right| \, \eta, \vartheta, A_{i j} \rund{\zeta + \vartheta_\sigma c_{i j}^{D, \sigma}} }: \\
&& \phantom{12} \left.\rund{\V_2^D \times \cz^{|1}}\right|_{U_{i j, D}} \mathop{\rightarrow}\limits^\sim \left.\rund{\V_2^D \times \cz^{|1}}\right|_{U_{j i, D}} \,~,
\end{eqnarray*}

which are obviously $\rund{\V_2^D \times \cz^\times}$-equivariant and fulfill $\rund{\Phi_{i j}^{D, 1}}^\vee = \rund{\omega_D | \vartheta, \eta}^* \Phi_{\widehat i \ \widehat j}^{2 (1 - g) - D, 1}$~. Furthermore $\rund{d \M_2^D}_v \partial_{\eta_\rho} = \beta^{D, \rho}(v)$ and $\rund{d \M_2^D}_v \partial_{\vartheta_\sigma} = \gamma^{D, \sigma}(v)$ in $H^1\rund{sT \X_v}_1$ for all $v \in V_{g, D}$~, and so we already have (i). \\

$n \rightarrow n + 1$~, $n \geq 2$~: Let $\M_n^D$ be given by local super charts with $\rund{V_n^D \times \cz^\times}$-equivariant transition morphisms

\[
\Phi_{i j}^{D, n}: \left.\rund{\V_n^D \times \cz^{|1}}\right|_{U_{i j, D}} \mathop{\rightarrow}\limits^\sim \left.\rund{\V_n^D \times \cz^{|1}}\right|_{U_{j i, D}}
\]

having $\rund{\Phi_{i j}^{D, n}}^\vee = \rund{\omega_D | \vartheta, \eta}^* \Phi_{\widehat i \ \widehat j}^{2 (1 - g) - D, n}$~. Then the $z$- and $\zeta$-component of $\Phi_{i j}^{D, n}$ are an even resp. odd element of $\O\rund{U_{i j, D}} \otimes \rund{\bigwedge \cz^{\left|r_D + s_D\right.} \left/ \rund{\bigwedge\nolimits^n \cz^{r_D + s_D}} \right.} \boxtimes \bigwedge \cz$~, the first invariant, the second $1$-homogeneous under the $\rund{\V_n^D \times \cz^\times}$-action.

{\it Local construction:} Taking preimages of the same homogeneity type of these components under the canonical projection

\[
\O_{V_{g, D}} \otimes \rund{\bigwedge \cz^{r_D + s_D} \left/ \rund{\bigwedge\nolimits^{n + 1} \cz^{r_D + s_D}}\right.} \boxtimes \bigwedge \cz \twoheadrightarrow \O_{V_{g, D}} \otimes \rund{\bigwedge \cz^{r_D + s_D} \left/ \rund{\bigwedge\nolimits^n \cz^{r_D + s_D}}\right.} \boxtimes \bigwedge \cz
\]

gives a strong $\rund{\V_{n + 1}^D \times \cz^\times}$-equivariant isomorphism

\[
\Psi_{i j}: \left.\rund{\V_{n + 1}^D \times \cz^{|1}}\right|_{U_{i j, D}} \mathop{\rightarrow}\limits^\sim \left.\rund{\V_{n + 1}^D \times \cz^{|1}}\right|_{U_{j i, D}} \\
\]

such that $\Psi_{i j} |_{\left.\rund{\V_n^D \times \cz^{|1}}\right|_{U_{i j, D}}} = \Phi_{i j}^{D, n}$~. Therefore

\[
\Psi_{k i} \circ \Psi_{j k} \circ \Psi_{i j} = \Id_{\left.\rund{\V_{n + 1}^D \times \cz^{|1}}\right|_{U_{i j k, D}} } + \alpha_{i j k}
\]

with an appropriate even section $\alpha_{i j k}$ of $\rund{\rund{\bigwedge^n \cz^{r_D + s_D}} \boxtimes sT^\rel \M_{g, D}}^{V_{g, D} \times \cz^\times}$~, expressed in the local super chart $i$~. Without restriction we may assume that $\Psi_{j i}~=~\Psi_{i j}^{- 1}$~, and then $\rund{\alpha_{i j k}}$ becomes a $2$-cocycle in the sheaf $\rund{\rund{\bigwedge^n \cz^{r_D + s_D}} \boxtimes sT^\rel \M_{g, D}}_0^{V_{g, D} \times \cz^\times}$ on $M_{g, D}$~, which is easily seen to be a direct sum of copies of $\rund{sT^\rel \M_{g, D}}_0$~, $T^\rel M_{g, D} \otimes L_{g, D}^*$ and $L_{g, D}$ and so is locally free. Therefore by \cite{GrauRem} theorem 10.5.5

\[
\rund{\pi_{M_{g, D}} }_{(2)} \rund{\rund{\bigwedge\nolimits^n \cz^{r_D + s_D}} \boxtimes sT^\rel \M_{g, D}}_0^{V_{g, D} \times \cz^\times} = 0 \,~,
\]

and so after a refinement of the atlas of $\M_{g, D}$~, locally in $V_{g, D}$~, $\rund{\alpha_{i j k}} = d \rund{\eps_{i j}}$ with an even $1$-cochain $\rund{\eps_{i j}}$ in $\rund{\rund{\bigwedge^n \cz^{r_D + s_D}} \boxtimes sT^\rel \M_{g, D}}^{V_{g, D} \times \cz^\times}$~. Let $\eps_{i j}$ be expressed in the local super chart $i$~. Then

\[
\Phi_{i j}^{D, n + 1} := \Psi_{i j} \circ \rund{\Id_{\left.\rund{\V_{n + 1}^D \times \cz^{|1}}\right|_{U_{i j, D}} } - \eps_{i j}} : \left.\rund{\V_{n + 1}^D \times \cz^{|1}}\right|_{U_{i j, D}} \mathop{\rightarrow}\limits^\sim \left.\rund{\V_{n + 1}^D \times \cz^{|1}}\right|_{U_{j i, D}}
\]

is $\rund{\V_{n + 1}^D \times \cz^\times}$-equivariant, and $\Phi_{k i}^{D, n + 1} \circ \Phi_{j k}^{D, n + 1} \circ \Phi_{i j}^{D, n + 1} =  \Id_{\left.\rund{\V_{n + 1}^D \times \cz^{|1}}\right|_{U_{i j k, D}} }$~. So locally in $\V_{n + 1}^D$~, $\rund{\Phi_{i j}^{D, n + 1} }$ are the transition morphisms of a super family $\M_{n + 1}^D \twoheadrightarrow \V_{n + 1}^D$ having $\left.\M_{n + 1}^D\right|_{\V_n^D} = \M_n^D$~.

{\it Global construction:} On the overlaps of two such local constructions their transition morphisms differ by a section of the sheaf $\rund{\pi_{M_{g, D}}}_{(1)} \rund{\rund{\bigwedge\nolimits^n \cz^{r_D + s_D}} \boxtimes sT^\rel \M_{g, D}}_0^{V_{g, D} \times \cz^\times}$~on $V_{g, D}$~, which is locally free by \cite{GrauRem} theorem 10.5.5. We see that the obtructions to define $\M_{n + 1}^D$ globally lie in 

\[
H^1\rund{V_{g, D}, \rund{\pi_{M_{g, D}}}_{(1)} \rund{\rund{\bigwedge\nolimits^n \cz^{r_D + s_D}} \boxtimes sT^\rel \M_{g, D}}_0^{V_{g, D} \times \cz^\times}} \,~,
\]

which vanishes by Cartan's theorem B since $V_{g, D}$ is Stein.






Finally since $\rund{\rund{\left.\omega_D \right| \vartheta, \eta}^* \Phi_{\widehat i \ \widehat j}^{2 (1 - g) - D, n + 1}}^\vee$ is also $\rund{\V_{n + 1}^D \times \cz^\times}$-equivariant and coincides with $\Phi_{i j}^{D, n + 1}$ on $\left.\rund{\V_n^D \times \cz^{|1}}\right|_{U_{i j, D}}$ after passing from $\Phi_{i j}^{D, n + 1}$ to

\[
\frac{1}{2} \rund{\Phi_{i j}^{D, n + 1} + \rund{\rund{\left.\omega_D \right| \vartheta, \eta}^* \Phi_{\widehat i \ \widehat j}^{2 (1 - g) - D, n + 1}}^\vee}
\]

we may assume that $\rund{\Phi_{i j}^{D, n + 1}}^\vee = \rund{\left.\omega_D \right| \vartheta, \eta}^* \Phi_{\widehat i \ \widehat j}^{2 (1 - g) - D, n + 1}$~. $\Box$ \\

We can immediately deduce:

\begin{theorem}[Versality of $\widetilde \M_{g, D}$ ] \label{versalitysuper}
{\bf Completeness:} Let $\N \twoheadrightarrow \W$ be a super family of compact super Riemann surfaces of type $(g, D)$~. If $\W^\#$ is a contractible Stein manifold then there exists a fiberwise biholomorphic super family morphism $(\Xi, \xi): \N \rightarrow \widetilde \M_{g, D}$~. More precisely we have

{\bf Local completeness:} For every isomorphism $\sigma: \Y_{w_0} := \pi_\N^{- 1}\rund{w_0} \mathop{\rightarrow}\limits^\sim \X_{v_0}$~, $w_0 \in \W^\#$ and $v_0 \in V_{g, D}$~, there exists a fiberwise biholomorphic super family morphism $(\Xi, \xi): \N \rightarrow \widetilde \M_{g, D}$ such that $\xi\rund{w_0} = v_0$ and $\Xi|_{\Y_{w_0}} = \sigma$~. $\xi$ will in general {\bf not} be unique.

{\bf Anchoring property:} For every fiberwise biholomorphic super family morphism \\
$(\Sigma, \sigma): \N|_{\W^\natural} \rightarrow \widetilde \M_{g, D}$~, $\W$ of nilpontency index $\geq 2$ and $\W^\#$ a Stein manifold, there exists a unique morphism $\xi: \W \rightarrow s\T_{g, D}$ such that $\xi|_{\W^\natural} = \sigma$ and $\xi$ can be extended to a fiberwise biholomorphic super family morphism $(\Xi, \xi): \N \rightarrow \widetilde \M_{g, D}$ such that $\Xi|_{\N|_{\W^\natural}} = \Sigma$~.
\end{theorem}

{\it Proof:} The {\it anchoring property} is evident from (i) in theorem \ref{constructsuperfam} and proposition \ref{versalitygeneral}.

{\it Local completeness:} Write $\W := (W, \S)$ and let $\m \lhd \S$ be the nilradical. For $n \in \nz \setminus \{0\}$ define the ringed space $\W_n := \rund{\W, \S \left/ \m^n\right.}$~. Then $\W_1 = W$ and either $\W_n = \W_{n + 1}^\natural$ or $\W_n = \W_{n + 1}$~. By the anchoring property of $\M_{g, D}$ there exists a fiberwise biholomorphic family morphism $\rund{\Xi_1, \xi_1}: \N|_W \rightarrow \M_{g, D}$ such that $\xi_1\rund{w_0} = v_0$ and $\left.\Xi_1\right|_{\Y_{w_0}} = \sigma$~. By the anchoring property of $\widetilde \M_{g, D}$ there exist fiberwise biholomorphic $\rund{\Xi_n, \xi_n}: \N|_{\W_n} \rightarrow \widetilde \M_{g, D}$~, $n \geq 2$~, such that $\left.\xi_n\right|_{\W_{n - 1}} = \xi_{n - 1}$ and $\left.\Xi_n\right|_{\N|_{\W_{n - 1}} } = \Xi_{n - 1}$~. So all $\rund{\Xi_n, \xi_n}$ glue together to a fiberwise biholomorphic family morphism $(\Xi, \xi): \N \rightarrow \widetilde \M_{g, D}$ since $\m$ is locally nilpotent.

{\it Completeness} now is evident from the completeness of $\M_{g, D}$ and local completeness of $\widetilde \M_{g, D}$~. $\Box$ \\

Property (ii) in theorem \ref{constructsuperfam} yields an obvious embedding

\[
\O\rund{s\T_{g, D}}_0^\times \hookrightarrow \Aut \ \widetilde \M_{g, D} \,~, \, f \mapsto \rund{\widetilde \Phi_{g, D, f}, \widetilde \varphi_{g, D, f}} \,~,
\]

where the $\widetilde \Phi_{g, D, f}$ in the standard local super charts of $\widetilde \M_{g, D}$ are given by $\rund{\id_{V_{g, D}}, z \left| \frac{1}{f} \eta, f \vartheta, f \zeta\right.}$ and extend the $\Phi_{\M_{g, D}, f}$ from (\ref{exact}) in section \ref{basics} as in general {\bf not} strong super family automorphism. Here $\O\rund{s\T_{g, D}}_0^\times$ is equipped with the twisted product $f \bullet h := \rund{f \circ \rund{\id_{V_{g, D}} \left| \frac{1}{h} \eta, h \vartheta\right.}} h$ for all $f, h \in \O\rund{s\T_{g, D}}_0^\times$~, not the usual one! \\

Infact we have broken the $\O\rund{V_{g, D}}^\times$-symmetry of $\M_{g, D}$~:

\begin{quote}
Let $f \in \O\rund{V_{g, D}}^\times \setminus \schweif{1}$ and assume $\Xi \in \Aut_\strong \ \M_2^D$ with $\Xi|_{\M_{g, D}} = \Phi_{\M_{g, D}, f}$ (notation as in the proof of theorem \ref{constructsuperfam}). Then both fiberwise biholomorphic super family morphisms $\rund{\Xi, \Id_{\V_2^D}}$ and $\rund{\left.\widetilde \Phi_{g, D, f}\right|_{\M_2^D}, \left.\widetilde \varphi_{g, D, f}\right|_{\V_2^D}}: \M_2^D \hookrightarrow \widetilde \M_{g, D}$ coincide with $\Phi_{\M_{g, D}, f}$ on $\M_{g, D}$~. Since $\rund{\V_2^D}^\natural = V_{g, D}$ and $\left.\widetilde \varphi_{g, D, f}\right|_{\V_2^D} \not= \Id_{\V_2^D}$ we have a contradiction to the anchoring property of $\widetilde \M_{g, D}$~.
\end{quote}

As a consequence one cannot expect an essentially stronger anchoring property for $\widetilde \M_{g, D}$~:

\begin{quote}
Let $f|_{\V_{2 n}^D} = 1$ but $f|_{\V_{2 n + 1}^D} \not= 1$ for some $f \in \O\rund{s\T_{g, D}}_0^\times$ and $n \in \schweif{1, \dots, \floor{\frac{r_D + s_D - 1}{2}} }$ (in other words $f$ congruent to $1$ modulo $\rund{\bigwedge^{2 n} \cz^{r_D + s_D}}$ but {\bf not} mo\-du\-lo $\rund{\bigwedge^{2 n + 1} \cz^{r_D + s_D}}$~). Then $\rund{\widetilde \Phi_{g, D, f}, \widetilde \varphi_{g, D, f}}$ and $\Id_{\widetilde \M_{g, D}} \in \Aut \ \widetilde \M_{g, D}$ coincide on $\M_{2 n}^D$ while already $\left.\widetilde \varphi_{g, D, f}\right|_{\V_{2 n + 2}^D} \not=~\Id_{\V_{2 n + 2}^D}$~.
\end{quote}

\begin{theorem} \label{superfamunique}
\item[(i)] Every fiberwise biholomorphic super family endomorphism of $\widetilde \M_{g, D}$ is infact an automorphism.
\item[(ii)] $\widetilde \M_{g, D} \twoheadrightarrow s\T_{g, D}$ is up to super family isomorphism uniquely determined by the completeness of theorem \ref{versalitysuper}, the infinitesimal universality (i) in theorem \ref{constructsuperfam}, and $s\T_{g, D}^\#$ being a contractible Stein manifold.
\end{theorem}

{\it Proof:} (i) Let $(\Xi, \xi)$ be a fiberwise biholomorphic super family endomorphism of $\widetilde \M_{g, D}$~. Then $\rund{\Xi|_{\M_{g, D}}, \xi^\#} \in \Aut \ \M_{g, D}$ by \ref{ordfamunique} (i), in particular $\xi^\# \in \Aut \ V_{g, D}$~. Furthermore \\
$(d \xi)_v: sT_v s\T_{g, D} \rightarrow sT_{\xi^\#(v)} s\T_{g, D}$ is an isomorphism by the infinitesimal universality of $\widetilde \M_{g, D}$~, and so $\xi \in \Aut \ s\T_{g, D}$ by the super inverse function theorem.

(ii) Let $\N \twoheadrightarrow \W = (W, \S)$ be another super family of compact super Riemann surfaces of type $(g, D)$ with the same completeness and infinitesimal universality property as $\widetilde \M_{g, D}$ and $\W^\#$ being a contractible Stein manifold. Then by the completeness of $\N$ and $\widetilde \M_{g, D}$ we obtain fiberwise biholomorphic super family morphisms

\[
\widetilde \M_{g, D} \mathop{\longrightarrow}\limits^{(\Xi, \xi)} \N \mathop{\longrightarrow}\limits^{(\Phi, \varphi)} \widetilde \M_{g, D} \,~,
\]

whose composition is an automorphism of $\widetilde \M_{g, D}$ by (i). Therefore $\varphi: \W \rightarrow s\T_{g, D}$ is a projection, and we have to show that it is infact an isomorphism. Furthermore by the infinitesimal universality of $\N$ and $\widetilde \M_{g, D}$ we see that $(d \xi)_v: sT_v s\T_{g, D} \rightarrow sT_{\xi^\#(v)} \W$ and $(d \varphi)_w: sT_w \W \rightarrow sT_{\varphi^\#(w)} s\T_{g, D}$ are isomorphisms for all $v \in V_{g, D}$ and $w \in W$~. This shows that first of all $\varphi^\#$ is locally biholomorpic. But $W$ is connected and $V_{g, D}$ simply connected, therefore $\varphi^\#$ must infact be globally biholomorphic, and so we may assume without restriction that $W = V_{g, D}$ and $\varphi^\# = \id_{V_{g, D}}$~. Now assume $k \in \nz \setminus \{0\}$ to be minimal such that $\H := \rund{\left.\m^k \left/ \m^{k + 1}\right.} \right/ \Im \psi_k \not= 0$~, where $\m \lhd \S$ denotes the nilradical and

\[
\psi_k: \O_{V_{g, D}} \otimes \bigwedge\nolimits^k \cz \hookrightarrow \m^k \left/ \m^{k + 1}\right. \,~, \, f \mapsto f \circ \varphi \,~.
\]

Then since $\H$ is coherent and graded there exist $v \in V_{g, D}$ and an even or odd $v$-linear projection $\lambda: \H_v \twoheadrightarrow \cz$ (which means $\lambda (h S) = h(v) \lambda S$ for all $h \in \rund{\O_{V_{g, D}}}_v$ and $S \in \H_v$~). Now easy calculations show that we obtain an even resp. odd derivation $\delta \in sT_v \W \setminus \{0\}$ by setting

\[
\delta (f \circ \varphi) := 0 \, \text{ and } \, \delta h := \lambda \pi_v h
\]

for all $f \in \rund{\O_{s\T_{g, D}} }_v$ resp. $h \in \m^k_v$ with the canonical graded projection $\pi: \m^k \twoheadrightarrow \m^k \left/ \m^{k + 1}\right.$~. Obviously $(d \varphi)_v \delta = 0$ in contradiction to $(d \varphi)_v$ being an isomorphism. Therefore all $\psi_k$ must be surjective, which shows that $\varphi$ is an isomorphism. $\Box$ \\

{\bf For the rest of this section we assume $g \geq 2$ and $D \in \schweif{3 - 2 g, \dots, - 1}$~.} Then $H^0\rund{sT \X_v}$~, the space of infinitesimal automorphisms of $\X_v$~, $v \in V_{g, D}$ arbitrary, is purely even and generated by the global super vector field $\zeta \partial_\zeta$~, and so one can easily deduce the super moduli space for type $(g, D)$~:

\begin{lemma} \label{uniquemod}
Let $\W$ be an admissible base and $\xi,~\tau:~\W~\rightarrow~s\T_{g, D}$ morphisms such that \\
$\xi^\#~=~\tau^\#$~. Then there exists a strong super family isomorphism $\Phi: \xi^* \widetilde \M_{g, D} \mathop{\rightarrow}\limits^\sim \tau^* \widetilde \M_{g, D}$ with $\Phi^\# = \id_{\rund{\xi^\#}^* M_{g, D}}$ iff there exists $f \in \O(\W)_0^\times$ such that $\tau = \rund{\id_{V_{g, D}} \left| \frac{1}{f} \eta, f \vartheta\right.} \circ \rund{\Id_\W, \xi}$~. Moreover, in this case $f$ can be chosen such that $\Phi$ is given by $\rund{\Id_\W, z \ | \ f \zeta}$ in the standard local super coordinates coming from $\widetilde \M_{g, D}$~, and this additional property determines $f$ uniquely.
\end{lemma}

{\it Proof:} `$\Leftarrow$': obvious by the action of $s\T_{g, D} \times \cz^\times$ on $\widetilde \M_{g, D}$~.

`$\Rightarrow$': {\it Uniqueness} is obvious.

{\it Existence:} locally in $\W$ by induction over the nilpotency index $k$ of $\W = (W, \S)$~. Let $\m \lhd \S$ be the nilradical.

$k = 1$ : Then $\W = W$ is an ordinary complex manifold, and so $\Phi = \Phi_{\xi^* \M_{g, D}, f}$ with some $f \in \O(W)^\times$~.

$k \rightarrow k + 1$ : Assume $\W$ is of nilpotency index $k + 1$ and define $\N := \xi^* \widetilde \M_{g, D}$~. Then by induction hypothesis there exists $f \in \O\rund{\W^\natural}_0^\times$ such that $\tau|_{\W^\natural} = \rund{\id_{V_{g, D}} \left| \frac{1}{f} \eta, f \vartheta\right.} \circ \rund{\Id_{\W^\natural}, \xi|_{\W^\natural}}$ and $\Phi|_{\left.\N\right|_{\W^\natural}}$ is given by $\rund{\left.\Id_\W~, z \right| f \zeta}$ in the standard local super coordinates coming from $\widetilde \M_{g, D}$~. Locally in $\W$ we can take $\widehat f \in \O(\W)_0^\times$ such that $\left.\widehat f \right|_{\W^\natural} = f$ and so by using the action of $s\T_{g, D} \times~\cz^\times$ on $\widetilde \M_{g, D}$ define $\widehat\tau := \rund{\id_{V_{g, D}} \left| \frac{1}{\widehat f} \eta, \widehat f \vartheta\right.} \circ \rund{\Id_\W, \xi}: \W \rightarrow s\T_{g, D}$ and the strong isomorphism $\Sigma: \N \mathop{\rightarrow}\limits^\sim {\widehat\tau}^* \widetilde \M_{g, D}$ by $\rund{\Id_\W, z \left| \widehat f \zeta\right.}$ in the standard local super coordinates coming from $\widetilde \M_{g, D}$~. Then $\widehat\tau$ and $\tau$ coincide on $\W^\natural$ and the fiberwise biholomorphic super family morphisms $\rund{\rund{\widehat\tau}^\smile_{\widetilde \M_{g, D}} \circ \Sigma, \widehat\tau}$ and $\rund{\breve \tau_{\widetilde \M_{g, D}} \circ \Phi, \tau}: \N \rightarrow \widetilde \M_{g, D}$ on $\N|_{\W^\natural}$~. By the anchoring property of $\widetilde \M_{g, D}$ we see that $\widehat \tau = \tau$~, and so $\Phi = \Sigma \circ \rund{\Id_\N + \delta}$ with a suitable even section $\delta$ of

\[
\rund{\pi_\N^\#}_* \rund{\m^k \boxtimes_{\O_W} sT^\rel \N|_W} = \m^k \boxtimes_{\O_W} \rund{\pi_\N^\#}_* sT^\rel \N|_W
\]

identified via the canonical graded isomorphism $\lambda_{\m^k, 0}$ of lemma \ref{coherent iso} (i) with $\H := \m^k$~, $F := sT^\rel \N|_W$~, and $\pi := \pi_\N^\#$~. Since we are in the case $g \geq 2$ and $D \in \schweif{3 - 2 g, \dots, - 1}$~, infact $\delta$ is given by $A \zeta \partial_\zeta$ for some $A \in \rund{\m^k}_0$~, and so $\Phi$ in the standard local super coordinates coming from $\widetilde \M_{g, D}$ by $\rund{\left.\Id_\W, z \right| \widetilde f \zeta}$~, $\widetilde f := \widehat f (1 + A)~\in~\O(\W)_0^\times$~. Furthermore also

\[
\tau = \rund{\id_{V_{g, D}} \left| \frac{1}{\widehat f} \eta, \widehat f \vartheta\right.} \circ \rund{\Id_\W, \xi} = \rund{\id_{V_{g, D}} \left| \frac{1}{\widetilde f} \ \eta, \widetilde f \vartheta\right.} \circ \rund{\Id_\W, \xi}
\]

since products of $A$ with odd elements of $\O(\W)$ vanish.

Finally $f$ is globally defined by its uniqueness. $\Box$ \\



In particular we obtain an exact sequence

\begin{eqnarray} \label{exact super}
1 \hookrightarrow \O\rund{s\T_{g, D}}_0^\times &\hookrightarrow& \Aut \ \widetilde \M_{g, D} \mathop{\longrightarrow}\limits^{{}^\#} \Aut \ M_{g, D} \,~. \\
\notag f &\mapsto& \widetilde \Phi_{g, D, f}
\end{eqnarray}

\begin{defin}

\[
\widetilde \Gamma_{g, D} := \schweif{\gamma \, \left| \, (\Gamma, \gamma) \in \Aut \ \widetilde \M_{g, D}\right.} \left/ \O\rund{s\T_{g, D}}_0^\times \right.
\]

is called the super modular group for type $(g, D)$~. It acts on the super orbifold $s\T_{g, D} \left/ \rund{s\T_{g, D} \times \cz^\times}\right.$~.
\end{defin}

\begin{theorem}[Super moduli space for type $(g, D)$ ] \label{supermod}
\item[(i)] $\rund{s\T_{g, D} \left/ \rund{s\T_{g, D} \times \cz^\times}\right.}\left/ \ \widetilde \Gamma_{g, D}\right.$ is the super moduli space for type $(g, D)$~, more precisely: Let $\W$ be a connected admissible base and $\xi, \tau: \W \rightarrow s\T_{g, D}$ morphisms. Then 
$\xi^* \widetilde \M_{g, D} \simeq~\tau^* \widetilde \M_{g, D}$ strongly iff there exist $f \in \O(\W)_0^\times$ and $\gamma \in \Aut \ s\T_{g, D}$ representing an element of $\widetilde \Gamma_{g, D}$ such that

\[
\gamma \circ \tau = \rund{\id_{V_{g, D}}^{\phantom{,}} \ \left| \ \frac{1}{f} \ \eta, f \vartheta\right.} \circ \rund{\Id_\W, \xi} \,~.
\]

\item[(ii)] We have isomorphisms

\[
\begin{array}{ccc}
\left.\rund{\Aut \ \widetilde \M_{g, D}}\right/ \O\rund{s\T_{g, D}}_0^\times & \mathop{\longrightarrow}\limits^{|_{\M_{g, D}} } & \rund{\Aut \ \M_{g, D}} \left/ \O\rund{V_{g, D}}^\times \right. \\
\downarrow & \circlearrowleft & \downarrow \\
\widetilde \Gamma_{g, D} & \mathop{\longrightarrow}\limits_{{}^\#} & \Gamma_{g, D}
\end{array} \,~,
\]

where the vertical maps are given by the canonical projections $(\Gamma, \gamma) \mapsto \gamma$~.
\end{theorem}

{\it Proof:} (i) `$\Rightarrow$': Let $\Phi: \xi^* \widetilde \M_{g, D} \mathop{\rightarrow}\limits^\sim \tau^* \widetilde \M_{g, D}$ be a strong isomorphism and $w_0 \in W := \W^\#$ arbitrary. Write $\N := \xi^* \widetilde \M_{g, D}$ and $\Y_{w_0} := \pi_\N^{- 1}\rund{w_0} = \X_{\xi^\#\rund{w_0}}$~. Then by the local completeness of $\widetilde \M_{g, D}$ and theorem \ref{superfamunique} (i) there exists $(\Gamma, \gamma) \in \Aut \ \widetilde \M_{g, D}$ such that $\gamma^\#\rund{\tau^\#\rund{w_0}} = \xi^\#\rund{w_0}$ and $\Gamma|_{\X_{\tau^\#\rund{w_0}} } = \rund{\Phi|_{\Y_{w_0}} }^{- 1}$~. Let $\widehat\Gamma := \left.\rund{\Id_\W, \Gamma}\right|_{\tau^* \widetilde \M_{g, D}}$ be the pullback of $\Gamma$ under $\tau$ and $\gamma \circ \tau$~, so the unique strong fiberwise biholomorphic morphism $\tau^* \widetilde \M_{g, D} \rightarrow (\gamma \circ \tau)^* \widetilde \M_{g, D}$ making the diagram

\[
\begin{array}{ccc}
\phantom{123456,} \tau^* \widetilde \M_{g, D} & \mathop{\longrightarrow}\limits^{\widehat\Gamma} & (\gamma \circ \tau)^* \widetilde \M_{g, D} \phantom{1234567890} \\
\breve \tau_{\widetilde \M_{g, D}} \downarrow & \circlearrowleft & \downarrow (\gamma \circ \tau)^\smile_{\widetilde \M_{g, D}} \\
\phantom{12345,} \widetilde \M_{g, D} & \mathop{\longrightarrow}\limits_\Gamma & \widetilde \M_{g, D} \phantom{1234567890}
\end{array}
\]

commutative. Define the fiberwise biholomorphic super family morphism

\[
(\Omega, \gamma \circ \tau) := (\gamma \circ \tau)^\smile_{\widetilde \M_{g, D}} \circ \widehat\Gamma \circ \Phi = \Gamma \circ \breve \tau_{\widetilde \M_{g, D}} \circ \Phi: \ \N \rightarrow \widetilde \M_{g, D} \,~.
\]

Then $\rund{\Omega|_{\N|_W}, (\gamma \circ \tau)^\#}$ and $\rund{\rund{\xi^\#}^\smile_{\M_{g, D}}, \xi^\#}: \N|_W = \rund{\xi^\#}^* \M_{g, D} \rightarrow \M_{g, D}$ coincide on $\Y_{w_0}$~, and so $(\tau \circ \gamma)^\# = \xi^\#$ by the anchoring property of $\M_{g, D}$ and since $W$ is connected. Since moreover $\pr_{M_g} \circ \Omega^\#$ and $\rund{\pr_{M_g} \circ \xi^\#}^\smile_{M_g}: \N^\# = \rund{\pr_{\T_g} \circ \xi^\#}^* M_g \rightarrow M_g$ coincide on $\Y_{w_0}^\#$ we even have $\Omega^\# = \rund{\xi^\#}^\smile_{M_{g, D}}$ by the anchoring property of $M_g$~. Therefore \\
$\widehat \Gamma \circ \Phi: \N \rightarrow (\gamma \circ \tau)^* \widetilde \M_{g, D}$ is a fiberwise biholomorphic super family morphism with \\ $\rund{\widehat \Gamma \circ \Phi}^\#~=~\id_{\N^\#}$~. So by lemma \ref{uniquemod} there exists $f \in \O(\W)_0^\times$ such that \\
$\gamma \circ \tau = \rund{\id_{V_{g, D}} \left| \frac{1}{f} \eta, f \vartheta\right.} \circ \rund{\Id_\W, \xi}$ and $\widehat\Gamma \circ \Phi$ is given by $\rund{\Id_\W, z \ | \ f \zeta}$ in the standard local super coordinates coming from $\widetilde \M_{g, D}$~.






`$\Leftarrow$': Easy exercise using lemma \ref{uniquemod} and the definition of $\widetilde \Gamma_{g, D}$~.

(ii) Commutativity of the diagram is evident. The first row is surjective by the anchoring property of $\widetilde \M_{g, D}$ and theorem \ref{superfamunique} (i) and injective by (\ref{exact super}). By corollary \ref{autstrong} the projection on the right is even an isomorphism. $\Box$ \\






\section{Supersymmetric super families} \label{superfamSUSY}

In this section we construct versal supersymmetric super families of compact super Riemann surfaces as sub super families of $\widetilde \M_{g, 1 - g}$ containing the versal supersymmetric families $\bigcup_{\eps \in \frac{1}{2} \zz^{2 g}} \M_g^\eps$ constructed in section \ref{ordfamSUSY}. {\bf Obviously we have to exclude the case $g = 1$~,} it will be treated separately in section \ref{remain}.

\[
r_{1 - g} = s_{1 - g} = \left\{\begin{array}{l}
0 \text{ if } g = 0 \\
2 (g - 1) \text{ if } g \geq 2
\end{array}\right.~.
\]

Let $s\T_g^0$ be the image of the embedding $(u, 0 | \chi, \chi): \T_g^{| 2 (g - 1)} \hookrightarrow s\T_{g, 1 - g}$ if $g \geq 2$ and $s\T_0^0 := s\T_{0, 1}$~, which is just a single point. Given the super family isomorphism $\rund{\widetilde \Omega_{g, D}, \widetilde \omega_{g, D}}: {\widetilde \M_{g, D}}^\vee \mathop{\rightarrow}\limits^\sim \widetilde \M_{g, 2 (g - 1) - D}$ by (iii) in theorem \ref{constructsuperfam}, obviously $s\T_g^0$ is the fixed point subsupermanifold of $\widetilde \omega_{g, 1 - g}$~, $\rund{s\T_g^0}^\# = \T_g^0$ and ${\widetilde \Omega_{g, D}}^\vee = {\widetilde \Omega_{g, 2 (1 - g) - D}}^{- 1}$~. By lemma \ref{respectSUSY} (ii) and theorem \ref{charactSUSY} we see that $\widetilde \M_g^0 := \left.\widetilde \M_{g, 1 - g} \right|_{s\T_g^0}$ is supersymmetric and the standard local super charts of $\widetilde \M_{g, 1 - g}$ give a standard atlas of $\widetilde \M_g^0$~. \\

Let $g \geq 2$ and $\eps \in \frac{1}{2} \zz^{2 g} \setminus \schweif{0}$~. Can we also find a complex subsupermanifold $s\T_0^\eps \simeq \T_g^{| 2 (g - 1)}$ of $s\T_{g, 1 - g}$ with body $\T_g^\eps$ such that $\left.\widetilde \M_{g, 1 - g}\right|_{s\T_0^\eps}$ is supersymmetric? The answer is yes: Taking $E_g \otimes \rund{\id_{M_g}, \id_{\cz^g} + \frac{\eps^\mu \lambda_\mu}{1 - g}}^* L_g$ instead of $E_g$ the construction procedure of theorem \ref{constructsuperfam} leads to a super family $\widehat \M_{g, D}^\eps \twoheadrightarrow s\T_{g, D}$~, a super family isomorphism \\ $\rund{\widehat \Omega_{g, D}^\eps, \widehat \omega_{g, D}^\eps}: \rund{\widehat \M_{g, D}^\eps}^\vee \mathop{\rightarrow}\limits^\sim \widehat \M_{g, 2 (1 - g) - D}^\eps$ and an embedding

\[
\O\rund{s\T_{g, D}}_0^\times \hookrightarrow \Aut \ \widehat \M_{g, D}^\eps \,~, \, f \mapsto \rund{\widehat \Phi_{g, D, f}, \widehat \varphi_{g, D, f}} \,~.
\]

On the other hand by the anchoring property of $\widetilde \M_{g, D}$ and theorem \ref{superfamunique} (i) there exist super family isomorphisms $\rund{\Delta_{g, D}^\eps, \delta_{g, D}^\eps}: \widetilde \M_{g, D} \mathop{\rightarrow}\limits^\sim \widehat \M_{g, D}^\eps$ such that ${\Delta_{g, D}^\eps}^\# = \rund{\id_{\cz^g} - D \frac{\eps^\mu \lambda_\mu}{1 - g}, \id_{M_g}}$~. Now define the super family isomorphism

\[
\rund{\widetilde \Omega_{g, D, \eps}, \widetilde \omega_{g, D, \eps}} := \rund{\Delta_{g, 2 (1 - g) - D}^\eps}^{- 1} \circ \widehat \Omega_{g, D}^\eps \circ \rund{\Delta_{g, D}^\eps}^\vee : \widetilde \M_{g, D}^\vee \mathop{\rightarrow}\limits^\sim \widetilde \M_{g, 2 (1 - g) - D} \,~.
\]

Then ${\widetilde \Omega_{g, D, \eps}}^\vee = {\widetilde \Omega_{g, 2 (1 - g) - D, \eps}}^{- 1}$ and without restriction $\left.\widetilde \Omega_{g, D, \eps} \right|_{\M_{g, D}^\vee} = \Omega_{g, D, \eps}$~. Furthermore $s\T_g^\eps := \rund{\delta_{g, 1 - g}^\eps}^{- 1} s\T_g^0$ is obviously the fixed point subsupermanifold of $\widetilde \omega_{g, 1 - g, \eps}$~, \\
$\rund{s\T_g^\eps}^\# = \T_g^\eps$~, and by the same reason as before $\widetilde \M_g^\eps := \left.\widetilde \M_{g, 1 - g} \right|_{s\T_g^\eps} \simeq \left.\widehat \M_{g, 1 - g}^\eps \right|_{s\T_g^0}$ is supersymmetric.

\begin{theorem}[Versality of $\dot\bigcup_{\eps \in \E_g} \widetilde \M_g^\eps$ ] 
\item[(i)] {\bf Completeness:} Let $\N \twoheadrightarrow \W$ be a supersymmetric super family of compact super Riemann surfaces of genus $g$~. If $\W^\#$ is a contractible Stein manifold then there exists a fiberwise biholomorphic super family morphism $(\Xi, \xi): \N \rightarrow \dot\bigcup_{\eps \in \E_g} \widetilde \M_g^\eps$~. More precisely we have

{\bf Local completeness:} For every supersymmetric isomorphism $\sigma: \Y_{w_0} := \pi_\N^{- 1}\rund{w_0} \mathop{\rightarrow}\limits^\sim \X_{v_0}$~, $w_0 \in \W^\#$~, $\eps \in \frac{1}{2} \zz^{2 g}$, and $v_0 \in \T_g^\eps$~, there exists a supersymmetric fiberwise biholomorphic super family morphism $(\Xi, \xi): \N \rightarrow \widetilde \M_g^\eps$ such that $\xi\rund{w_0} = v_0$ and $\Xi|_{\Y_{w_0}} = \sigma$~.

and the

{\bf Anchoring property:} For every supersymmetric fiberwise biholomorphic super fa\-mi\-ly morphism $(\Sigma, \sigma): \N|_{\W^\natural} \rightarrow \widetilde \M_g^\eps$~, $\eps \in \frac{1}{2} \zz^{2 g}$~, $\W$ of nilpontency index $\geq 2$ and $\W^\#$ a Stein manifold, there exists a unique morphism $\xi:~\W~\rightarrow s\T_g^\eps$ such that $\xi|_{\W^\natural} = \sigma$ and $\xi$ can be extended to a supersymmetric fiberwise biholomorphic super family morphism \\
$(\Xi, \xi): \N \rightarrow~\widetilde \M_g^\eps$ such that $\Xi|_{\N|_{\W^\natural}} = \Sigma$~.
\item[(ii)] {\bf Infinitesimal universality:} $\rund{d \widetilde \M_{g, 1 - g}}_v sT_v \ s\T_g^\eps = H^1\rund{sT^+ \X_v}$ for all $\eps \in \frac{1}{2} \zz^{2 g}$ and $v \in \T_g^\eps$~.
\end{theorem}

{\it Proof:} (i) {\it Anchoring property:} {\it Uniqueness} is obvious by the anchoring property of $\widetilde \M_{g, 1 - g}$~.

{\it Existence:} By maybe passing from $\widetilde \M_{g, 1 - g}$ to $\widehat \M_{g, 1 - g}^\eps$ we may assume $\eps = 0$~. Let \\
$\rund{\Xi, \xi}: \N \rightarrow \widetilde \M_{g, 1 - g}$ be given by theorem \ref{versalitysuper} anchoring property. In a standard atlas of $\N$ the identities in the local super charts glue together to a strong isomorphism $\Theta: \N \mathop{\rightarrow}\limits^\sim \N^\vee$~. Therefore

\[
\rund{\widehat\Xi, \widehat\xi} := \widetilde \Omega_{g, 1 - g} \circ \Xi^\vee \circ \Theta
\]

is another fiberwise biholomorphic super family morphism $\N \rightarrow \widetilde \M_{g, 1 - g}$ having \\
$\left.\widehat\Xi\right|_{\N|_{\W^\natural}} = \Xi|_{\N|_{\W^\natural}} = \Sigma$ since $\Sigma$ is already supersymmetric and $\sigma$ maps to $s\T_g^0$~. Therefore $\widehat\xi = \xi$ by the anchoring property of $\widetilde \M_{g, 1 - g}$~. On the other hand $\widehat\xi = \widetilde \omega_{g, 1 - g} \circ \xi$~, and so $\xi$ must map to $s\T_g^0$~. Furthermore let $\Xi$ in standard atlasses of $\N$ and $\widetilde \M_g^0$ be given by $\Xi_{i j}$~. Then all $\left.\Xi_{i j}\right|_{\W^\natural \times \cz^{|1}}$ are already self-dual, and so by taking $\frac{1}{2} \rund{\Xi_{i j} + \Xi_{i j}^\vee}$ instead of $\Xi_{i j}$ we obtain $\Xi_{i j}^\vee = \Xi_{i j}$~, and so $\Xi$ itself becomes supersymmetric.

The rest follows as in the proof of theorem \ref{versalitysuper} using the anchoring property and completeness of $\dot\bigcup_{\eps \in \frac{1}{2} \zz^{2 g}} \M_g^\eps$~.

(ii) same as for theorem \ref{versalfamSUSY} (ii). $\Box$ \\

\begin{cor}[Local uniqueness of supersymmetries] \label{SUSYunique}
Let $\N \twoheadrightarrow \W$ be a super fa\-mi\-ly of compact super Riemann surfaces of genus $g \not= 1$ with two supersymmetries $\D$ and $\D'$~. Then locally in $\W$ there exists a strong supersymmetric super family isomorphism $\Phi: (\N, \D) \mathop{\rightarrow}\limits^\sim \rund{\N, \D'}$ with $\Phi^\# = \id_{\N^\#}$~.
\end{cor}

In the end corollary \ref{SUSYunique1} will show that the statement is also true for genus~$1$~. \\

{\it Proof:} trivial if $g = 0$ by corollary \ref{SUSYunique0}. So assume $g \geq 2$~. By the completeness of $\dot\bigcup_{\eps \in \frac{1}{2} \zz^{2 g}} \widetilde \M_g^\eps$ we may assume that there exists - locally in $\W$ - a fiberwise biholomorphic supersymmetric $(\Xi, \xi): (\N, \D) \rightarrow \M_g^0$~. By corollary \ref{SUSYunique0} locally in $W$ there exists a strong supersymmetric isomorphism $\Psi: \rund{\N|_W, \D'} \mathop{\rightarrow}\limits^\sim \rund{\N|_W, \D}$ with $\Psi^\# = \id_{\N^\#}$~, and so

\[
\rund{\Xi|_{\N|_W} \circ \Psi, \xi^\#}: \rund{\N|_W, \D'} \rightarrow \M_g^0
\]

is fiberwise biholomorphic and supersymmetric. Therefore by the anchoring property of $\dot\bigcup_{\eps \in \frac{1}{2} \zz^{2 g}} \widetilde \M_g^\eps$ locally in $\W$ there exists a fiberwise biholomorphic supersymmetric \\
$(\Omega, \omega): \rund{\N, \D'} \rightarrow \widetilde \M_g^0$ with $\Omega^\# = \Xi^\#$~. Write $\xi = (\varphi, \kappa): \W \rightarrow s\T_g^0$ with $\varphi: \W \rightarrow \T_g^0$ and $\kappa \in \O(\W)_1^{\oplus 2 (g - 1)}$~. Then by lemma \ref{uniquemod} locally in $\W$ there exists $f \in \O(\W)_0^\times$ such that $\omega = \rund{\varphi \ \left| \ \frac{1}{f} \kappa, f \kappa\right.}$~, and since $\omega$ maps to $s\T_g^0$ we see that $\rund{1 - f^2} \kappa = 0$~. By maybe taking $\widetilde \Phi_{g, 1 - g, -1} \circ (\Omega, \omega)$ instead of $(\Omega, \omega)$ we may assume without restriction that $f^\# + 1$ has no zero, therefore $(f - 1) \kappa = 0$ and so $\omega = \xi$~. We see that

\[
\Phi := \rund{\pi_\N, \Omega}^{- 1} \circ \rund{\pi_\N, \Xi}: (\N, \D) \rightarrow \rund{\N, \D'} \,~,
\]

where $\rund{\pi_\N, \Omega}: \rund{\N, \D'} \mathop{\rightarrow}\limits^\sim \xi^* \widetilde \M_g^0$ and $\rund{\pi_\N, \Xi}: (\N, \D) \mathop{\rightarrow}\limits^\sim \xi^* \widetilde \M_g^0$ are the strong supersymmetric isomorphisms associated to $(\Omega, \xi)$ resp. $(\Xi, \xi)$~, is a strong supersymmetric isomorphism with $\Phi^\# = \id_{\N^\#}$~. $\Box$ \\

$\dot\bigcup_{\eps \in \frac{1}{2} \zz^{2 g}} \widetilde \M_g^\eps$ is again a maximal sub super family of $\widetilde \M_{g, 1 - g}$ admitting a supersymmetry:

\begin{quote}
trivial if $g = 0$. For $g \geq 2$ let $\W \hookrightarrow s\T_{g, 1 - g}$ be an embedding of an admissible base such that its image contains $\dot\bigcup_{\eps \in \frac{1}{2} \zz^{2 g}} s\T_g^\eps$ and $\D$ a supersymmetry on $\left.\widetilde \M_{g, 1 - g}\right|_{\W}$~. Then $\W^\# = \dot\bigcup_{\eps \in \frac{1}{2} \zz^{2 g}} \T_g^\eps$ by the maximality of $\dot\bigcup_{\eps \in \frac{1}{2} \zz^{2 g}} \M_g^\eps$~, and by corollary \ref{SUSYunique0} there exists a strong supersymmetric isomorphism \\
$\Psi: \rund{\dot\bigcup_{\eps \in \frac{1}{2} \zz^{2 g}} \M_g^\eps, \D} \mathop{\rightarrow}\limits^\sim \rund{\dot\bigcup_{\eps \in \frac{1}{2} \zz^{2 g}} \M_g^\eps, \D'}$~, $\D'$ denoting the usual supersymmetry on $\dot\bigcup_{\eps \in \frac{1}{2} \zz^{2 g}} \widetilde \M_g^\eps$~, with identity as body. By the anchoring property of $\dot\bigcup_{\eps \in \frac{1}{2} \zz^{2 g}} \widetilde \M_g^\eps$ there exists a fiberwise biholomorphic super family morphism

\[
(\Xi, \xi): \left.\widetilde \M_{g, 1 - g}\right|_\W \rightarrow \mathop{\dot\bigcup}\limits_{\eps \in \frac{1}{2} \zz^{2 g}} \widetilde \M_g^\eps
\]

such that $\xi^\# = \id_{\dot\bigcup_{\eps \in \frac{1}{2} \zz^{2 g}} \T_g^\eps}$ and $\Xi|_{\dot\bigcup_{\eps \in \frac{1}{2} \zz^{2 g}} \widetilde \M_g^\eps} = \Psi$~. By the infinitesimal universality of $\widetilde \M_{g, 1 - g}$ and equality of dimensions $\rund{d \ \xi|_{s\T_g^\eps}}_v: sT_v s\T_g^\eps \rightarrow sT_v s\T_g^\eps$ is an automorphism for all $v \in \bigcup_{\eps \in \frac{1}{2} \zz^{2 g}} \T_g^\eps$~, and so $\xi|_{\dot\bigcup_{\eps \in \frac{1}{2} \zz^{2 g}} s\T_g^\eps} \in \Aut \ \dot\bigcup_{\eps \in \frac{1}{2} \zz^{2 g}} s\T_g^\eps$ by the super inverse function theorem.

Now $\Xi$ induces a strong isomorphism $\left.\widetilde \M_{g, 1 - g}\right|_\W \simeq \xi^* \widetilde \M_{g, 1 - g}$ with identity as body, and so by lemma \ref{uniquemod} for every $v \in \dot\bigcup_{\eps \in \frac{1}{2} \zz^{2 g}} \T_g^\eps$ there exists an open neighbourhood $U \subset V_{g, {1 - g}}$ and $f \in \O\rund{\left.s\T_{g, 1 - g}\right|_U}$ such that $\xi$ and $\rund{\id_{V_{g, 1 - g}} \left| \frac{1}{f} \eta, f \vartheta\right.}$ coincide on $\W|_U$~. Therefore $\xi$ is an embedding, and so $\W = \dot\bigcup_{\eps \in \frac{1}{2} \zz^{2 g}} s\T_g^\eps$~.
\end{quote}

However, if $g \geq 2$~, since $\dot\bigcup_{\eps \in \frac{1}{2} \zz^{2 g}} \widetilde \M_g^\eps$ is not invariant under the automorphisms $\Phi_{g, 1 - g, f}$~, $f \in \O\rund{s\T_{g, 1 - g}}_0^\times$~, of $\widetilde \M_{g, 1 - g}$ it cannot be a biggest sub super family admitting a supersymmetry, and so such a biggest one does {\bf not} exist at all. \\

{\bf For the rest of this section we assume $g \geq 2$} and describe the super moduli space of all supersymmetric compact super Riemann surfaces of genus $g$ as a quotient of $\dot\bigcup_{\eps \in \E_g} s\T_g^\eps$~. As a direct consequence of the exact sequence (\ref{exact super}) in section \ref{superfam} we observe

\begin{lemma}[$\rund{s\T_{g, 1 - g} \times \cz^\times}$-equivariance of $\Delta_{g, 1 - g}^\eps$ ] \label{equivariant} For all $\eps \in \frac{1}{2} \zz^{2 g} \setminus \schweif{0}$
\item[(i)] there exists an automorphism $\phantom{1} \widehat{} \phantom{1}$ of $\O\rund{s\T_{g, 1 - g}}_0^\times$ such that \\
$\Delta_{g, 1 - g}^\eps \circ \widetilde \Phi_{g, 1 - g, f} \circ \rund{\Delta_{g, 1 - g}^\eps}^{- 1} =~\widehat \Phi_{g, 1 - g, \widehat f}$ for all $f \in \O\rund{s\T_{g, 1 - g}}_0^\times$~, and
\item[(ii)] $\Delta_{g, 1 - g}^\eps$ is uniquely determined up to composition with $\widetilde \Phi_{g, 1 - g, f}$~, $f \in \O\rund{s\T_{g, 1 - g}}_0^\times$~.
\end{lemma}

From lemma \ref{equivariant} (i) we easily deduce that $\widetilde \Phi_{g, 1 - g, - 1}$ restricts to a supersymmetric super family automorphism of $\widetilde \M_g^\eps$ of order $2$ for all $\eps \in \frac{1}{2} \zz^{2 g}$~. Putting all restrictions together we obtain an exact sequence

\begin{equation} \label{pm 1 SUSY}
1 \hookrightarrow \schweif{\pm 1}^{\frac{1}{2} \zz^{2 g}} \hookrightarrow \Aut_\SUSY \ \mathop{\dot\bigcup}\limits_{\eps \in \frac{1}{2} \zz^{2 g}} \widetilde \M_g^\eps \ \mathop{\longrightarrow}\limits^{{}^\#} \ \Aut \ \mathop{\dot\bigcup}\limits_{\eps \in \frac{1}{2} \zz^{2 g}} M_g \,~,
\end{equation}


indeed:

\begin{quote}
Let $(\Psi, \psi) \in \Aut_\SUSY \ \widetilde \M_g^0$ with $\Psi^\# = \id_{M_g}$~. Then $\psi^\# = \id_{\T_g^0}$~, and the strong supersymmetric isomorphism $\rund{\pi_{\widetilde \M_g^0}, \Psi}: \widetilde \M_g^0 \mathop{\rightarrow}\limits^\sim \psi^* \widetilde \M_g^0$ associated to $(\Psi, \psi)$ has body $\id_{M_g}$~. Therefore by lemma \ref{uniquemod} there exists $f \in \O\rund{s\T_g^0}_0^\times$ such that \\
$\psi = \rund{\id_{V_{g, 1 - g}} \left| \frac{1}{f} \chi, f \chi\right.}$ and $\rund{\pi_{\widetilde \M_g^0}, \Psi}$ is given by $\rund{\left.\Id_{s\T_g^\eps}, z \right| f \zeta}$ in the standard local super coordinates coming from $\widetilde \M_{g, 1 - g}$~. Since $\psi$ maps to $s\T_g^0$ we see that $\rund{1 - f^2} \chi = 0$ and since $\Psi|_{\M_g^0}$ is strong and supersymmetric even $f^\# = \pm 1$~. By maybe passing from $\Psi$ to $\widetilde \Phi_{g, 1 - g, - 1} \circ \Psi$ we may assume that $f^\# = 1$~, and so $(1 - f) \chi = \chi$~, which implies that $\Psi = \rund{\pi_{\widetilde \M_g^0}, \Psi}$ is strong. Finally since $\Psi$ is even supersymmetric we obtain $f = 1$~.
\end{quote}


\begin{prop} \label{restrict action}
Taking suitable representatives for each $\eps \in \frac{1}{2} \zz^{2 g}$ separately and re\-stric\-ting to $\dot\bigcup_{\eps \in \frac{1}{2} \zz^{2 g}} \widetilde \M_g^\eps$ gives an embedding

\[
\left.\rund{\Aut \ \widetilde \M_{g, 1 - g}}\right/ \O\rund{s\T_{g, 1 - g}}_0^\times \hookrightarrow \left.\rund{\Aut_\SUSY \ \mathop{\dot\bigcup}\limits_{\eps \in \frac{1}{2} \zz^{2 g}} \widetilde \M_g^\eps}\right/ \schweif{\pm 1}^{\frac{1}{2} \zz^{2 g}} \,~.
\]

\end{prop}

{\it Proof:} {\it Welldefinedness:} Let $(\Psi, \psi) \in \Aut \ \widetilde \M_{g, 1 - g}$ and $\eps \in \frac{1}{2} \zz^{2 g}$~. Then by the maximality of $\dot\bigcup_{\eps \in \frac{1}{2} \zz^{2 g}} \M_g^\eps$ there exists $\eps' \in \frac{1}{2} \zz^{2 g}$ such that $\psi^\#\rund{\T_g^\eps} = \T_g^{\eps'}$~. We want to show that a suitable representative of $[\Psi]$ in $\left.\rund{\Aut \ \widetilde \M_{g, 1 - g}}\right/ \O\rund{s\T_{g, 1 - g}}_0^\times$ maps $\widetilde \M_g^\eps$ supersymmetrically onto $\widetilde \M_g^{\eps'}$~. By lemma \ref{equivariant} after maybe passing to $\widehat \M_{g, 1 - g}^{\eps'}$ we may assume $\eps' = 0$ without restriction. Furthermore by corollary \ref{SUSYunique0} taking $\widetilde \Phi_{g, 1 - g, h} \circ \Psi$ instead of $\Psi$~, $h \in \O\rund{V_{g, 1 - g}}^\times$ appropriate, we may assume that $\Psi|_{\M_g^\eps}: \M_g^\eps \mathop{\rightarrow}\limits^\sim \M_g^0$ is already supersymmetric. Therefore by the anchoring property of $\widetilde \M_g^0$ there exists a fiberwise biholomorphic supersymmetric morphism $(\Xi, \xi): \widetilde \M_g^\eps \rightarrow \widetilde \M_g^0$ such that $\xi^\# = \psi^\#$ and $\Xi = \Psi|_{\M_g^\eps}$~. Since $\xi^\#$ and $(d \xi)|_{\T_g^\eps}$ are isomorphisms so is also $\xi$ by the super inverse function theorem. Finally by lemma \ref{uniquemod} there exists $f \in \O\rund{s\T_{g, 1 - g}}_0^\times$ such that $\Xi = \left.\widetilde \Phi_{g, 1 - g, f} \circ (\Psi, \psi)\right|_{\widetilde \M_g^\eps}$~.

Conversely let $\widetilde \Phi_{g, 1 - g, f}$ map $\widetilde \M_g^\eps$ supersymmetrically onto itself, $f \in \O\rund{s\T_{g, 1 - g}}_0^\times$~. Then by the exact sequence (\ref{pm 1 SUSY}) we see that $f|_{s\T_g^\eps} = \pm 1$~. \\




{\it Injectivity:} Let $\Psi \in \Aut \ \widetilde \M_{g, 1 - g}$ coincide with $\widetilde \Phi_{g, 1 - g, \pm 1}$ on some $\widetilde \M_g^\eps$~. Then since the first row of (\ref{restrict}) in section \ref{ordfamSUSY} is an embedding we see that $\Psi^\# = \id_{M_{g, 1 - g}}$~, and so $\Psi = \Phi_{g, 1 - g, f}$ for some $f \in \O\rund{s\T_{g, 1 - g}}_0^\times$ by (\ref{exact super}) in section \ref{superfam}. $\Box$ \\

So since by theorem \ref{supermod} (ii) the map ${}^\#$ on the left is an isomorphism, the canonical projection $\left.\rund{\Aut \ \widetilde \M_{g, 1 - g}}\right/ \O\rund{s\T_{g, 1 - g}}_0^\times \twoheadrightarrow \widetilde \Gamma_{g, 1 - g}$ yields a faithful action of $\widetilde \Gamma_{g, 1 - g}$ on the super orbifold $\dot\bigcup_{\eps \in \frac{1}{2} \zz^{2 g}} \rund{\left.s\T_g^\eps \right/ \schweif{\pm 1}}$ such that

\[
\begin{array}{ccc}
\phantom{1} \widetilde \Gamma_{g, 1 - g} & \hookrightarrow & \Aut \ \mathop{\dot\bigcup}\limits_{\eps \in \frac{1}{2} \zz^{2 g}} \rund{\left.s\T_g^\eps \right/ \schweif{\pm 1}} \\
{}^\# \downarrow & \circlearrowleft & \downarrow {}^\# \\
\phantom{1} \Gamma_{g, 1 - g} & \mathop{\hookrightarrow}\limits_{| \ {}_{\bigcup_{\eps \in \frac{1}{2} \zz^{2 g}} \T_g^\eps}} & \Aut \ \mathop{\bigcup}\limits_{\eps \in \frac{1}{2} \zz^{2 g}} \T_g^\eps
\end{array} \,~,
\]

where the first row is given by taking suitable representatives for each $\eps \in \frac{1}{2} \zz^{2 g}$ separately and restricting to $\dot\bigcup_{\eps \in \frac{1}{2} \zz^{2 g}} s\T_g^\eps$~. For every $\eps \in \frac{1}{2} \zz^{2 g}$ the image of the normalizer $\widetilde \Gamma_g^\eps$ of $\left.s\T_g^\eps \right/ \schweif{\pm 1}$ in $\widetilde \Gamma_{g, 1 - g}$ is precisely the normalizer $\Gamma_g^\eps$ of $\T_g^\eps$ in $\Gamma_{g, 1 - g}$~.

\begin{theorem} The super moduli space of all supersymmetric compact super Riemann surfaces of genus $g$ is

\[
\left.\rund{\mathop{\dot\bigcup}\limits_{\eps \in \frac{1}{2} \zz^{2 g}} \rund{\left.s\T_g^\eps \right/ \schweif{\pm 1}} } \right/ \widetilde \Gamma_{g, 1 - g} = \mathop{\dot\bigcup}\limits_{\eps \in \E_g} \rund{\left.\rund{\left.s\T_g^\eps \right/ \schweif{\pm 1}}\right/ \widetilde \Gamma_g^\eps} \,~,
\]

more precisely: Let $\W$ be a connected admissible base and $\xi, \tau: \W \rightarrow \dot\bigcup_{\eps \in \frac{1}{2} \zz^{2 g}} s\T_g^\eps$ morphisms. Then there exists a strong supersymmetric isomorphism $\xi^* \widetilde \M_{g, 1 - g} \simeq~\tau^* \widetilde \M_{g, 1 - g}$ iff $\gamma \circ \tau = \xi$ for some $\gamma \in \Aut \ \dot\bigcup_{\eps \in \frac{1}{2} \zz^{2 g}} s\T_g^\eps$ representing an element of $\widetilde \Gamma_{g, 1 - g}$ in $\Aut \ \dot\bigcup_{\eps \in \frac{1}{2} \zz^{2 g}} \rund{\left.s\T_g^\eps \right/ \schweif{\pm 1}}$~.
\end{theorem}

{\it Proof:} `$\Rightarrow$': similar to the proof of theorem \ref{supermod} (i) `$\Rightarrow$'. Let $\Phi: \N := \xi^* \widetilde \M_g^0 \mathop{\rightarrow}\limits^\sim \tau^* \widetilde \M_g^\eps$ be a strong supersymmetric isomorphism. Now by proposition \ref{restrict action} we may assume (notation as in the proof  of theorem \ref{supermod} (i) `$\Rightarrow$') that $\Gamma|_{\widetilde \M_g^\eps}: \widetilde \M_g^\eps \rightarrow \widetilde \M_g^0$ is supersymmetric. Therefore also $\widehat \Gamma \circ \Phi: \N \rightarrow~(\gamma \circ \tau)^*~\widetilde \M_g^0$ is supersymmetric, which implies $f = \pm 1$~.

`$\Leftarrow$': trivial by proposition \ref{restrict action}. $\Box$

\section{The remaining types} \label{remain}

In this final section we study the remaining types $(g, D)$ {\bf not} treated in section \ref{superfam}.

\begin{prop} \label{completenessgeneral} Let $\M \twoheadrightarrow V^{|q}$~, $V$ a complex manifold, be a super family of complex supermanifolds, $W$ a complex manifold, and $\varphi: W \rightarrow V$ holomorphic such that $\rund{\pi_{\varphi^* \M}^\#}_{(l)} sT^\rel \varphi^* \M = 0$ for all $l \geq 2$~. Then equivalent are:

\item[(i)] $\varphi^* (d \M): \varphi^* sT V^{|q} \rightarrow \rund{\pi_{\varphi^* \M}^\#}_{(1)} sT^\rel \varphi^* \M$ is a projection,

\item[(ii)] for every super family $\N \twoheadrightarrow \W$ with $\W^\# = W$~, $\W$ of nilpotency index $\geq 2$~, and fiberwise biholomorphic super family morphism $(\Sigma, \sigma): \N|_{\W^\natural} \rightarrow \M$ such that $\sigma^\# = \varphi$ there exists, locally in $\W$~, a fiberwise biholomorphic super family morphism $\Xi: \N \rightarrow \M$ such that $\Xi|_{\N|_{\W^\natural}} = \Sigma$~.
\end{prop}

{\it Proof:} (i) $\Rightarrow$ (ii): same as the local construction in the proof of proposition \ref{versalitygeneral} (i) $\Rightarrow$ (ii). Now use the surjectivity of $\lambda_{\m^{k - 1}, 1}$ from lemma \ref{coherent iso} (ii).

(ii) $\Rightarrow$ (i): similar to the surjectivity in the proof of proposition \ref{versalitygeneral} (ii) $\Rightarrow$ (i). Let $\beta~\in~H^0 \rund{W, \rund{\pi_{\varphi^* \M}^\#}_{(1)} sT^\rel \varphi^* \M}_1$~. After shrinking $W$ and taking a suitable atlas of $\varphi^* \M$ with local super charts $\Omega_i^{|n}$~, $\Omega_i \subset \cz^m$ open, and transition morphisms $\Psi_{i j}: \Omega_{i j}^{|n} \mathop{\rightarrow}\limits^\sim \Omega_{j i}^{|n}$ we may assume that $\beta$ is represented by the $1$-cocycle $\rund{\eps_{i j}}$~, each $\eps_{i j}$ an odd section of $sT^\rel \varphi^* \M$ on the overlap of the two local super charts $i$ and $j$ and expressed in the local super chart $i$~. Therefore

\[
\rund{\Id_{\cz^{0|1}}, \Psi_{i j}} \circ \rund{\Id_{\cz^{0|1} \times \Omega_{i j}^{|n}} + \eta \eps_{i j}} : \cz^{0|1} \times \Omega_{i j}^{|n} \mathop{\rightarrow}\limits^\sim \cz^{0|1} \times \Omega_{j i}^{|n}
\]

are the transition morphisms of a super family $\N \twoheadrightarrow W^{|1}$ with $\N|_W = \varphi^* \M$ and $(d \N)|_W \partial_\eta = \beta$~. By (ii) there exists, locally in $W$~, a fiberwise biholomorphic super family morphism $(\Xi, \xi): \N \rightarrow \M$ such that $\xi^\# = \varphi$ and $\Xi|_{\N|_W} = \breve \varphi_\M$~. Therefore $\xi = \varphi + \eta \delta$ with some $\delta \in \varphi^* \rund{sT V^{|q}}_1$ and so $\rund{\varphi^* (d \M)} \delta = (d \N)|_W \partial_\eta = \beta$~. Same calculation for $\beta \in H^0 \rund{W, \rund{\pi_{\varphi^* \M}^\#}_{(1)} sT^\rel \varphi^* \M}_0$ regarding $\bigwedge \cz$ as purely even. $\Box$ \\

{\bf From now on let $g \geq 1$ and $D \in \schweif{- 4 g + 4, \dots, - 2 g + 2} \cup \schweif{0, \dots, 2 g - 2}$~.} We can use proposition \ref{completenessgeneral} to construct locally complete families of compact super Riemann surfaces of type $(g, D)$ : \\

If $g \geq 2$ and $D \in \schweif{- 4 g + 4, \dots, - 2 g + 2}$ then $0 \leq \deg \rund{T^\rel M_{g, D} \otimes L_{g, D}^*} \leq 2 g - 2$~, and $s_D := h^1\rund{X_v, \left.L_{g, D}\right|_{X_v}} = g - 1 - D$ is independent of $v \in V_{g, D}$~. In this case define $r_D := g$~. \\

If $g \geq 2$ and $D \in \schweif{0, \dots, 2 g - 2}$ then $r_D := h^1\rund{T X_v \otimes \left.L_{g, D}^*\right|_{X_v}} = 3(g - 1) + D$ is independent of $v \in V_{g, D}$~, and now define $s_D := g$~. \\

Finally if $g = 1$ and $D = 0$ define $r_D := s_D := 1$~.

\begin{theorem}[Existence of locally complete super families] \label{superfamcomplete} For every $v \in V_{g, D}$ there exists an open neighbourhood $V \subset V_{g, D}$ of $v$ and a super family $\widetilde \M \twoheadrightarrow V^{\left|r_D + s_D\right.}$ of compact super Riemann surfaces such that $\left.\widetilde \M\right|_V = \left.\M_{g, D}\right|_V$ and

\[
\varphi^* \rund{d \widetilde \M}: \varphi^* sT V^{\left|r_D + s_D\right.} \rightarrow \rund{\pi_{\varphi^* M_{g, D}}}_{(1)} sT^\rel \varphi^* \M
\]

is a projection for every complex manifold $W$ and holomorphic map $\varphi: W \rightarrow V$~.

{\bf Local completeness of $\widetilde \M$~:} For every super family $\N \twoheadrightarrow \W$ of compact super Riemann surfaces and every isomorphism $\sigma: \Y_{w_0} := \pi_\N^{- 1}\rund{w_0} \mathop{\rightarrow}\limits^\sim \X_{v_0}$~, $v_0 \in V$ and $w_0 \in \W^\#$~, there exists, locally in $\W$ around $w_0$~, a fiberwise biholomorphic super family morphism \\
$(\Xi, \xi): \N \rightarrow \widetilde \M$ such that $\xi\rund{w_0} = v_0$ and $\Xi|_{\Y_{w_0}} = \sigma$~.
\end{theorem}

For the proof we need

\begin{lemma} \label{projection} Let $\pi: M \twoheadrightarrow V$ be a holomorphic family of compact Riemann surfaces of genus $g$~, $V$ a complex manifold, and $L \twoheadrightarrow M$ a holomorphic line bundle of degree $\geq 0$~. Then $\pi_{(2)} L = 0$~, and locally in $V$ there exists\- a projection $\rho: \O_V^{\oplus g} \twoheadrightarrow \pi_{(1)} L$ of $\O_V$-modules such that for any holomorphic map $\varphi: W \rightarrow V$~, $W$ another complex manifold, the pullback $\varphi^* \rho: \O_W^{\oplus g} \rightarrow \rund{\pi_{\varphi^* M}}_{(1)} {\breve \varphi_M}^* L$ of $\rho$ under $\varphi$ and $\breve \varphi_M$ is again a projection.
\end{lemma}



Infact, since $\pi$ is proper, $\pi_{(1)} L$ is coherent by the direct image theorem 10.4.6 of \cite{GrauRem}. \\

{\it Proof:} After shrinking $V$ let $n \in \Gamma^\hol(V, M)$ and $E \twoheadrightarrow M$ be the holomorphic line bundle associated to the divisor $- (\deg L + 1) n(V)$~. $\deg (L \otimes E) = - 1$~, so \cite{GrauRem} theorem 10.5.5 tells us that locally $\pi_{(1)} (L \otimes E) \simeq \O_V^{\oplus g}$ and $\pi_{(2)} (L \otimes E) = 0$~. Furthermore, the canonical embedding $\iota: L \otimes E \hookrightarrow L$ as holomorphic sections of $L$ having a zero at $n(V)$ of order at least $\deg L + 1$ induces an exact sequence

\[
\pi_{(1)} (L \otimes E) \mathop{\longrightarrow}\limits^\rho \pi_{(1)} L \rightarrow \pi_{(1)} (L / \Im \iota) \rightarrow \pi_{(2)} (L \otimes E) \rightarrow \pi_{(2)} L \rightarrow \pi_{(2)} (L / \Im \iota) \,~,
\]

where $\pi_{(1)} (L / \Im \iota) = \pi_{(2)} (L / \Im \iota) = 0$ by Cartan's theorem B since $L / \Im \iota$ is coherent with support $n(V)$~. Therefore also $\pi_{(2)} L = 0$~, and $\rho$ is sujective.

The rest follows easily from the obvious equivariance of the construction of $\rho$ under $\varphi$ and $\breve \varphi_M$~. $\Box$ \\

{\it Proof of theorem \ref{superfamcomplete}:} similar to the proof of theorem \ref{constructsuperfam}. For $n \in \nz \setminus \schweif{0}$ define the ringed spaces

\[
\V_n := \rund{V, \O_V \otimes \bigwedge \cz^{r_D + s_D} \left/ \rund{\bigwedge\nolimits^n \cz^{r_D + s_D}}\right.} = \rund{\schweif{0}, \bigwedge \cz^{r_D + s_D} \left/ \rund{\bigwedge\nolimits^n \cz^{r_D + s_D}}\right.} \times V \,~.
\]

Under suitable shrinking of $V$ as a neighbourhood of $v$ we construct super families $\M_n \twoheadrightarrow \V_n$ such that $\left.\M_n\right|_{\V_{n - 1}} = \M_{n - 1}$ by induction on $n$~. Finally $\widetilde \M := \M_{r_D + s_D + 1}$ will have the desired properties. \\

We start with $V := V_{g, D}$ and $\M_1 := \M_{g, D}$~. Choose an atlas of $\M_1$ with local super charts $U_i^{|1}$~, $U_i \subset V_{g, D} \times \cz$ open, and transition morphisms $\Phi_{i j}: U_{i j}^{|1} \mathop{\rightarrow}\limits^\sim U_{j i}^{|1}$~. \\

If $g \geq 2$ and $D \in \schweif{- 4 g + 4, \dots, - 2 g + 2}$ then after shrinking $V$ let $\beta^\rho := \pi e_\rho$~, $\rho = 1, \dots, g$~, where the projection $\pi: \O_V^{\oplus g} \twoheadrightarrow \rund{\pi_{M_{g, D}} }_{(1)} L$ is given by lemma \ref{projection} with \\
$L := T^\rel M_{g, D}~\otimes~L_{g, D}^*$~. $\rund{\pi_{M_{g, D}} }_{(1)} L_{g, D}$ is globally free of rank $s_D$ by \cite{GrauRem} theorem 10.5.5 and Grauert's theorem, so let $\rund{\gamma^1, \dots, \gamma^{s_D}}$ be a global frame.

If $g \geq 2$ and $D \in \schweif{0, \dots, 2 g - 2}$ then $\rund{\pi_{M_{g, D}} }_{(1)} \rund{T^\rel M_{g, D} \otimes L_{g, D}^*}$ is globally free of rank $r_D$ again by \cite{GrauRem} theorem 10.5.5 and Grauert's theorem, so let $\rund{\beta^1, \dots, \beta^{r_D}}$ be a global frame. After shrinking $V$ let $\gamma^\sigma := \pi e_\sigma$~, $\sigma = 1, \dots, g$~, with the projection $\pi: \O_V^{\oplus g} \twoheadrightarrow \rund{\pi_{M_{g, D}} }_{(1)} L$ given by lemma \ref{projection} with $L := L_{g, D}$~.

Finally if $g = 1$ and $D = 0$ let $\beta := \pi 1$ with the projection $\pi: \O_V \twoheadrightarrow \rund{\pi_{M_{g, D}} }_{(1)} L$ given by lemma \ref{projection} with $L := T^\rel M_{1, 0} \otimes L_1^* \simeq L_1^*$ and $\gamma := \pi 1$ with the projection $\pi: \O_V \twoheadrightarrow \rund{\pi_{M_{g, D}} }_{(1)} L$ given by lemma \ref{projection} with $L := L_1$~. \\

After refining the atlas of $\M_{g, D}$ and again shrinking $V$~, $\beta^\rho$~, $\rho = 1, \dots, r_D$~, and $\gamma^\sigma$~, $\sigma = 1, \dots, s_D$~, are given by cocycles $\rund{b_{i j}^\rho \zeta \partial_z}$ resp. $\rund{c_{i j}^\sigma \partial_\zeta}$~, all $b_{i j}^\rho, c_{i j}^\sigma \in \O\rund{U_{i j}}$~. \\

We define $\M_2$ by the local super charts $\left.\rund{\V_1 \times \cz^{|1}}\right|_{U_i}$ with transition morphisms

\begin{eqnarray*}
&& \Phi_{i j}^2 := \rund{\Id_{\rund{\schweif{0}, \bigwedge \cz^{r_D + s_D} \left/ \rund{\bigwedge\nolimits^1 \cz^{r_D + s_D}}\right.}}, \Phi_{i j}} \circ \rund{\Id_{\left.\rund{\V_2 \times \cz^{|1}}\right|_{U_{i j}} } + \eta_\rho b_{i j}^\rho \zeta \partial_z + \vartheta_\sigma c_{i j}^\sigma \partial_\zeta} : \\
&& \phantom{12} \left.\rund{\V_2 \times \cz^{|1}}\right|_{U_{i j}} \mathop{\rightarrow}\limits^\sim \left.\rund{\V_2 \times \cz^{|1}}\right|_{U_{j i}} \,~.
\end{eqnarray*}

$n \rightarrow n + 1$ : As in the induction step of the proof of theorem \ref{constructsuperfam} we see that the obstructions to construct - locally in $\V_{n + 1}$ - a super family $\M_{n + 1} \twoheadrightarrow \V_{n + 1}$ with $\left.\M_{n + 1}\right|_{\V_n} = \M_n$ lie in

\[
\rund{\pi_{M_{g, D}} }_{(2)} \rund{\rund{\bigwedge\nolimits^n \cz^{r_D + s_D}} \boxtimes sT^\rel \M_{g, D}}_0 \,~,
\]

which is zero by lemma \ref{projection} and \cite{GrauRem} theorem 10.5.5.

{\it Local completeness of $\widetilde \M$} now follows from the anchoring property of $\M_{g, D}$ and induction over the nilpotency index of $\W$ using proposition \ref{completenessgeneral}. $\Box$ \\





However, because of the presence of special divisors there is {\bf no} chance to construct a versal super family of type $(g, D)$~, even if one allows arbitrary admissible bases:

\begin{theorem}[Non-existence of versal super families] \label{noversalfam} There exists {\bf no} super family $\M \twoheadrightarrow \V$ of type $(g, D)$ being

{\bf complete:} for every super family $\N \twoheadrightarrow \W$ of compact super Riemann surfaces of type $(g, D)$ there exists locally in $\W$ a fiberwise biholomorphic super family morphism $\N \rightarrow \W$~.


and

{\bf infinitesimal universal:} $(d \M)_v: sT_v \V \rightarrow H^1\rund{sT \X'_v}$~, $\X'_v := \pi_\M^{- 1}(v)$~, is an isomorphism for all $v \in \V^\#$~. \\



\end{theorem}

{\it Proof:} Assume to the contrary that $\M$ is a family with these properties. By standard Riemann surface theory we know that $h^1\rund{sT \X_v}_1$ is {\bf not} locally independent of $v \in V_{g, D}$~. On the other hand by \cite{GrauRem} theorem 10.5.4 all sets $\schweif{v \in V_{g, D} \ \left| \ h^1\rund{sT \X_v}_1 \geq d\right.}$~, $d \in \nz$~, are analytic, and so there exists an embedding $\varphi: U \hookrightarrow V_{g, D}$~, $U \subset \cz$ an open neighbourhood of $0$~, such that $h^1\rund{sT \X_{\varphi(u)} }_1$ is independent of $u$ on $U \setminus \{0\}$ but smaller than $h^1\rund{sT \X_{\varphi(0)}}_1$~. Therefore by lemma \ref{projection} and since $\rund{\pi_{\varphi^* M_{g, D}}}_{(1)} {\breve \varphi_{\M_{g, D}} }^* \rund{sT^\rel \M_{g, D}}_1$ is coherent, we see that after shrinking $U$ there exists \\


\[
\beta \in H^0\rund{U, \rund{\pi_{\varphi^* M_{g, D}} }_{(1)} {\breve \varphi_{\M_{g, D}} }^* sT^\rel \M_{g, D}}_1
\]


such that $\beta(0) \not= 0$ in $H^1\rund{sT \X_{\varphi(0)}}$ and $\beta|_{U \setminus 0} = 0$~. After again shrinking $U$ we construct, as in the proof of lemma \ref{completenessgeneral} (ii) $\Rightarrow$ (i), a super family $\N \twoheadrightarrow U^{|1}$ of compact super Riemann surfaces having $\N|_U = \varphi^* \M_{g, D}$ with transition morphisms

\[
\rund{\Id_{\cz^{0|1}}, \Psi_{i j}} \circ \rund{\Id_{\cz^{0|1} \times \Omega_{i j}^{|1}} + \eta \eps_{i j}} : \cz^{0|1} \times \Omega_{i j}^{|1} \mathop{\rightarrow}\limits^{\sim} \cz^{0|1} \times \Omega_{j i}^{|1} \,~,
\]

$\rund{\eps_{i j}}$ a $1$-cocycle in ${\breve \varphi_{\M_{g, D}} }^* \rund{sT^\rel \M_{g, D}}_1$ representing $\beta$~, each $\eps_{i j}$ expressed in the pullback of the local super chart $i$ of $\M_{g, D}$ under $\breve \varphi_{\M_{g, D}}$~.

Now by completeness of $\M$ there exists, after shrinking $U$~, a fiberwise biholomorphic super family morphism $(\Xi, \xi): \N \rightarrow \M$~. Let $u \in U \setminus \schweif{0}$ be arbitrary. Then $(d \M)_{\varphi(u)} (d \xi)_u \partial_\eta = (d \N)_u \partial_\eta = \beta(u) = 0$~, and so $(d \xi)_u \partial_\eta = 0$ by the infinitesimal universality of $\M$~. We see that $\xi|_{\rund{U \setminus \schweif{0}}^{|1}}$ is independent of $\eta$~, and so pullback of the atlas of $\M$ gives an atlas of $\N$~, whose transition morphisms, when restricted to $\N|_{(U \setminus \schweif{0})^{|1}}$~, are independent of $\eta$~. But this implies that the transition morphisms of this atlas must be completely independent of $\eta$~, and so $\beta(0) = (d \N)_0 \partial_\eta = 0$~. Contradiction! $\Box$ \\

Now, in the end, allthough the construction of a versal super family for type $(1, 0)$ is {\bf not} possible by theorem \ref{noversalfam}, we construct a versal supersymmetric super family of compact super Riemann surfaces of genus~$1$~. We may assume $\E_1 = \schweif{0, \eps_0}$ with $\eps_0 \in \frac{1}{2} \zz^2 \setminus \schweif{0}$~. Then $h^1\rund{sT^+ \X_v}_1 = 1$ for all $v \in \T_1^0$~, and $H^1\rund{sT \X_v}_1 = 0$ for all $v \in \T_1^{\eps_0}$~. Therefore similar to the proof of theorem \ref{constructsuperfam}, $n = 2$~, one can construct a supersymmetric super family $\widetilde \M_1^0 \twoheadrightarrow s\T_1^0 := \rund{\T_1^0}^{|1}$ such that $\left.\widetilde \M_1^0\right|_{\T_1^0} = \M_1^0$ and

\[
\rund{d \widetilde \M_1^0}_v : sT_v s\T_1^0 \rightarrow H^1\rund{sT^+ \X_v}
\]

is an isomorphism for all $v \in \T_1^0$~. Furthermore we can write $\widetilde \M_1^0$ explicitly as

\[
\Pr_{s\T_1^0}: \left.\rund{s\T_1^0 \times \cz^{1|1}} \right/ \spitz{\widetilde S, \widetilde T} \twoheadrightarrow s\T_1^0
\]

with commuting

\[
\widetilde S := \rund{\left.\Id_{s\T_1^0}, z + 1 \, \right| \, \zeta} \text{ and } \widetilde T := \rund{\left.\Id_{s\T_1^0}, z + u + \zeta \chi \, \right| \, \zeta + \chi}~\in~\Aut_\strong \rund{s\T_1^0 \times \cz}
\]

and standard supersymmetry. Hereby $u$ and $\chi$ denote the even resp. odd super coordinate on the base $s\T_1^0$~. Similar to proposition \ref{versalitygeneral} (i) $\Rightarrow$ (ii) and theorem \ref{versalitysuper} one proves

\begin{theorem}[Versality of $\widetilde \M_1^0 \dot\cup \M_1^{\eps_0}$ ] \label{versalitysuperSUSY1}
\item[(i)] {\bf Completeness:} Let $\N \twoheadrightarrow \W$ be a su\-per\-sym\-me\-tric super family of compact super Riemann surfaces of genus $1$~. If $\W^\#$ is a contractible Stein manifold then there exists a supersymmetric fiberwise biholomorphic super family morphism $(\Xi, \xi): \N \rightarrow \widetilde \M_1^0 \dot\cup \M_1^{\eps_0}$~. More precisely we have again

{\bf Local completeness:} For every supersymmetric isomorphism $\sigma: \Y_{w_0} := \pi_\N^{- 1}\rund{w_0} \mathop{\rightarrow}\limits^\sim \X_{v_0}$~, $w_0 \in \W^\#$ and $v_0 \in \T_1^0$ resp. $\T_1^{\eps_0}$~, there exists a supersymmetric fiberwise biholomorphic family morphism $(\Xi, \xi): \W \rightarrow \widetilde \M_1^0$ resp. $\M_1^{\eps_0}$ such that $\xi\rund{w_0} = v_0$ and $\Xi|_{\Y_{w_0}} = \sigma$~.

and the

{\bf Anchoring property:} For every fiberwise biholomorphic super family morphism \\
$(\Sigma, \sigma): \N|_{\W^\natural} \rightarrow \widetilde \M_1^0$ resp. $\M_1^{\eps_0}$~, $\W$ of nilpontency index $\geq 2$ and $\W^\#$ a Stein manifold, there exists a unique morphism $\xi: \W \rightarrow s\T_1^0$ resp. $\T_1^{\eps_0}$ such that $\xi|_{\W^\natural} = \sigma$ and $\xi$ can be extended to a supersymmetric fiberwise biholomorphic super family morphism $(\Xi, \xi): \N \rightarrow \widetilde \M_1^0$ resp. $\M_1^{\eps_0}$ such that $\Xi|_{\N|_{\W^\natural}} = \Sigma$~.
\item[(ii)] {\bf Infinitesimal universality:} $d \rund{\widetilde \M_1^0 \dot\cup \M_1^{\eps_0}}_v: sT_v \rund{s\T_1^0 \dot\cup \T_1^{\eps_0}} \rightarrow H^1\rund{sT^+ \X_v}$ is an isomorphism for all $v \in \T_1^0 \cup \T_1^{\eps_0}$~.
\end{theorem}

\begin{cor} \label{SUSYunique1}
Corollary \ref{SUSYunique} is true also for supersymmetric super families of compact super Riemann surfaces of genus $1$~.
\end{cor}

{\it Proof:} Let $\N \twoheadrightarrow \W$ be a super family of compact super Riemann surfaces of genus $1$ with two supersymmetries $\D$ and $\D'$~. By theorem \ref{versalitysuperSUSY1} we may assume without loss of generality that $\N = \xi^* \rund{\widetilde \M_1^0 \dot\cup \M_1^{\eps_0}}$ with some $\xi:~\W~\rightarrow s\T_1^0 \dot\cup \T_1^{\eps_0}$ and $\D$ is the pullback of the supersymmetry on $\widetilde \M_1^0 \dot\cup \M_1^{\eps_0}$ by $\breve \xi_{\widetilde \M_1^0 \dot\cup \M_1^{\eps_0}}$~. If $\xi$ maps to $\T_1^{\eps_0}$ then the assertion is true by corollary \ref{SUSYunique0}. So we may assume that $\xi$ maps to $s\T_1^0$~. In this case write $\xi = (h | \kappa)$ with $h \in \O(\W)_0$ such that $h^\#$ maps to $\T_1^0$ and $\kappa \in \O(\W)_1$~. Then

\[
\N = \left.\rund{\W \times \cz^{1|1}} \right/ \spitz{\sigma, \tau}
\]

with $\sigma := \rund{\left.\Id_\W~, z + 1 \, \right| \, \zeta}$ and $\tau := \rund{\left.\Id_\W~, z + h + \zeta \kappa \, \right| \, \zeta + \kappa}$~. Let $\widetilde \D$ and $\widetilde \D'$ denote the pullbacks of $\D$ resp. $\D'$ under the canonical strong projection $\W \times \cz^{1|1} \twoheadrightarrow \N$~. Then $\widetilde \D$ is the standard supersymmetry and $\widetilde \D'$ generated by $(\zeta + \lambda) \partial_z + F \partial_\zeta$ with suitable $\lambda \in \O(\W \times \cz)_1$ and $F \in \O(\W \times \cz)_0^\times$~. $\sigma$ and $\tau$ are supersymmetric also with respect to $\widetilde \D'$~, so $\lambda$ and $F$ are $1$-periodic, and

\begin{eqnarray*}
\lambda &=& \lambda(z - h) + \rund{\rund{\lambda' \lambda + F}(z - h) - 1} \kappa \,~, \\
F &=& F(z - h) + \rund{F \lambda}'(z - h) \kappa \,~,
\end{eqnarray*}

where $\lambda'$ and $\rund{F \lambda}'$ denote the derivatives w.r.t. the $z$-coordinate of $\cz$~. Therefore $\lambda$ and $F$ are infact independent of $z$~, and $(F - 1) \kappa = 0$~. So locally in $\W$ we can take $\sqrt{F} \in \O(\W)_0^\times$ such that $\sqrt{F}^\# + 1$ has no zeros and so $\rund{\sqrt{F} - 1} \kappa = 0$~. We obtain a strong supersymmetric isomorphism

\[
\rund{\left.\Id_\W, z + \frac{\zeta \lambda}{\sqrt{F} + 1} \, \right| \, \sqrt{F} \zeta - \frac{\lambda}{\sqrt{F} + 1} }: \rund{\W \times \cz^{1|1}, \widetilde \D} \mathop{\rightarrow}\limits^\sim \rund{\W \times \cz^{1|1}, \widetilde \D'}
\]

with body $\id_{W \times \cz}$ commuting with $\sigma$ and $\tau$ and so inducing a strong supersymmetric isomorphism $\Phi: (\N, \D) \mathop{\rightarrow}\limits^\sim \rund{\N, \D'}$ with $\Phi^\# = \Id_{\N^\#}$~. $\Box$

\end{document}